%% file: LERW160305.tex
\documentclass[11pt]{article}
\usepackage{amssymb,latexsym}
\usepackage[dvips]{graphicx} 
\usepackage{url}
\usepackage{picture}
\usepackage{color}

\newtheorem{prp}{Proposition}
\newtheorem{lem}[prp]{Lemma}\newtheorem{thm}[prp]{Theorem}

\newtheorem{cor}[prp]{Corollary}
\newenvironment{prf}{\begin{trivlist}\item[\emph{Proof.}]}{\end{trivlist}
  \medskip\par}
 
\newenvironment{rem}{\begin{trivlist}\item[\emph{Remarks.}]}{\end{trivlist}  \medskip\par}
\def\prpb{\begin{prp}}\def\prpe{\end{prp}}
\def\lemb{\begin{lem}}\def\leme{\end{lem}}
\def\thmb{\begin{thm}}\def\thme{\end{thm}}
\def\corb{\begin{cor}}\def\core{\end{cor}}
\def\prfb{\begin{prf}}\def\prfe{\end{prf}}
\def\remb{\begin{rem}}\def\reme{\end{rem}}
\def\prpa#1{\label{p:#1}}\def\prpu#1{Proposition~\ref{p:#1}}

\def\thma#1{\label{t:#1}}\def\thmu#1{Theorem~\ref{t:#1}}
\def\cora#1{\label{c:#1}}\def\coru#1{Corollary~\ref{c:#1}}
\def\seca#1{\label{s:#1}}\def\secu#1{Section~\ref{s:#1}}

\def\itmb{\begin{enumerate}}\def\itme{\end{enumerate}}
\def\itdb{\begin{itemize}}\def\itde{\end{itemize}}
\def\ittb{\begin{description}}\def\itte{\end{description}}
\def\eqnb{\begin{equation}}\def\eqne{\end{equation}}
\def\arrb#1{\begin{array}{#1}}\def\arre{\end{array}}
\def\tabb#1{\par\noindent\begin{tabular}{#1}}
\def\tabe{\end{tabular}\par\noindent}
\def\eqna#1{\label{e:#1}}\def\eqnu#1{(\ref{e:#1})}

\def\QED{\relax\ifmmode\let\@tempa\relax\ifcase\@eqcnt\def\@tempa{& & &}\or
  \def\@tempa{& &}\else\def\@tempa{&}\fi\@tempa $\Box$ \else\hfill $\Box$ \fi}
\def\DDD{\relax\ifmmode\let\@tempa\relax\ifcase\@eqcnt\def\@tempa{& & &}\or
 \def\@tempa{& &}\else\def\@tempa{&}\fi\@tempa $\Diamond$
 \else\hfill $\Diamond$ \fi}

\def\Rom#1{\uppercase\expandafter{\romannumeral#1}}
\def\dsp{\displaystyle}

\def\liminf{\displaystyle \mathop{\underline{\lim}}\limits}
\def\limsup{\displaystyle \mathop{\overline{\lim}}\limits}

\def\Ccomb#1#2{\setbox0=\hbox{$\displaystyle\mathrm{C}$}\setbox1=\hbox{%
$\scriptstyle #1$}\kern \wd1{\mathrm{C}}_{\kern -1.05\wd0\kern -0.99\wd1{#1}
 \kern 1.15\wd0{#2}}}

\def\vec#1{{\mathbf{#1}}} \def\vec#1{\overrightarrow{#1}}
\def\clvec#1#2#3{\def\clvecone{#3}\left(\arrb{c} \dsp #1\\ \dsp #2
 \ifx\clvecone\empty\else\\ \dsp #3\fi\arre\right)}

\def\le{\leqq} \def\leq{\leqq} \def\geq{\geqq}

\def\reals{{\mathbb R}}\def\twreals{{\mathbb R^2}}

\def\integers{{\mathbb Z}}\def\pintegers{{\mathbb Z}_+}
\def\nintegers{{\mathbb N}}

\def\prb#1{\def\prbone{#1}
  \ifx\prbone\empty{\mathrm{P}}\else{\mathrm{P[\;}}#1{\mathrm{\;]}}\fi}
\def\prbseq#1#2{\def\prbseqone{#2}
  \ifx\prbseqone\empty{\mathrm{P}}_{#1}\ignorespaces
  \else{\mathrm{P}}_{#1}{\mathrm{[\;}}#2{\mathrm{\;]}}\fi}

\def\EEseq#1#2{\def\EEseqone{#2}
  \ifx\EEseqone\empty{\mathrm{E}}_{#1}\else
 {\mathrm{E}}_{#1}{\dsp\mathrm{[\;}}#2{\mathrm{\;]}}\fi}

\def\VVseq#1#2{\def\VVseqone{#2}
  \ifx\VVseqone\empty{\matrm{V}}_{#1}\else
 {\mathrm{V}}_{#1}{\dsp\mathrm{[\;}}#2{\mathrm{\;]}}\fi}

\def\sg{Sierpi\'{n}ski gasket}
\def\parr{\par\noindent}
\def\eqsb{\begin{eqnarray*}}\def\eqse{\end{eqnarray*}}

\makeatletter
\@addtoreset{equation}{section}
\makeatother

\title{
 A family of self-avoiding random walks interpolating the loop-erased random walk 
and a self-avoiding walk on the Sierpi\'nski gasket
}
\author{
Kumiko Hattori
\\ \and
Noriaki Ogo 
\\ \and
Takafumi Otsuka 
}

\date{\today}
\begin{document}
\maketitle

\footnotetext{Department of Mathematics and Information Sciences,
 Tokyo Metropolitan University, Hachioji, Tokyo 192-0397, Japan.}








\begin{center}
ABSTRACT
\end{center}
We show that the  `erasing-larger-loops-first' (ELLF) method, which was first 
introduced for erasing loops from the simple random walk on the \sg 
, does work  also for non-Markov random walks, in particular,  
self-repelling walks  to construct a new family of self-avoiding walks 
on the \sg . 
The one-parameter family constructed in this method continuously 
connects the loop-erased random walk and a 
self-avoiding walk 
which has the same asymptotic behavior as the `standard' self-avoiding walk.
We prove the existence of the scaling limit and study some path properties: 
The exponent $\nu$ governing the short-time behavior of the 
scaling limit varies continuously in the parameter. The limit process 
is almost surely self-avoiding, while it has  
path Hausdorff dimension $1/\nu $, which is strictly greater than 1.

\parr

\vspace*{1in}\par

\bigskip\par
\noindent\textit{MSC2010 Subject Classifications:}
60F99, 60G17, 28A80, 37F25, 37F35 
\bigskip\par
\noindent\textit{Corresponding author:} 
Kumiko Hattori, 
khattori@tmu.ac.jp
\parr
Department of Mathematics and Information Sciences,
 Tokyo Metropolitan University, Hachioji, Tokyo 192-0397, Japan.
\parr tel: +81 42 677 2475


\section{Introduction}
\seca{Intro}

The self-avoiding walk (SAW) and the loop-erased random walk (LERW)
are two typical examples of non-Markov random walks on graphs.
The self-avoiding walk is defined by the uniform measure on  
self-avoiding paths of a given length.
In this paper we call this model the `standard' self-avoiding walk
(`standard' SAW),  
for we shall deal with a family of different walks whose paths are self-avoiding. 
The loop-erased random walk is a 
random walk obtained by erasing loops from the simple random 
walk in chronological order (as soon as each loop is made).  
Although the LERW has self-avoiding paths, 
it has a different distribution from that of the `standard' SAW.

Two of the basic questions concerning random walks are:
\par\noindent
(1) What is the asymptotic behavior of the walk as the number of steps
tends to infinity?
To be more specific,  
if $X(n)$ denotes the location of the walker starting at the origin after $n$ steps,
does the mean square displacement show a power behavior?  In other words, does 
the following hold in some sense? 
\[E[|X(n)|^2] \sim n^{2\nu },\]
where $|X(n)|$ denotes the Euclidean distance from the starting point and 
$\nu$ is a positive constant.  If it is the case, what is the value of 
the displacement exponent $\nu$?

\par\noindent
(2) Does the walk have a scaling limit?
A scaling limit is the limit as the edge length of the graph tends to $0$.  
To give some examples, Brownian motion on ${\mathbb Z}^d$ and Brownian motion on the Sierpi\'nski gasket are obtained as the scaling limit of the simple random walk on their  
respective graph approximations.   
The displacement exponent $\nu$ governs also the short-time behavior of 
the scaling limit.

Question (1) originated from the problem of the end-to-end distance of 
long polymers.  Since no two monomers can occupy the same place, 
a self-avoiding walk is expected to model polymers.
There have been many works, not only mathematical works, but also 
computer simulations and heuristics aimed at answering  the question, however,  for 
`standard' self-avoiding walk on 
${\mathbb Z}^d$ with $d=2, 3, 4$, it is not solved rigorously yet.
Question (2) 
for ${\mathbb Z}^d$,  $d=2, 3, 4$ has not been given a rigorous answer yet, either,   
while for ${\mathbb Z}^d$ with $d>4$ the answers are given; the scaling limit is the 
$d$-dimensional Browinan motion and $\nu=1/2$. 
The difficulties for  $d=2, 3, 4$ lie in the strong self-avoiding effect 
in low dimensions.
For what is known about `standard' self-avoiding walks on  ${\mathbb Z}^d$, see \cite{MS}.

The situation is quite different for the LERW on ${\mathbb Z}^d$.  
The existence of the scaling limit has been proved for all $d$, 
and the asymptotic behavior has been studied in terms of 
the growth exponent (the reciprocal of the displacement 
exponent).  For $d=2$ Schramm-Loewner 
evolution (SLE) has played an essential role.  
For some further discussion of the LERW on ${\mathbb Z}^d$, 
see \cite{LSW}, \cite{schramm}, \cite{Kozma} and \cite{lawler2}.

 The Sierpi\'nski gasket provides a space which is `low-dimensional', but 
permits rigorous analysis.  For this fractal space, 
the displacement exponent $\nu$ of the `standard' SAW is 
obtained in \cite{HK}.  The scaling limit is studied in \cite{hh} and it is 
proved that the same $\nu$ governs the short-time behavior of 
the limit process $X_t$, that is, there exist positive constants $C_1$ and 
$C_2$ such that 
\[C_1 \leq \frac{E[|X_t |]}{t^{\nu}}\leq C_2  \]
holds for small enough $t$ (\cite{HHH}). 
As for the LERW, the scaling limit was obtained by two groups independently, using different 
methods (\cite{STW}, \cite{HM}). 

SLE mentioned above is a profound theory, which goes far beyond 
the investigation of the scaling limit of the LERW on ${\mathbb Z}^2$.
It is a unified theory for a variety of random curves in ${\mathbb R}^2$ that 
involves a parameter $\kappa$,  and different values of $\kappa$ correspond
to different models. $\kappa =2$ corresponds to the scaling limit of the LERW and 
$\kappa =8/3$ is conjectured to be the scaling limit of the SAW.   
Thus, SLE  is expected to connect the SAW and 
the LERW on ${\mathbb R}^2$.

There arises a natural question:  Is it possible to construct a model that 
connects the SAW and the LERW on the Sierpi\'nski gasket continuously in some parameter?
In this case we cannot use SLE, for which the conformal invariance of models in  
${\mathbb R}^2$ plays  an essential role. 
  
In this paper, we construct a one-parameter family of self-avoiding random walks 
on the Sierpi\'nski gasket continuously  connecting  the LERW and a 
SAW  which has the same asymptotic behavior as the `standard' SAW.
We prove the existence of the scaling limit and show some path properties: 
The exponent $\nu$ governing the short-time behavior of the 
scaling limit varies continuously in the parameter. The limit process 
is almost surely self-avoiding, while it has  
path Hausdorff dimension $1/\nu $, which is strictly greater than 1.

The main ingredients for the model are the one-parameter family of self-repelling walks 
on the Sierpinski gasket studied in \cite{HHH} and \cite{HH},  
and  the `erasing-larger-loops-first' (ELLF) method employed 
in the study of the LERW \cite{HM}.
A self-repelling walk is a walk that is discouraged, if not prohibited, 
to return to points it has visited before.  There have been 
a variety of  models on ${\mathbb Z}$.  See, for example, the survey paper  \cite{hofkonig} and the references therein.  
The model we use here is unique in the way of discouraging returns; 
penalties are given for backtracks and sharp turns, rather than for revisits 
to the same points or the same edges.   

For the `standard' LERW on graphs, 
the uniform spanning tree proves to be a powerful tool (\cite{STW}). 
By `standard', we mean the loops are erased chronologically as first 
introduced by G. Lawler  (\cite{Lawler}).  On the other hand, 
\cite{HM} constructed a LERW on the Sierpi\'{n}ski gasket by ELLF, that is, 
by erasing loops in 
descending order of size of loops and proved that the resulting LERW has the same 
distribution as that of the `standard' LERW.  
The uniform spanning tree is powerful in the sense that
it can be used on any graph, however,  this tool is valid only 
for loop erasure from {\it simple} random walks.  
We prove that ELLF does work also for other kinds of random 
walks on some fractals, in particular, for self-repelling walks on the 
Sierpi\'nski gasket, 
for the method is based on self-similarity.
Thus, our construction is performed by erasing loops  
from the family of self-repelling walks by the ELLF method.

In \secu{SRW}, we describe the set-up and recall the family of self-repelling walks 
introduced in \cite{HHH} and \cite{HH} in a more concise manner.  In \secu{LE}, we describe the ELLF method of loop-erasing in a more organized manner than \cite{HM},  
and apply it to the self-repelling walks to obtain a new family of 
self-avoiding walks interpolating LERW and SAW.  In \secu{SL} we study 
the scaling limit.  In \secu{Path} we prove some 
properties of the limit process concerning the short-time behavior. 
In \secu{CR}, we give the conclusion and some remarks.


\section{Self-repelling walk on the pre-\sg s}
\seca{SRW}

Let us first recall the definition of the pre-\sg s, that is, graph approximations of 
the Sierpi\'nski gasket which is a fractal with Hausdorff dimension $\log 3/\log 2$.
Let $O=(0,0)$, $\dsp a=(\frac{1}{2}, \frac{\sqrt{3}}{2})$, $b=(1,0)$ and 
define 
$F'_{0}$ to be the graph that consists of the three vertices and 
the three edges of $\triangle Oab$.
Define similarity maps  $f_i : {\mathbb R}^2 \to  {\mathbb R}^2$, $i=1,2,3$ by
\[f_1(x)=\frac{1}{2}x,\ f_2(x)=\frac{1}{2}(x+a), \ f_3(x)=\frac{1}{2}(x+b),\]
and a recursive sequence of graphs 
$\{F'_{N}\}_{N=0}^{\infty }$ by
\[F'_{N+1}=f_1(F'_{N}) \cup f_2(F'_{N}) \cup f_3(F'_{N}).
\]
Let $F_N$ be the union of $F'_N$ and its reflection with respect to the $y$-axis, and 
let $G_N$ and $E_N$ be the sets of the vertices and of the edges of $F_N$, respectively.
$F_3$ is shown in Fig. 1. 

\begin{figure}[htb]
\begin{center}
\input{Fig2a.tex}
\\[1\baselineskip]
\caption{$F_3$}
\end{center}
\end{figure}

Let ${\cal T}_M$ be the set of all upward (closed and filled) triangles which are translations of 
 $2^{-M} \triangle Oab$ and whose vertices are in $G_M$; an 
element of ${\cal T}_M$ is called a {\bf $2^{-M}$-triangle}. 

For each $N\in \pintegers =\{0,1,2, \ldots \}$, 
denote the set of finite paths  on $F_N$ starting from $O$, 
not hitting any vertices in $G_0$ other than 
$O$ on the way,  and stopped 
at the first hitting time of $a$ by 
\[W_N = \{\ w=(w(0), w(1), \cdots , w(n)):
\  w(0) =O, \  w(n)=a, \ w(i) \in G_N,  \] 
\[ \{w(i), w(i+1)\} \in E_N, \ w(i) \not\in G_0\setminus \{O\}, \  0 \leq i \leq n-1,  
\ n \in {\mathbb N}\ \}.\]
For a path $w=(w(0), w(1), \cdots , w(n)) \in W_N$, denote the number of steps by $\ell (w) :=n$.

If we assign probability $(1/4)^{\ell (w)-1}$ to each $w\in W_N$, then 
this determines a probability measure on paths which is 
the same as that induced by the simple random walk on $F_N$ starting from $O$ and stopped at the first  
hitting time of $a$ conditioned that the walk does not hit any vertices in $G_0\setminus \{O\}$ 
on the way.  The factor $(1/4)^{-1}$ comes from 
this conditioning. 

We shall assign  probabilities such that they give  
random walks whose revisits to the same points are discouraged.  
First let us start with paths in $W_1$.  The idea is that we
give a penalty to $w \in W_1$ every time it makes a sharp turn or 
a backtrack at $G_1\setminus G_0$,  or revisits $O$.  We realize it 
by using  $N(w)$, the total number of  sharp turns  and  
backtracks, and $M(w)$, the total number of revisits to $O$,  
and by  
assigning probability $u^{N(w)+M(w)} x_u^{\ell (w)-1}$,  
where $u$ is a parameter taking values in $[0,1]$ and $x_u$ is a positive constant 
determined so that the sum of the probabilities over $W_1$ equals to $1$.
This is a natural way to define a self-repelling walk on $F_1$: 
If $u=1$, then  we have $x_1=1/4$ and the simple random walk given above,  and 
if $u=0$, then  the probability is supported on a set of self-avoiding paths.
On a general $W_N$, we  define the probability recursively.

To give a precise definition, we shall make some preparations.   
For a path $w\in \bigcup _{N=1}^{\infty } W_N$ and $A \subset {\mathbb R}^2$,
 we define the  hitting time of $A$ 
 by
\[T_A(w)=\inf  \{j \geq 0 :\ w(j) \in A\},\]
where we set $\inf \emptyset =\infty $. 
For 
$w\in W_N$ and $0 \leq M \leq N$, we shall define a recursive sequence 
$\{T_i^M(w)\}_{i=0}^{m}$ of
{\bf hitting times of  $G_M$} as follows:
Let 
 $T_{0}^{M}(w)=0$, and for $i\geq 1$, let
\[T_{i}^{M}(w)=\inf \{j>T_{i-1}^{M}(w) :\  w(j)\in G_{M}\setminus
\{w(T_{i-1}^{M}(w))\}\},\]
here we take  $m$ to be the smallest integer such that 
$T_{m+1}^{M}(w)=\infty $. Then 
$T_{i}^{M}(w)$ is the time (steps) taken for the path $w$ to 
hit vertices in $G_{M}$ for the $(i+1)$-th time, 
under the condition that if $w$ hits the same vertex in $G_{M}$ 
more than once in a row, we count it only `once'.

For each $M \in \integers_{+}$, we define a coarse-graining map
$\displaystyle Q_{M}: \bigcup _{N=M}^{\infty} W_N \to W_M$ by setting $(Q_{M}w)(i)=w(T_{i}^{M}(w))$ 
for $i=0,1,2,\ldots,m$, where $m$ is as above.
Note that 
\[Q_{K}\circ Q_{M}=Q_{K}, \ \ \ \mbox{  if  } \ \ K \leq M\]
holds and that 
if $w \in W_N$ and $M \le N$, then
$Q_M w \in W_M$.

For $w\in W_1$, 
 define the reversing number $ 
N(w) $ and the revisiting
number $M(w)$ by
\eqnb
\eqna{N}
N(w) = \sharp \{1\leq  i\leq \ell (w)-1
\ :\  \overrightarrow{w(i-1)w(i)} \cdot 
\overrightarrow{w(i) w(i+1)} 
< 0 ,\ w(i) \notin G_0 
\ \}, 
\eqne
\eqnb
\eqna{M}
M(w)
 = \sharp \{ 1 \leq  i  \leq \ell (w)
\ :\ w(i) =O \},
\eqne
where $\vec{a}\cdot\vec{b}$ denotes the inner product of $\vec{a}$ and
$\vec{b}$ in $\reals^2$. 

For $x >0$ and $0 \leq u \leq 1$, define 
\eqnb \eqna{phi}
\Phi (x,u)=\sum _{w \in W_1}  u^{N(w)
+M(w)} \ x^{\ell (w)}.
\eqne
For each $u$, within the radius of convergence $r_u>0$ as a power series in $x$,
we have 
the following explicit form of $\Phi$ given in \cite{HHH}:
\[\Phi(x,u)=\frac{x^2 \{1+(1+u)x-u(1-u^2)x^2
+2(1-u)^2 u^2 x^3\}}{ (1+ux)(1-2ux)-4u^2 x^2 \{1+2(1-u^2)x^2
-2u(1-u)^2x^3 \} }.
\]
.
\prpb
\prpa{fp}

(Proposition 2.3 in \cite{HHH})

\itmb
\item[(1)]
For each $u \in [0,1]$, there is a unique fixed point $x_u$
of the mapping $\Phi(\cdot ,u): (0, r_u) \rightarrow 
(0, \infty) $,
that is, $\Phi(x_u,u)=x_u,$ $x_u>0$.
As a function in $u$, $x_u$ is continuous and strictly decreasing
on $[0,1]$.
\item[(2)]
Let $\dsp \tilde{\lambda}_u=\frac{\partial \Phi}{\partial x} 
(x_u,u)$. Then $\tilde{\lambda}_u >2$ and 
$\tilde{\lambda}_u$ is continuous in $u$.
\itme
\prpe

In the two extreme cases, we know that $\dsp x_0=\frac{\sqrt{5}-1}{2}$, 
$\dsp \tilde{\lambda}_0=\frac{7-\sqrt5}2$,
and $\dsp x_1=\frac{1}{4}$, $\tilde{\lambda} _1 =5$.

To define a family of probability measures $\{P_N^u, \ u\in [0,1] \}$ on each 
$W_N$,  we consider decompositions of a path based on the self-similarity and the 
symmetries of the pre-\sg s.  
Assume $w\in W_N$ and $0\leq M< N$ and denote $\tilde{w}=Q_Mw$.
Since the pair of adjacent $2^{-M}$-triangles including $\tilde{w}(i-1)$, $\tilde{w}(i)$
and $\tilde{w}(i+1)$ is similar to $F_{N-M}$,  
there is a unique decomposition
\eqnb \eqna{decomposition1}
(\tilde{w}; w_1, \cdots , w_{\ell (\tilde{w})}), \ \tilde{w} \in W_M, \ w_i \in 
W_{N-M}, \ i=1, \cdots , \ell (\tilde{w})
\eqne
such that 
the path segment  $(w(T_{i-1}^M(w)), w(T_{i-1}^M(w)+1)), \cdots , w(T_{i}^M(w)))$ of 
$w$ is identified with $w_i \in W_{N-M}$ by appropriate similarity, 
rotation, translation and reflection so that $w(T_{i-1}^M(w))$ is identified with 
$O$ and $w(T_{i}^M(w))$ with $a$. 
We shall use this kind of identification throughout this paper.
We illustrate a simple example of the decomposition in Fig. 2.

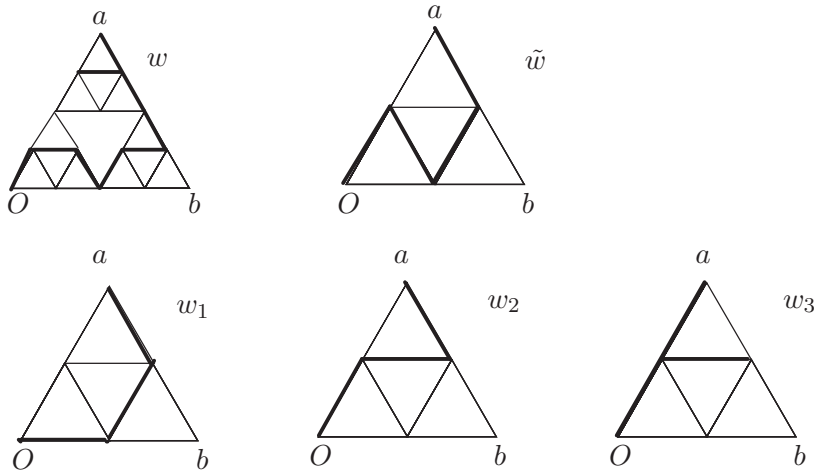
\begin{figure}[htb]
\begin{center}
\input{decompositionsmall2.tex}
\\[1\baselineskip]
\caption{$w, \tilde{w}, w_1, w_2, w_3$}
\end{center}
\end{figure}

 
First, for each $w\in W_1$,
let 
\eqnb \eqna{defprob1}
P_{1}^u(w)=u^{N(w)+M(w)}\ x_u^{\ell (w)-1},
\eqne
and define $P_N^u$ on
$W_N$ recursively by
\eqnb \eqna{defprob2}
P_N^u(w)=P_{N-1}^u(\tilde{w})\ \prod_{i=1}^{\ell (\tilde{w})}
P_1^u(w_i),\eqne
where $(\tilde{w}; w_1, \cdots , w_{\ell (\tilde{w})}) $ is the decomposition of 
$w\in W_N$ with $M=N-1$ given in \eqnu{decomposition1}. 
Denote the image measure of $P_N^u $ 
induced by the mapping $Q_M $ by $P_N^u \circ Q_M^{-1}$. 
 $P_{N}^u$ is self-similar in the 
sense that 
$P_N^u \circ Q_M^{-1} = P_M^u$.

$(W_N, \{ P_N^u\}_{u\in [0,1]})$
defines a family of self-repelling walks $Z_N^u$ on $F_N$ 
such that 
\eqnb
\eqna{Z}
Z_N^u(w)(i)=w(i),\ \  i=0, \cdots , \ell (w), \ \ w\in W_N.\eqne

In \cite{HHH}, it is proved that for each $u$, the sequence 
 $\{Z_N^u(\tilde{\lambda} _u^N\ \cdot \ )\}_{N=1}^{\infty }$ 
of time-scaled self-repelling walks 
converges to a continuous process as $N \to \infty$.  The 
one-parameter family of the limit processes $\{Z^u(\ \cdot \ ), u\in [0,1]\}$ 
continuously interpolates 
a self-avoiding process  ($u=0$) and Brownian motion ($u=1$) on the \sg.

In the next section, we erase loops from this family of  self-repelling walks  to 
obtain a one-parameter family of self-avoiding walks. 
For this purpose, 
we introduce an auxiliary family of self-repelling walks.  Let
\[F_N^V=F_ N\cup (F'_N+b),\]
which consists of three adjoining copies of $F'_N$, and 
let $G_N^V$ and $E_N^V$ be the sets of the vertices and of the edges 
of $F_N^V$, respectively.  Denote the set of finite paths on $F_N^V$ 
starting from $O$ and stopped at the first hitting time of $a$ by 
\[V_N^0=
\{\ w=(w(0), w(1), \cdots , w(n)):
\  w(0) =O, \  w(n)=a, \ w(i) \in G_N^V,  \] 
\[ \{w(i), w(i+1)\} \in E_N^V, \ w(i) \neq a, \  0 \leq i \leq n-1,  
\ n \in {\mathbb N}\ \}.\]
and define $\ell (w)$, $T_i^M(w)$ and $Q_M$ in the same way as for $W_N$.
Let 
\[V_N = \{\ w \in V_N^0 :
\  Q_0w =(O, b, a) \ \}.\]
Paths in $V_N$ are allowed to leak into the `interior' of the third 
copy of $F'_N$.
A path $w\in V_N$ defined in this way consists of two parts, $(w(0), w(1), \cdots , w(T_1^0(w)))$ and 
$(w(T_1^0(w)), w(T_1^0(w)+1),  \cdots ,$
\parr  $w(T_2^0(w)))$, and they can be identified with 
some $w', w'' \in W_N$, respectively.
Define a probability measure ${P}'^u_N$ on $V_N$ by 
\[{P}'^u_N[w]=P_N^u[w'] \cdot P_N^u[w''],\]
where $P_N^u$ is defined in \eqnu{defprob1} and \eqnu{defprob2}.  
 $(V_N, \{{P}'^u_N\}_{u\in [0,1]})$ defines a family of  self-repelling random walks  ${Z}'^u_N$ on $F_N^V$ 
such that 
\eqnb
\eqna{ZZ}
{Z}'^u_N(w)(i)=w(i), i=0, \cdots , \ell (w), \ \ w\in V_N.\eqne
This is a family of self-repelling walks that hit $b$ 'once' in the sense that 
$ Q_0w =(O,b,a)$.



\section{Loop erasure by the erasing-larger-loops-first rule}
\seca{LE}

For  $(w(0), w(1), \cdots , w(n))\in W_N\cup V_N$, 
if  there are $c \in G_N$,  $i$ and $ j$, $0\leq i < j \leq n$ such that $w(i)=w(j)=c$ and 
$w(k)\neq c$ for any $i<k<j$,  
we call the path segment $[w(i),w(i+1), \ldots , w(j)]$ a {\bf loop formed at} ${\bf c}$ and 
define its {\bf diameter} by $d=\sup _{i\leq k_1<k_2 \leq j} |w(k_1)-w(k_2)|$, where 
$|\ \cdot \ |$ denotes the Euclidean distance. 
Note that a loop can be a part of another larger loop formed at some other 
vertex. 
By definition the paths in $W_N\cup V_N$ do not have any loops 
with diameter greater than or equal to $1$. 
Let $\Gamma _N$ be the set of loopless paths on $F_N$ from $O$ to $a$:
\[\Gamma_0=\{(O,a), (O, b, a)\},\]
\[\Gamma _N= \{\ (w(0), w(1), \cdots , w(n)) \in W_N \cup V_N:
\  w(i) \neq w(j), \ 0\leq  i< j\leq n,\ n\in \nintegers \  \} .
\]
Note that any loopless path from $O$ to $a$ is confined in $\triangle Oab$.

\vspace{0.5cm}\parr
{\bf Loop erasure on} ${\bf F^V_1}$

We shall now describe the  loop-erasing procedure for paths in  $W_1\cup V_1$:
\itmb
\item[(i)] Erase all the loops formed at $O$;
\item[(ii)] Progress one step forward along the path, and 
erase all the loops at the new position;
\item[(iii)] Iterate this process, taking another step forward along the path and erasing the 
loops there, until reaching $a$.
\itme

To be precise, for $w \in W_1\cup V_1$, define the recursive sequence $\{s_i\}_{i=0}^{n}$, 
\[s_0=\sup \{j:w(j)=O \},\]
\[s_i=\sup \{j:w(j)=w(s_{i-1}+1)\}.\]
If $s_i>s_{i-1}+1$, then 
$[w(s_{i-1}+1), w(s_{i-1}+2), \ldots , w(s_{i}-1), w(s_{i})]$
forms a loop or multiple loops at 
$w(s_{i-1}+1)=w(s_{i})$, so we erase this part by removing  
$w(s_{i-1}+1), w(s_{i-1}+2), \ldots , w(s_{i}-2)$, and $w(s_{i}-1)$.
If $w(s_n)=a$,  then 
we have obtained a loop-erased path,
\[L w =[w(s_0), w(s_1), \ldots ,w(s_n)]\in 
\Gamma _1, \]
where $L: W_1\cup V_1 \to \Gamma _1$ is the loop-erasing operator.
\QED

\vspace{0.5cm}\par
Fig. 3 shows all the possible loopless paths from $O$ to $a$ on $F_1$. 
Here only the parts in $\triangle Oab$ are shown, for 
any path cannot go into the other triangles without making a loop.

\begin{figure}[htb]
\begin{center}
\input{Fig4star.tex}
\\[1\baselineskip]
\caption{Loopless paths from $O$ to $a$ on $F_1$}
\end{center}
\end{figure}
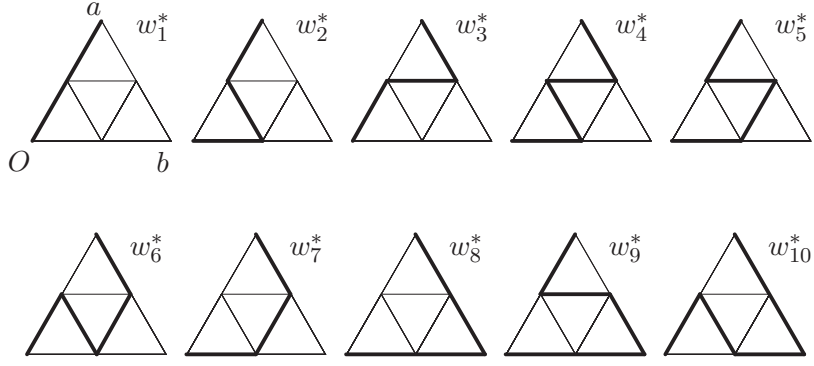

Note that $w \in W_1$ implies $Lw \in W_1\cap \Gamma _1$, but that 
 $w \in V_1$ can result in  $Lw \in W_1 \cap \Gamma _1$, with $b$ being erased 
 together with a loop.
So far, our loop-erasing procedure is the same as the chronological method defined for paths on 
${\mathbb Z}^d $ in \cite{Lawler}. 

For a general $N$, 
we erase loops from the largest scale loops down, repeatedly applying the 
loop-erasing procedure on $F^V_1$.

\vspace{0.5cm}\parr
{\bf First step of the induction -- erasing largest scale loops}

We shall illustrate the first step of loop erasure.
Decompose a  path $w\in W_N \cup V_N$ into $(Q_1w; w_1, \cdots , w_{\ell (Q_1w)})$, 
$w_i \in W_{N-1} \cup V_{N-1}$ 
$i=1, \cdots , \ell(Q_1w)$ as in \eqnu{decomposition1}. 
Fig. 4(a)  shows  $w \in W_N \cup V_N$ and Fig. 4(b) shows $Q_1w$.  
Erase all the loops in chronological order 
from   $Q_1w \in W_1\cup V_1$ to obtain $LQ_1w$ as in Fig. 4(c), then restore the original fine structures to the remaining parts as shown in Fig. 4(d). 
That is, if we write 
\[L Q_1w =[w(T_0^1),  w(T_{s_1}^1), \ldots ,w(T_{s_n}^1)], \  n=\ell (LQ_1w),\]
for each $i$, 
fit the path segment $w_{s_i+1}=(w(T_{s_i}^1),  w(T_{s_i}^1+1), \cdots , w(T_{s_i+1}^1) )$
between $w(T_{s_i}^1)$ and $w(T_{s_{i+1}}^1)$ of $LQ_1w$.
We call the path obtained at this stage 
$\tilde{L}w$.
Notice that in this stage all the loops with diameter greater than $1/2$ have been erased.
Let $\hat{Q}_1w=LQ_1w$. This completes the first induction step.
\QED

\begin{figure}[htb]
\begin{center}
\input{Fig3a.tex}
\\[1\baselineskip]
\caption{The loop-erasing procedure: (a) $w$, (b) $Q_1w$, (c) $LQ_1w=\hat{Q}_1w$, (d) $\tilde{L}w$}
\end{center}
\end{figure}
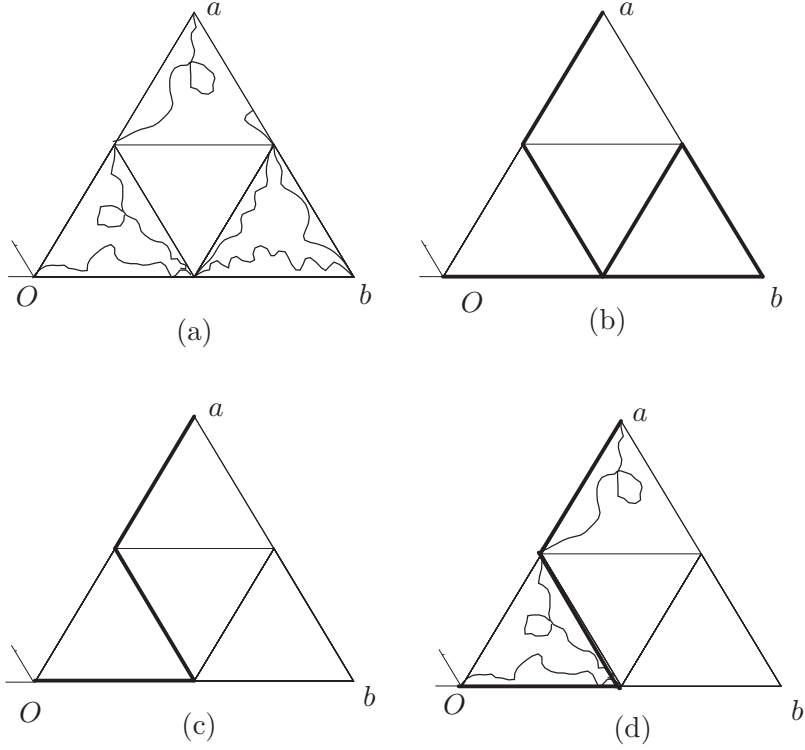
\vspace{0.5cm}\parr

\vspace{0.5cm}\par
The idea is  to repeat a similar procedure within each $2^{-1}$-triangle to erase all 
loops with diameter greater than $1/4$, and then within 
each $4^{-1}$-triangle, and so on, until there remain no loops. 
To describe next induction steps more precisely,  
we make some preparations.
For $w \in  W_N$ and $M \leq N$, we shall define 
the  sequence $(\Delta_1, \ldots , \Delta_k)$ of  the  $2^{-M}$-triangles 
$w$ `passes through',  and 
their exit times $\{T_i^{ex,M}(w)\}_{i=1}^{k}$ 
as a subsequence of  $\{T_i^{M}(w)\}_{i=1}^{m}$ as follows:
Let  $T_0^{ex, M}(w)=0$. 
There is a unique element of ${\cal T}_M$ that contains $w(T_0^M)$ and $w(T_1^M)$, 
which we denote by  $\Delta _1$. 
For $i\geq 1$, define 
\[J(i)=\min \{j \geq 0\ : \ j<m,\ T_j^M(w)>T_{i-1}^{ex, M}(w),\ 
w(T_{j+1}^M(w))\not\in \Delta _i\},\] 
if the minimum exists, 
otherwise $J(i)=m$.   
Then define $ T_i^{ex, M}(w)=T_{J(i)}^{M}(w)$, and let $\Delta _{i+1}$ be
the unique $2^{-M}$-triangle that contains both $w(T_i^{ex, M})$ and $w(T_{J(i)+1}^{M})$. 
By definition, we see that  $\Delta _{i} \cap 
\Delta _{i+1 }$ is a one-point set $\{w(T_i^{ex, M})\}$, for $i=1, \ldots , k-1$. 
We denote the sequence of these triangles by 
$\sigma _M(w)=( \Delta _1, \ldots , \Delta _k)$, 
and call it the ${\bf 2^{-M}}${\bf -skeleton} of $w$.
We call the sequence  
\par\noindent
$\{T_i^{ex, M}(w)\}_{i=0, 1, \ldots , k}$ 
{\bf exit times} from the triangles in the skeleton.  
For each $i$, there is an $n=n(i)$ such that 
$T_{i-1}^{ex, M}(w)=T^M_{n}(w)$.  If $T_{i}^{ex, M}(w)=T^M_{n+1} (w)$,
we say that $\Delta _i \in \sigma _M(w)$ is  {\bf Type 1},  
and if  $T_{i}^{ex, M}(w)=T^M_{n+2}(w)$,  {\bf Type 2}.
If $w\in \Gamma _N$ and $M\leq N$, then 
its $2^{-M}$-skeleton is a collection of distinct $2^{-M}$-triangles 
and each of them is either Type 1 or Type 2.  
Assume $w\in W_N\cup V_N$ and $M\leq N$.
For each $\Delta $ in $\sigma _M(w)$, the {\bf path segment of} ${\bf w} $ {\bf in}
 ${\bf \Delta }$ is 
defined by 
\eqnb
\eqna{segment}
w|_{\Delta }=[w(n), \ T^{ex, M}_{i-1}(w) \leq n \leq T^{ex, M}_{i}(w)].
\eqne 
Note that the definition of $T^M_i$ allows a path segment $w|_{\Delta }$ to leak into the 
neighboring 
$2^{-M}$-triangles.  

If $Q_Mw \in \Gamma _M$,  then 
$w|_{\Delta } \in W_{N-M}$ or $w|_{\Delta } \in V_{N-M}$ (identification implied),  
according to the type of $\Delta \in \sigma_M(w)$, 
where the entrance to $\Delta $ is identified with $O$ and the exit with $a$. 
This means that each $w $ 
satisfying  $Q_Mw \in \Gamma _M$
can be decomposed uniquely to 
\eqnb 
\eqna{decomposition2}
(\sigma_M(w); \ w|_{\Delta_1}, \cdots , w|_{\Delta_k}), \ 
w|_{\Delta _i } \in W_{N-M}\cup V_{N-M},\ i=1, \cdots , k.
\eqne
Conversely, 
given a collection of distinct $2^{-M}$-triangles  $\{\Delta _i\} _{i=1}^{k}$  
such that 
$O\in \Delta_1$, 
$a \in \Delta_k$, 
$\Delta _i$ and $\Delta _{i+1}$ are neighbors,  and  
$w'_i \in W_{N-M}\cup V_{N-M}$, $i=1, \cdots , k$, then  
we can assemble them to obtain a unique element $w$ of $W_N\cup V_N$.

We call a loop $[w(i), w(i+1), \cdots , w(i+i_0)]$ a ${\bf 2^{-M}}${\bf -scale loop} whenever 
there exists an $M \in \pintegers $ such that  
\[\min \{N' : w(i)=w(i+i_0)\in G_{N'}\}=M,\ d\geq 2^{-M},\]
where $d$ is the diameter of the loop. 

Using above as a base step, we shall now describe the induction step of our operation:

\vspace{0.5cm}\parr
{\bf Induction step}

Let $w\in W_N \cup V_N$. 
For $1\leq M\leq N$,  
assume that all of the $2^{-1}$ to $2^{-M}$-scale loops have been erased from  
$w$, and denote by $w' \in W_N\cup V_N$ the path obtained at this stage.
Then $Q_Mw' \in \Gamma _M$.

\itmb
\item[1)]
Since $Q_Mw' \in \Gamma _M$, we have the decomposition of $w'$: 
 $(\sigma _M(w') ; \ w'_1, \cdots w'_k )$,
$w'_i \in W_{N-M} \cup V_{N-M}$ as
 given in \eqnu{decomposition2}.

\item[2)] 
From each $w_i' $, 
erase $2^{-1} $-scale 
loops (largest scale loops) according to the base step procedure above to obtain $\tilde{L}w_i' \in W_{N-M}\cup V_{N-M}$ and 
$\hat{Q}_{1} w_i' \in \Gamma_1$. 

\item[3)]
Assemble $(\sigma _M(w'); \ \tilde{L}w'_1, \cdots ,\tilde{L} w'_k)$ and 
 $(\sigma_M(w');\ \hat{Q}_{1} w'_1, \cdots , \hat{Q}_{1} w'_k)$ 
to obtain $w'' \in W_N \cup V_N$ and
 $\hat{Q}_{M+1}w \in \Gamma_{M+1} $, respectively. $w''$ has no 
$2^{-1}$ to $2^{-(M+1)}$-scale loops. 
\itme
\QED

\vspace{0.5cm}\par
We then continue this operation until we have erased all of the loops to have $Lw=\hat{Q}_Nw \in \Gamma_N$.
In this way, the loop erasing operator $L$ defined for $W_1\cup V_1$ has been extended to 
$L:\bigcup _{N=1}^{\infty } (W_N\cup V_N) \to \bigcup _{N=1}^{\infty }\Gamma _N$ 
with $L(W_N\cup V_N)=\Gamma _N$.
Notice that the operation described above is essentially a repetition of loop-erasing  for $W_1\cup V_1$.
$\hat{Q}_M$ is a map from $\bigcup_{N=M} ^{\infty} (W_N\cup V_N)$ 
to $\Gamma _M$.
In the induction step, we observe that $\hat{Q}_{M+1}w=Q_{M+1}w''$.  
Although  it may occur that 
$\sigma _{M+1}(w'') \neq \sigma _{M+1}(w')$ because of the erasure of 
$2^{-(M+1)}$-scale loops,  
it holds that $\sigma _M(w'')=\sigma _M(w')$, which
can be extended to $\sigma _K(w')=\sigma _K(w'')$ for any $K\leq M$.
We remark that  the procedure implies that for any  $w\in W_N \cup V_N$,

\eqnb
\eqna{sigmainv}  
\sigma _K(\hat{Q}_{M}w)=\sigma _K (\hat{Q}_{K}w) \ \mbox{ for any } N\geq M \geq K .\eqne
In particular, 
\eqnb
\eqna{sigmainv2}  
\sigma _K(Lw)=\sigma _K (\hat{Q}_{K}w) \ \mbox{ for }  K \leq N.\eqne
i.e., in the process of loop-erasing, once loops of $2^{-K}$-scale and greater 
have been erased, the $2^{-K}$-skeleton does not change any more.  However, it should be noted 
that the types of the triangles can change from Type 2 to Type 1.


We induce measures $\hat{P}^u_N=P^u_N\circ L^{-1}$ and ${\hat{P}}'^u_N={P}'^u_N\circ L^{-1}$, which satisfy 
$\hat{P}^u_N[\Gamma_N]=1$ and ${\hat{P}}'^u_N [\Gamma _N]=1$.  
For $w^*_1, \cdots , w^*_{10}$ shown in Fig. 3, 
denote 
\eqnb \eqna{pq}
p_i= \hat{P}^u_1 [w^*_i] =P^u_1[w: Lw=w^*_i], \ \ 
q_i= {\hat{P}}'^u_1 [w^*_i] ={P}'^u_1[w: Lw=w^*_i].
\eqne 
$p_i$ and $q_i$ can be obtained as explicit functions of $u$ and $x_u$ by direct, but lengthy calculations, which are shown in Appendix.  
In the case that $u=1$ (the ordinary loop-erased random walk), we have 
$x_1=1/4$, $p_1=1/2$, 
$p_2=p_3=p_7=2/15$, $p_4=p_5=p_6=1/30$, $p_8=p_9=p_{10}=0$, 
$q_1=1/9$, $q_2=q_3=11/90$, $q_4=q_5=q_6=2/45$, $q_7=8/45$, $q_8=2/9$ 
and $q_9=q_{10}=1/18$ as in \cite{HM}.
For $u=0$, we have $p_1=x_0$, $p_7=x_0^2$ and $p_i=0$ otherwise, 
and $q_1=x_0^4$, $q_2=q_3=x_0^3$, $q_8=x_0^2$ and $q_i=0$ otherwise, 
with $x_0=(\sqrt{5}-1)/2$ as
in \cite{hhk1}.  

 $\hat{P}^u_N$ and  ${\hat{P}}'^u_N$ define two families of walks on $F_N$ obtained 
by erasing loops from $Z^u_N$ and  ${Z}'^u_N$, respectively. 
We remark that that $\displaystyle \frac{2}{3}
\hat{P}^1_N+\frac{1}{3}{\hat{P}}'^1_N$ equals to the `standard' LERW studied 
in \cite{STW}.

An important observation is that in the process of erasing loops from 
$Z^u_{N+1}$, if we stop at the point where we have obtained
$\hat{Q}_NZ^u_{N+1}$, it is nothing but the procedure for obtaining $LZ^u_N$ from
$Z^u_N$.  The same holds also for ${Z}'^u_{N+1}$.  This can be expressed as:
\eqnb
\eqna{qhat}
P_{N+1}^u[\{v':\hat{Q}_Nv'=v\}]=\hat{P}_N^u[v],\ \ 
{P}'^u_{N+1}[\{v':\hat{Q}_Nv'=v\}]={\hat{P}}'^u_N[v].
\eqne    
In this stage what is left to do for obtaining $LZ^u_{N+1}$ from 
$\hat{Q}_NZ^u_{N+1}$ is a sequence of loop-erasing from $Z^u_1$ or ${Z}'^u_1$.
This combined with \eqnu{qhat} leads to a `decomposition' of 
LERW measures.  For $w\in \Gamma _{N+1}$,  
 \begin{eqnarray*}
\hat{P}_{N+1}^u[\ w\ ]&=&
\sum _{v\in \Gamma_N} P^u_{N+1}[\{v':  Lv'=w\} |\ \hat{Q}_N v'=v\ ]
 \ P^u_{N+1}[\{v': \hat{Q}_Nv'=v\} ]
\\
&=&
\sum _{v \in \Gamma_N} \ (\prod _{i=1}^{k} \hat{P}^{*u}_{1}[\ w_i\ ])
\ \hat{P}_{N}^u[\ v\ ],
\end{eqnarray*}
 where  $\sigma _N (v)=(\Delta _1, \cdots , \Delta_k)$, 
$w_i=v|_{\Delta _i}$ (identification implied),  
$\hat{P}^{*u}_{1}=\hat{P}^u_1$ if $\Delta _i$ is  Type 1, and 
$\hat{P}^{*u}_{1}={\hat{P}}'^u_1$ if $\Delta _i$ is  Type 2. 
A similar decomposition holds also for ${\hat{P}}'^u_{N+1}$. 
This is the key to the recursion relations of generating functions 
defined below. 
 
For $w\in \Gamma _N$, let us denote  the number of $2^{-N}$- triangles of Type 1, 
(the path passes two of the vertices) and
those of Type 2 (the path passes all three vertices) in $\sigma _N(w)$  by $s_1(w)$ and $s_2(w)$, respectively.
Note that $\ell (w)=s_1 (w)+2 s_2(w)$. 
Define two sequences, $\{\hat{\Phi} _N\}_{N \in {\mathbb N}}$ and $\{\hat{\Theta}_N\}_{N \in {\mathbb N}}$, of 
generating functions by:  
\[\hat{\Phi}_N (x,y)=\sum _{w\in \Gamma _N }\hat{P}_N^u(w)x^{s_1(w)}y^{s_2(w)},\]
\[\hat{\Theta}_N (x,y)=\sum _{w\in \Gamma _N }{\hat{P}}'^u_N(w)x^{s_1(w)}y^{s_2(w)},\ \ \ x,y\geq 0.\]
For simplicity, we shall denote $\hat{\Phi}_1 (x,y)$ and $\hat{\Theta}_1 (x,y)$ by 
$\hat{\Phi} (x,y)$ and $\hat{\Theta}(x,y)$ and omit writing $u$-dependence 
explicitly. 
Similar to Proposition 3 in \cite{HM}, 
we have

\prpb
\prpa{recursion}

The above generating functions satisfy the following recursion relations for all $N \in {\mathbb N}$ :
\[\hat{\Phi}(x,y)= p_1 x^2 +(p_2+p_3) xy +p_4y^2 +(p_5+p_6) x^2y + p_7x^3,
\]
\[\hat{\Theta} (x,y)=q_1 x^2 +(q_2+q_3) xy +q_4y^2 +(q_5+q_6) x^2y + q_7x^3+q_8x^2y+
(q_9+q_{10})xy^2,\]

\[\hat{\Phi}_{N+1}(x,y)=\hat{\Phi}_{N}(\hat{\Phi} (x,y), \hat{\Theta}(x,y)),\]
\[\hat{\Theta}_{N+1}(x,y)=\hat{\Theta }_{N}(\hat{\Phi} (x,y), \hat{\Theta}(x,y)),\]
where $p_i= \hat{P}^u_1 [w^*_i] $ and $q_i= {\hat{P}}'^u_1 [w^*_i]$, $i=1,2, \cdots ,10$. 
\prpe

\vspace{0.5cm}\par
Define the mean matrix by
\eqnb
\eqna{matrix}
{\bf M}=\left [
\begin{array}{cc}
\frac{\partial}{\partial x}\Phi (1,1) 
& \frac{\partial}{\partial y}\Phi (1,1) 
\\
\frac{\partial}{\partial x}\Theta (1,1) 
& \frac{\partial}{\partial y}\Theta (1,1)
\end{array}
\right ].
\eqne
It is a strictly positive matrix, and the 
larger eigenvalue $\lambda =\lambda (u)$ is a continuous 
function of $u$, satisfying $2<\lambda <3$.   

Let $Z^u_N$ and ${Z}'^u_N$ be as in \eqnu{Z} and \eqnu{ZZ}.  
The loop-erasing procedure together with the structure of the Sierpi\'nski gasket 
 leads to (Proposition 4 in \cite{HM})
\parr
\prpb
\prpa{traverse} 
Let $M\leq N$.
Conditioned on $\sigma _M(LZ^u_N)=(\Delta _1, \ldots , \Delta _k)$ and the type of 
each element of the skeleton, the traverse times of the triangles 
\[T_i^{ex, M}(LZ^u_N)-T_{i-1}^{ex, M}(LZ^u_N), \ \ i=1,2, \ldots , k\]
are independent.  Each of them  has the same distribution as either 
$T^{ex, N-M}_1(LZ^u_{N-M})$ or \par\noindent 
$T^{ex, N-M}_1(L{Z}'^u_{N-M})$, according to whether 
$\Delta _i$ is of Type 1 or Type 2. 
\prpe


\vspace{0.5cm}\parr
\section{The scaling limit}
\seca{SL}

In this section, we investigate the limit of the loop-erased 
self-repelling walks constructed 
in \secu{LE} as the edge length tends to $0$.
Since it is easier to deal with continuous functions from the beginning, 
we regard $F_N$'s and $F^V_N$'s as closed subsets of ${\mathbb R}^2$ made up of all the points on their edges.  
We define the {\bf {\sg}}  by 
$F=cl(\cup_{N=0}^{\infty } F_{N})$, where $cl$ denotes closure.  
We start with a larger space $F^V=cl(\cup_{N=0}^{\infty } F^V_{N})$ and let 
\[
C=\{w\in C([0,\infty)\rightarrow F^V) \ :\  w(0)=O,\ 
\lim_{t\rightarrow\infty} w(t)=a\}\,. 
\]
$C$ is a complete separable 
metric space with the metric
\[
d(u,v)=\sup_{t\in[0,\infty)}\, |u(t)-v(t)|\,,\ 
u,v\in C,
\]
where  $|x-y|$, $x, y \in \twreals$,  denotes the Euclidean distance.
Hereafter, 
for $w \in \bigcup _{N=1}^{\infty }(W_N\cup V_N)$, we define  
\[w(t)=a,\ \ \ t\geq \ell (w),\]
and 
interpolate the path linearly, 
\[w(t)=(i+1-t)w(i)+(t-i)w(i+1),\ \   i \leq t <i+1,\ \ i \in \pintegers \]
so that we can regard $w$ as a continuous function on $[0, \infty)$. 
We shall regard $W_N$, $V_N$ and $\Gamma _N$ as subsets of  $C$.
Hitting times, $\{ T_i^M(w)\}_{i=1}^{m}$ are defined for $w\in C$ as in the previous sections, although the 
infimum is taken over continuous time:
\[T_0^M(w)=0,\ 
\ T_i^M(w)=\inf \{t>T_{i-1}^M(w)\ : \ w(t) \in G_M\setminus \{w(T_{i-1}^M(w))\} \}.\]
Notice that the condition $\lim _{t\to \infty} w(t)=a$ makes $\{T^M_i(w) \}_{i=0}^{m}$ a finite 
sequence.    

For $N \in \integers_{+}$, we define a coarse-graining  map
$Q_{N}:C \rightarrow C$ by $(Q_{N}w)(i)=w(T_{i}^{N}(w))$ 
for $i=0,1,2,\ldots, m$,  and by using linear
interpolation 
\[ (Q_{N}w)(t) = \left\{ \begin{array}{ll} (i+1-t)\ (Q_{N}w)(i)&+(t-i)\
(Q_{N}w)(i+1),\\
 & i \leq t < i+1,\ i=0,1,2,\ldots, m-1, \\
 a, & t \geq m. \end{array} \right. \]
We  define also   the $2^{-M}$-{\bf skeleton}, 
$\sigma _M(w)$  (a sequence of $2^{-M}$-triangles 
$w$ passes through), the {\bf exit times}  $\{ T_i^{ex, M}\}_{i=1}^{k}$ and  
types of triangles in a similar way 
to their counterparts in \secu{LE}.  
The loop-erasing operator is regarded as $L:\bigcup_{N=1}^{\infty}(W_N\cup V_N) \to 
\bigcup_{N=1}^{\infty}\Gamma_N$. 
$\hat{Q}_N$'s are as in \secu{LE} with resulting paths in $\Gamma_N$. 
$P_N^u$, ${P}'^u_N$, 
$\hat{P}_N^u$ and ${\hat{P}}'^u_N$ are regarded as probability measures on $C$.

In order to consider an almost sure limit, we shall couple walks on 
different pre-\sg s.
Let 
\[f_a(t)=\left\{  \begin{array}{ll}ta
& 0\leq t\leq 1,\\
a &t >1,  \end{array} \right. \]
where $\displaystyle a=(\frac{1}{2}, \frac{\sqrt{3}}{2})$, and 
\[\Omega '=\{  \omega=(\omega_0, \omega_1, \omega_2, \cdots )\ :\ \omega_0=f_a,\ \omega_N\in \Gamma _N,\   
\omega_{N-1} \triangleright \omega_{N}, \ N \in {\mathbb N} \},
\]
where $\omega_N  \triangleright \omega_{N+1}$ means that there exists a $v \in W_{N+1}\cup V_{N+1}$ such 
that $\hat{Q}_Nv=\omega_N$ and $Lv=\hat{Q}_{N+1}v=\omega_{N+1}$.  Namely, 
$v$ is a path obtained by adding a finer, $2^{-(N+1)}$ - scale structure 
(not loopless yet)  to $\omega_N$, 
and erasing  $2^{-(N+1)}$ - scale loops from $v$ gives $\omega_{N+1}$.
We assumed $\omega_0=f_a$ here, for we can deal with the case $\omega_0=f_b$
with 
\[f_b(t)=\left\{  \begin{array}{ll}tb
& 0\leq t\leq 1,\\
b+(t-1)(a-b)
& 1< t \leq 2\\
a &t >2,  \end{array} \right. \]
where $\displaystyle b=(1,0)$, in a similar way.

Define the projection onto the first $N+1$ elements by 
\[\pi _N \omega =(\omega_0, \omega_1, \ldots , \omega_N). \]
For each $u \in [0,1]$, define a probability measure $\tilde{P}_N$ on $\pi _N\Omega '$ by 
\[\tilde{P}_N[(\omega_0, \omega_1, \ldots , \omega_N)]=P_N^u[\ v :\ \hat{Q}_i v=\omega_i , i=0, \ldots , N\ ],\]
where $P^u_N$ is defined in \secu{SRW}. 
Although $\tilde{P}_N$ depends on $u$, we shall not write the $u$-dependence 
explicitly for simplicity.  
The following consistency condition is a direct consequence of the  loop-erasing procedure: 
\eqnb
\eqna{consistency}
\tilde{P}_N [(\omega_0, \omega_1, \ldots , \omega_N) ]=\sum _{\omega '} \tilde{P}_{N+1} [(\omega_0, \omega_1, \ldots , \omega_N, \omega ' ) ],
\eqne 
where the sum is taken over all possible $\omega ' \in \Gamma _{N+1}$ such that $\omega _N\triangleright \omega '$.



By virtue of \eqnu{consistency} and Kolmogorov's extension theorem for a projective limit, 
there is a probability measure $P$ on $\Omega _0=
C^{\nintegers}=C \times C \times \cdots \ $  such that
\[P[\ \Omega '\ ]=1,\]
\[  P\circ \pi _N^{-1} =\tilde{P}_{N}, \ N\in {\mathbb Z}_+, \]
where $\pi _N$ denotes the projection onto the first $(N+1)$ elements also here.

Define $Y_N :\ \Omega _0 \to C$ 
by 
\[Y_N((\omega_0, \omega_1, \ldots )) =\left\{  \begin{array}{ll} \omega _N
& \mbox{ if } (\omega_0, \omega_1, \ldots ) \in \Omega ' ,\\
f_a & \mbox{ otherwise }.  \end{array} \right. \]
Then $Y_N$ is an $F$-valued process $Y_N(\omega , t )$ on 
$(\Omega _0,  {\cal B}, P)$, 
where ${\cal B}$ is the Borel
algebra on $\Omega _0$ generated by the cylinder sets.  Then we have 
$P\circ Y_N^{-1}=\hat{P}_N^u$.
For $N\geq M$ and $\Delta \in {\cal T}_M$, denote the the path segment of $Y_N$ in 
$\Delta $ by $Y_N|_{\Delta}$ as the continuous version of \eqnu{segment}.

For $w \in \bigcup_{N=1}^{\infty}\Gamma _N$ and $j=1,2$, denote by $S_j^M(w)$ the number of $2^{-M}$-triangles of Type $j$ in $\sigma _M(w)$, 
and  ${\bf S}^M(w)=(S_1^M(w), S_2^M(w))$.
Note that if $w\in \Gamma _N$, then $\ell (w)=S_1^N(w)+2 S_2^N(w)$.

Let ${\bf S}=(S_1,S_2)$ and  ${\bf S}'=(S'_1,S'_2)$ be $({\mathbb Z}_+)^2$-valued random variables on 
$(\Omega _0,  {\cal B}, P)$ 
with the same distributions as those of $(S_1^1, S_2^1)$ under $\hat{P}_1^u$ and  under  ${\hat{P}}'^u_1$, respectively.   

\prpb
\prpa{brpr}
Fix arbitrarily $v \in \Gamma _M$, and let $\sigma _M(v)=(\Delta _1, \ldots , \Delta _k)$.
For each $i$, $1\leq i \leq k$,   
under the conditional probability $P[\ \ \cdot \ \ |Y_M=v ]$, 
$\{{\bf S}^{M+N}(Y_{M+N}|_{\Delta _i}),\ N=0,1,2, \cdots  \}$ is a two-type 
supercritical branching process, with the types of children corresponding  to the types of triangles.
The offspring distributions born from a Type 1 triangle 
and from a Type 2 triangle are equal to those of 
${\bf S}$ and ${\bf S}'$, respectively.  
If $\Delta _i$ is  Type 1, the process starts in state $(1,0)$, and  
if $\Delta _i$ is  Type 2,  in state $(0,1)$.

\begin{enumerate} 
\item[(1)]
The generating functions for the offspring distributions are 
\[E[\ x^{S_1}y^{S_2}\  ]=\hat{\Phi} (x,y),\]
\[E [\ x^{S'_1}y^{S'_2}\ ]=\hat{\Theta}  (x,y),\]
where $E$ is the expectation with regard to $P$. 
\item[(2)] 
Let ${\bf M}$ be the mean matrix given by \eqnu{matrix}.

Then
\[E[\ {\bf S}^{M+N}(Y_{M+N}|_{\Delta _i})\  |\ Y_M=v\  ]=
{\bf S}^M(v|_{\Delta _i}) {\bf M}^N.\]
\item[(3)]
$P[S_1+S_2 \geq 2]=P[S_1'+S'_2 \geq 2]=1$ (non-singularity).
\item[(4)]
$E[\ S_i\log  S_i\ ]<\infty ,\ \ E[\ S'_i\log  S'_i\ ]<\infty ,
\ i=1,2.$
\end{enumerate}
\prpe

\prpu{brpr} suggests that we should consider 
the time-scaled processes: 
\[X_N(\ \omega ,\ \cdot \ )=Y_N(\ \omega ,\  \lambda^N  \cdot \ ), \ \ N \in {\mathbb Z}_+ ,\]
where $\lambda $ is the  larger eigenvalue of the mean matrix.

\prpb
\prpa{X}
For $M\leq N$, the following holds:
\[\sigma _M(X_N)=\sigma _M(X_M)=\sigma _M(Y_M), \ \ \mbox{ a.s.}\]
\parr
and 
\eqnb \eqna{YY}
X_N(T_i^{ex, M}(X_N))=X_M(T_i^{ex, M}(X_M))=Y_M(T_i^{ex, M}(Y_M)), \ \ \mbox{ a.s.}
\eqne

\prpe
\vspace{0.5cm}\par
Note that 
if $\sigma _M(X_N)=(\Delta _1, \cdots , \Delta _k)$, then
\[T_j^{ex, M}(X_N)=\lambda ^{-N}\sum _{i=1}^{j}(S_1^N(X_N|\Delta _i)+2S_2^N(X_N|\Delta _i)), \ \ 
1\leq j \leq k.\]


Let ${\bf u}={}^t(u_1, u_2)$ and ${\bf v}=(v_1, v_2)$ be the right and left positive eigenvectors associated with $\lambda $ 
such that $({\bf u}, {\bf v})=1$ and $({\bf u}, {\bf 1})=1$.
\prpb
\prpa{supbr}
Fix arbitrarily $v \in \Gamma _M$, and let $\sigma _M(v)=(\Delta _1, \ldots , \Delta _k)$.
For each $i$, $1\leq i \leq k$,   
under the conditional probability $P[\ \ \cdot \ \ |Y_M=v ]$, 
we have the following:
\itmb
\item[(1)]
$\{\lambda ^{-(M+N)}{\bf S}^{M+N}(X_{M+N}|_{\Delta _i}),\ N=0,1,2, \ldots \}$ converges a.s. 
as $N \rightarrow \infty$
to a ${\mathbb R}^2$-valued random variable ${\bf S}^{*M,i}=(S_1^{*M,i}, S_2^{*M,i})$.
\item[(2)]
$\{ {\bf S}^{*M,i}$, $i=1, \cdots , k\}$ are independent.
\item[(3)]
There are random variables $B_1$ and $B_2$ such that 
${\bf S}^{*M,i}$ is equal in distribution to $\lambda ^{-M}B_1{\bf v}$ 
if $\Delta _i$ is of Type 1, and equal in distribution to 
$\lambda ^{-M}B_2{\bf v}$ 
if $\Delta _i$ is of Type 2. 

\item[(4)] 
\[P[B_i>0]=1,\ \    
E[B_i]=u_i,\ \ i=1,2.\] 
$B_1$ and $B_2$ have strictly positive probability density functions.

\item[(5)]
 
The Laplace transform of  $B_i$, $i=1,2$    
\[g _i(t)=E[\exp ( -tB_i)], \ \ t\in {\mathbb C}\]
are entire functions on ${\mathbb C}$ and are the solution to 
\[g _1 (\lambda  t)=\hat{\Phi}( g _1( t), g _2 ( t)),\  \  
g _2 (\lambda  t)=\hat{\Theta} ( g _1( t), g _2( t)),\
g _1(0)=g _2(0)=1.
\]

\itme
\prpe

(1)--(4) in \prpu{supbr} are the straightforward consequences of general limit theorems
for supercritical multi-type branching processes (Theorem 1 and Theorem 2 in V.6 of \cite{AN}).
$P[B_i>0]=1$ is a consequence of $\hat{\Phi}$ and $\hat{\Theta}$ having no terms with 
degree smaller than $2$.  
For the existence of the Laplace transform on the entire ${\mathbb C}$, 
we need careful study of the recursions.  We omit the details here, since they are 
lengthy and similar to the proof of Proposition 4.5 in \cite{hhk2}.

 \vspace{0.5cm}\par
Let $T^{*M}_{i}=\sum_{j=1}^{i}(S_1^{*M,j}+2S_2^{*M,j})$, then  
\eqnb
\eqna{limT}
\lim _{N \to \infty}
T_j^{ex, M}(X_N)=T_j^{*M}.
\eqne
By virtue of \prpu{X} and \prpu{supbr}, we can prove the almost sure uniform convergence
for $X_N$. The proof here closely follows the argument of \cite{BP}. 

\thmb
\thma{asconv}
$X_N$ converges uniformly in $t$ a.s. as $N \rightarrow \infty$ to a 
continuous process $X$. 
\thme

\prfb
Choose $\omega \in \Omega'$ such that the following holds for all 
$M \in {\mathbb Z}_+$: 
$\dsp \lim _{N \to \infty} T_i^{ex, M}(X_N)=T_i^{*M}$ exists 
and $T_i^{*M}-T_{i-1}^{*M}>0$
for all  $1\leq i \leq k$, where $k=k_M$ denotes the number of triangles in $\sigma _M(Y_M)$.  
Let $R=T_1^{*0}+ \varepsilon$,
where $\varepsilon >0$ is arbitrary.  
It suffices to show that $X_N(\omega ,t)$ converges uniformly 
in $t\in [0,R]$.  In fact, if $t>R$, $X_N(t)=a$ for a large enough $N$.

Fix $M \in \pintegers $ arbitrarily.  
By expressing the arrival time at $a$ as the sum of traversing times of $2^{-M}$-triangles, 
we have 
$T^{ex, M}_{ k}(X_N)=T_1^{ex, 0}(X_N)$ a.s..  
Letting $N \rightarrow \infty$, 
we have $T^{*M}_{k}=T_{1}^{*0}$ a.s..   

The choice of $\omega$ implies that 
there exists an $N_1=N_1(\omega ) \in \nintegers$
such that
\eqnb
\eqna{max}
\max_{1 \leq i \leq k}|T_i^{ex, M} (X_N)-T_i^{*M}| \leq 
\min_{1 \leq i \leq k}
(T_i^{*M}-T_{i-1}^{*M}), \ \ |T_{k}^{ex, M}(X_N)-T_{k}^{*M}|<\varepsilon 
,\eqne   
for $N \geq N_1$.

If $0\le t< T^{*M}_{k}$, then choose
$j \in \{1, \cdots , k\}$ such that 
$T^{*M}_{j-1} \leq t < T^{*M}_{j}$.
Then \eqnu{max} implies  that $T_{j-2}^{ex, M}(X_N)\leq t \leq T_{j+1}^{ex, M}(X_N)$,
for $N \geq N_1$.  Since \prpu{X} shows 
\eqnb \eqna{XX}
X_N(T_j^{ex, M}(X_N))=X_M(T_j^{ex, M}(X_M)), \eqne
for all $N$ with $N\geq M$, we have 
\[|X_N(T_j^{ex, M}(X_N))-X_N(t)| \leq 3 \cdot 2^{-M}. \]
Otherwise, if $ T^{*M}_{k}\le t \le  T^{*M}_{k}+\varepsilon =R$, then 
let $j=k$.
Since $T_{k-1}^{ex, M}(X_N) \leq t$,
\[|X_N(T_j^{ex, M}(X_N))-X_N(t)| \leq  2^{-M}. \]

Therefore, if  $N,N' \geq N_1$, then for any $t \in [0,R]$,
\begin{eqnarray*}
\lefteqn
{|X_N(t)-X_{N'}(t)|
}\\
&\leq&
|X_N(T_j^{ex, M}(X_N))-X_{N}(t)|+|X_{N'}(T_j^{ex, M}(X_{N'}))-X_{N'}(t)|\\
&&
+|X_N(T_j^{ex, M}(X_N))-X_{N'}(T_j^{ex, M}(X_{N'}))|
\\
& \leq &
6 \cdot 2^{-M},
\end{eqnarray*}
where the third term in the middle part 
is $0$ by \eqnu{XX}.
Since $M$ is arbitrary, we have the uniform convergence.

\QED
\prfe

\prpu{supbr} (5) implies that $E[\exp t B_i ] <\infty $ for $t>0$, which 
leads to:   

\prpb
\prpa{stay}
\[P[ \mbox{There exist }  t_0<t_1 \mbox{ such that } X(t)=X(t_0)\neq a \mbox{ for all }  t \in [t_0, t_1] ]=0.\]
\prpe

The proof is similar to that in \cite{hh}.

\prpb
\prpa{exittime}

The following holds for all $M\in \pintegers$ almost surely: 
\itmb
\item[(1)]
$\sigma _M(X)=\sigma _M(X_M)$,
\item[(2)]
$X(T^{*M}_i)=X_M(T_i^{ex, M}(X_M))$, 
\item[(3)] Let $\sigma _M(X_M)=(\Delta_1, \cdots , \Delta_{k_M})$.
If $T_{i-1}^{*M}<t<T_i^{*M}$, then $X(t) \in \Delta _i \setminus G_M$, 
for all $1\leq i\leq k_M$.
In particular, $T_i^{*M}=T_i^{ex, M}(X)=T_i^M(X)$.
\itme
\prpe

\prfb
(1) and (2) are direct consequences of \prpu{X}, \eqnu{limT} and \thmu{asconv}. 

To prove (3),  let $v_i=X(T_i^{*M}) $, $i=1, \cdots , k_M$ and we first prove that if $T_{i-1}^{*M}<t<T_i^{*M}$,
then $X(t) \not\in \{v_{i-1},  v_i\}$, by showing none of the following events  $A_j$,  $j=1,2,3,4$ 
has positive probability.

\parr
$A_1$ : There exists  $t_1$, $T_{i-1}^{*M}< t_1<T_i^{*M}$ such that 
$X(t)=v_i$ for all $t_1<t \leq T_i^{*M}$ holds for some $i\in \{1, \cdots , k_M\}$.
\parr
$A_2$ : There exists $t_1$, $T_{i-1}^{*M}< t_1 < T_i^{*M}$ such that 
$X(t)=v_{i-1}$ for all $T_{i-1}^{*M} \leq t<t_1$ holds for some $i\in \{1, \cdots , k_M\}$.
\parr
$A_3$: There exist  $t_1$ and $t_2$,  $T_{i-1}^{*M}<  t_1<t_2<  T_i^{*M}$ such that 
$X(t_1)=v_i$ and $X(t_2) \neq v_i$ holds for some $i\in \{1, \cdots , k_M\}$.
\parr
$A_4$: There exist  $t_1$ and $t_2$,  $T_{i-1}^{*M}<  t_1<t_2<  T_i^{*M}$ such that 
$X(t_1)\neq v_{i-1}$ and $X(t_2) =v_{i-1}$ holds for some $i\in \{1, \cdots , k_M\}$.

\prpu{stay} guarantees that $P[A_1]=P[A_2]=0$. Since $X$ is the uniform limit of a sequences of
self-avoiding walks, we have $P[A_3]=P[A_4]=0$.

Let $\sigma =(\Delta_1,\cdots , \Delta_{k_M})$ be a sequence such that 
$P[\sigma_M(X)=\sigma ]>0$.  Let $\Delta _i$ be one of the triangles in $\sigma$,  
and denote the third vertex of $\Delta _i$ (neither the exit or entrance) by $v_i^*$.
We prove that 
the probability that $X$ hits  $v_i^*$ 
at  some  $T_{i-1}^{*M}<t<T_i^{*M}$ is zero. 
We can take a decreasing sequence of triangles $\{\Delta_i^{(K)}\}_{K=M}^{\infty}$ 
such that $\Delta_i^{(M)}=\Delta _i$, $\Delta_i^{(K)} \in {\cal T}_K$ (a $2^{-K}$-triangle), 
$\Delta_i^{(K)} \supset \Delta_i^{(K+1)}$,  $\ \bigcap_{K=M}^{\infty}\Delta_i^{(K)}
=\{v_i^{*}\}$.
Denote $\tilde{p}=\max \{\sum _{i=5}^{10}p_i, \sum _{i=5}^{10}q_i\} <1$, where 
$p_i$ and $q_i$ are defined by \eqnu{pq}. 
For any $K$, with $K\geq M$, (1) implies 
\begin{eqnarray*}
P[\ \Delta _i^{(K)} \in \sigma _K(X) \ |\ \sigma _M(X)=\sigma  \ ] 
&=&\ 
P[\ \Delta _i^{(K)} \in \sigma _K(X_K) \ |\ \sigma _M(X_K)=\sigma  \ ] \\ 
&\leq&
\tilde{p}^{K-M}.  
\end{eqnarray*}
Thus it follows that 
\[P[\ \Delta _i^{(K)} \in \sigma _K(X) \mbox{ for all } 
K\geq M\ |\ \sigma _M(X)=\sigma  \ ] =0\]
and 
\[P[\ \Delta _i^{(K)} \in \sigma _K(X) \mbox{ for all } 
K\geq M \mbox{ for some  } 1\leq i\leq k_M \ |\ \sigma _M(X)=\sigma  \ ] =0,\]
therefore,
\[P[\ \Delta _i^{(K)} \in \sigma _K(X) \mbox{ for all } 
K\geq M \mbox{ for some  } i \in \{1, \cdots , k_M \}]=0.\]
This implies that the probability that $X$ hits any `third' vertex of  
the triangles in its skeleton  is zero.  This completes the proof of (3).

\QED

This proposition further leads to ;

\prfe

\thmb
\thma{SA}

\itmb
\item[(1)]
$X$ is almost surely self-avoiding in the sense that
\[P[\ X(t_1) \neq X(t_2), \ 0\leq t_1\leq t_2 \leq T_a(X)\  ] =0,\]
where $T_a(X)=\inf \{t>0 : X(t)=a\}=T_1^{*0}$. 
\item[(2)]
The Hausdorff dimension of the path $X([0, T_a(X)])$ is almost surely equal to $\log \lambda /\log 2$, which is a continuous function of $u$.
\itme
\thme

(1) is a consequence of \prpu{stay} and \prpu{exittime}.
To calculate the Hausdorff dimension, we use the fact that 
if a path $w$ is self-avoiding, then it holds that 
\[\tilde{\sigma} _1(w) \supset \tilde{\sigma} _2(w)  \supset \tilde{\sigma} _3 (w) \supset \cdots \to w ,\]
in the Hausdorff metric, where $\tilde{\sigma}_M (w)$ is the union of all the closed $2^{-M}$-triangles in $\sigma _M(w)$. 
We could call the sample path a `random graph' directed recursive construction, 
for the numbers of similarity maps are random variables.  
We obtain the Hausdorff dimension by applying Thoerem 4.3 in \cite{KH}  
\QED


\section{Path properties of the limit process}
\seca{Path}

In this section we study some more sample path properties of the limit process.  We assume $0<u\leq 1$,  
for the case of $u=0$ is considered in  \cite{HH}.
We shall not explicitly write $u$-dependence as in the previous section.

Let 
\[\nu =\nu (u)=\frac{\log 2}{\log \lambda}.\]
Recall, from \prpu{supbr} (5) that  
\[g_i(t)=E[\ \exp (-tB_i)\ ],\ \ i=1,2\]
satisfy the functional equations:
\[g_1(\lambda t)=\hat{\Phi} (g_1(t), g_2(t)), \ 
g_2(\lambda t)=\hat{\Theta} (g_1(t), g_2(t)).\]
Let 
\[ h_i(t)=-t^{-\nu} \log g_i(t).\] 
 
The proof of the following proposition uses the explicit forms of  $\hat{\Phi}$ and $\hat{\Theta}$, but it basically 
follows those of  \cite{BP} and \cite{Kuma}.

\prpb
\prpa{LT}
There exist positive constants $C_1$, $C_2$ and $t_0$ such that
\[C_2 \leq h_i(t) \leq C_1\ \ i=1,2\]
hold for all  $t\geq t_0$.
\prpe

\prfb
We prove the upper bound for $i=1$.  
Combining  $g_1(\lambda t)=\hat{\Phi}(g_1(t), g_2(t))$ and the fact that  
$\hat{\Phi}(x,y)$ contains the term 
$p_1x^2$, we have $g_1(\lambda t)\geq p_1g_1(t)^2$, which implies 
$h_1(\lambda t)\leq \frac{a_2}{2}\  t^{-\nu}+h_1(t)$, where $a_2=-\log p_1>0$.
By induction, we have $h_1(\lambda ^n t)\leq a_2t^{-\nu}+h_1(t)$,
for any $t>0$ and $n\in {\mathbb N}$.  Fix $t_1>0$ arbitrarily. 
Since $h(t)$ is continuous for $t>0$,  $b_1:=\max _{t\in [t_1, \lambda t_1]} h_1(t)$
exists. 
For $t>\lambda t_1$, there is a positive integer $m$ and $s \in (t_1, \lambda t_1]$ 
such that $t=\lambda ^ms$.  Then 
$h_1(t)=h_1(\lambda ^ms)\leq a_2s ^{-\nu}+h(s)\leq a_2t_1^{-\nu}+b_1=:C_1$.
Thus we have  $h_1(t)\leq C_1$ for any $t\geq t_1$.
The proof for $i=2$ is similar, with the use of the term $q_4 y^2$ in 
$\hat{\Theta}(x,y)$. Note that $q_4>0$ for $u>0$. Take the larger $C_1$. 

To show the lower bound, first note that for $x,y \in [0, 1]$,
$\max \{\hat{\Phi} (x, y), \hat{\Theta} (x,y)\}\leq \max \{x^2, y^2\}$,
which leads to ${\hat{\Phi}(x, y)+\hat{\Theta}(x,y)}\leq 2(x+y)^2$.
Let $g(t):=g_1(t)+g_2(t)$, then $g(\lambda t)\leq 2g(t)^2$.
$\tilde{h}(t):=-t^{-\nu}\log g(t)$ satisfies $\tilde{h}(\lambda t)\geq 
(-(1/2)\log 2-\log g(t)) t^{-\nu}=-t^{-\nu}(1/2)\log 2 +\tilde{h}(t)$.
By induction, we have $\tilde{h}(\lambda ^n t) \geq  t^{-\nu}
(-\log 2-\log g(t))$.
Since $-\log g(t) \to \infty$ as $t \to \infty$, we can take $t_2>0$ such 
that $-\log 2-\log g(t)>1$ for all  $t\geq t_2$, which implies $\tilde{h}(\lambda ^nt)
\geq  t^{-\nu}$ for all $t\geq t_2$.
In a similar way to the proof above, we can show that
for any 
$t\geq t_2$, $\tilde{h}(t)\geq (1/2)t_2^{-\nu}=:C_2$, thus  
$h_i(t)=-t^{-\nu}\log g_i(t)\geq -t^{-\nu}\log g(t)=\tilde{h}(t) \geq C_2$ holds 
for both $i=1,2$.
Let $t_0=\max \{t_1, t_2\}$.  
\QED

\prfe

We now use a Tauberian theorem of exponential type. 
The following theorem, Corollary A.17 from \cite{tets} has a most suitable 
form for our purpose. 
\thmb
\thma{Taubelian}
Assume P is a Borel probability measure supported on $[0,\infty)$, 
and denote its Laplace transform by 
\[g(s)=\int_0^{\infty}e^{-s\xi}P[d\xi ] , \ \ s>0.\]
If there are constants  $C_1>0$, $C_2>0$ and $0<\nu <1$ 
such that 
\[-C_1 \leq \liminf _{s \to \infty }s^{-\nu } \log g(s)
\leq \limsup _{s \to \infty }s^{-\nu } \log g(s)
\leq -C_2,\]
then there exist $C_3>0$ and $C_4>0$ such that 
\[-C_3 \leq \liminf _{x\to 0} x^{\nu/(1-\nu )} \log P[[0,x]] 
\leq \limsup _{x\to 0} x^{\nu/(1-\nu )} \log P[[0,x]] 
\leq -C_4, \ \ x>0.
\]

\thme

\vspace{0.5cm}\par

Let $\tilde{B}_i=(v_1+2v_2)B_i, \ i=1,2,$
where ${\bf v}=(v_1, v_2)$ is the positive left eigenvector 
corresponding to $\lambda$ introduced just before \prpu{supbr}  in \secu{SL}.
Then \prpu{LT} and \thmu{Taubelian} lead to

\corb
\cora{PWx}
There exist positive constants 
 $C_5, C_6$, and $x_0$ such that 
\[ e^{-C_5 x^{-\frac{\nu}{1-\nu } }} \leq P[\ \tilde{B}_i \leq x \ ] \leq  
e^{-C_6 x^{-\frac{\nu}{1-\nu }}},\ \ \ i=1,2\]
hold for any $x \leq x_0$.
\core
 
\parr
{\bf Remark} In \cite{Jones},  supercritical multi-type branching processes
are studied and detailed results on the tail behavior of the limit processes 
are given, but our case does not satisfy the conditions for 
his results.

\prpb
\prpa{PXd}
There exist positive constants 
 $C_7, C_8$ and $K$ such that 
\[e^{-C_7 (\delta t^{-\nu})^{1/(1-\nu)} } \leq P[\ |X(t)|\geq \delta \ ]
\leq  P[\ \sup _{0\leq s \leq t} |X(s)|\geq \delta \ ]
\leq e^{-C_8 (\delta t^{-\nu})^{1/(1-\nu)} }
,\ \ \ i=1,2\]
hold for $\delta t^{-\nu } \geq K$.  
\prpe

\prfb

For an arbitrarily given $0<\delta <1$, 
take $N \in {\mathbb N}$ such that $2^{-N}< \delta \leq 2^{-N+1}$ holds. 
Recall that if $\Delta _1$, the first element of $\sigma _N(X)$, 
is of Type 1, $T^{ex, N}_1(X)$ has the same distribution as 
that of $\lambda ^{-N}\tilde{B}_1$, and if of Type 2,  
the same distribution as 
that of $\lambda ^{-N}\tilde{B}_2$.
For $i=1, 2$ denote by $A_i$ the event that  $\Delta _1$ 
is of Type $i$. 

For the upper bound, since  $\sup _{0\leq s \leq t} |X(s)|\geq \delta $ implies $T^{ex,N}_1(X)<t$,
\begin{eqnarray*}
P[\ \sup _{0\leq s\leq t}|X(s)|\geq \delta \ ]&\leq & P[\ T^{ex,N}_1(X)<t \ ] \\ 
&=&\ 
P[\ \tilde{B}_1<\lambda ^Nt \ ]\ P[\ A_1\ ] +
P[\ \tilde{B}_2<\lambda ^Nt \ ]\ P[\ A_2\ ]\\
&\leq&
 e^{-C_6 (\lambda ^N t)^{-\frac{\nu}{1-\nu }}}\\
&\leq &
 e^{-C_8 (\delta t^{- \nu} ) ^{1/(1-\nu )}}, 
\end{eqnarray*}
where we assumed that $\lambda ^Nt\leq x_0$ in the second 
inequality and set $C_8=2^{-1/(1-\nu)}C_6$. 

For the lower bound, since  $T^{ex, N-1}_1< t$ implies 
$|X(t)|\geq \delta $,  
we can show that there exists a $C_7>0$ such that 
\[P[\ |X(t)|\geq \delta \ ] \geq  e^{-C_7(\delta t^{- \nu} ) ^{1/(1-\nu )}}\] 
holds  
for  $\lambda ^{N-1}t\leq x_0$. 
Take $K=2x_0^{-\nu}.$
\QED
\prfe

\parr


\thmb
\thma{exp}
For any $p>0$, 
there are positive constants $C_9$ and $C_{10}$ such that 
\[C_9\leq \liminf _{t\to 0}\frac{E[|X(t)|^p]}{t^{p \nu }}\leq 
\limsup _{t\to 0}\frac{E[|X(t)|^p]}{t^{p \nu }}
\leq C_{10}.\]

\thme

\prfb
\prpu{PXd} implies that the following holds for small enough $t$:

\begin{eqnarray*}
\frac{1}{p}E[|X(t)|^p]&=&\int_0^1\delta ^{p-1}  P_i[|X(t)|\geq \delta ]\ d\delta   
\geq
\int _{Kt^{\nu }}^1
\delta ^{p-1}  P_i[|X(t)|\geq \delta ]\ d\delta   \\
&\geq&
\int _{Kt^{\nu }}^1
\delta ^{p-1} e^{-C_{7}(\delta t^{-\nu})^{1/(1-\nu)} } d\delta 
=t^{p\nu} \int _K^{  t^{-\nu}}y^{p-1}e^{-C_7y^{1/(1-\nu)}}\ dy  \\
&\geq &
\displaystyle \frac{1}{2} t^{p\nu} \int _K^{\infty}y^{p-1}e^{-C_7y^{1/(1-\nu)}}\ dy
=C_9t^{p\nu}.
\end{eqnarray*}

\begin{eqnarray*}
\frac{1}{p}E[|X(t)|^p]&=&\int_0^{Kt^{\nu}} \delta ^{p-1} 
 P_i[|X(t)|\geq \delta ]\ d\delta  +
\int_{Kt^{\nu}} ^1\delta ^{p-1} 
 P_i[|X(t)|\geq \delta ]\ d\delta  \\
&\leq&
\int_0^{Kt^{\nu}} \delta ^{p-1} 
 d\delta  +
t^{p\nu} \int _K^{\infty} y^{p-1}e^{-C_8y^{1/(1-\nu)}}\ dy
=C_{10} t^{p\nu}
 .\end{eqnarray*}
\QED
\prfe

\coru{PWx} and \prpu{PXd} lead to 
a law of the iterated logarithm.  Since the argument is similar to that in \cite{HHH}, 
we just give the statement below:  

\thmb
\thma{loglog}
There are positive constants $C_{11}$ and $C_{12}$  such that 
\[C_{11} \leq \limsup _{t\to 0} \frac{|X(t)|}{\psi (t)}\leq C_{12},  \mbox{ a.s.},\]
where $\psi (t)=t^{\nu}(\log \log (1/t))^{1-\nu}$.
\thme

\section{Conclusion and remarks}
\seca{CR}

We constructed a one-parameter family of self-avoiding walks that interpolates the SAW and the LERW on the Sierpi\'nski gasket, and proved that the scaling limit exists.  
The exponent that governs the short-time behavior and equals to the reciprocal of the path Hausdorff dimension 
is a continuous function of the parameter.  
Our construction has proved that the ELLF method does work for non-Markov
random walks as well as the simple random walk.  

Although we restricted ourselves to  $u \in [0,1]$ above,   
all the results hold also for $u>1$, that is, for self-attracting walks.  
By numerical calculations we observe that 
$\lambda $ is a decreasing function of $u$ and conjecture that as  $u\to \infty$,  
$x^*=\lim_{u \to \infty} ux_u$,  $p^*_i=\lim_{u\to \infty }p_i(x_u, u)$ and 
 $q^*_i=\lim_{u\to \infty }q_i(x_u, u)$ exist with $x^* \sim 0.351$, 
$p^*_1 \sim 0.206$ , $p^*_2 \sim 0.124$, $p^*_3 \sim 0.206 $,  
$p^*_4 \sim 0.352$, $p^*_5 \sim 0.083$, 
$p^*_6 \sim 0$, $p^*_7 \sim 0.029$,
$q^*_1\sim 0.345$, 
$ q^*_2 \sim 0.034,$ $ q^*_3 \sim 0.242,$ 
$q^*_4 \sim 0.097$, $q^*_5 \sim 0.208$,  $q^*_7\sim 0.073$
and $q^*_i \sim 0$ otherwise.

\vspace{0.5cm}\parr
{\Large\bf Acknowledgments}
\vspace{0.2cm}\parr
One of the authors (K. Hattori) would like to thank Ben Hambly for suggesting the 
problem of an interpolation between the SAW and the LERW, and Tetsuya Hattori for helpful discussion.
The authors would like to thank Wataru Asada for valuable discussion.  

\vspace{0.5cm}\parr
{\Large\bf Appendix }

\vspace{0.2cm}\parr
We show $p_i$, $q_i$ introduced in \secu{LE}
as explicit functions of $x_u$ and $u$.
Define for $ i=1,2, \cdots , 10$ 
\[p_i(x,u)={\sum}  _i \ u^{N(w)+M(w)}x^{\ell (w)-1},\ \ 
q_i(x,u)={\sum}'_i \ u^{N(w)+M(w)}x^{\ell (w)-1},\  x, u \geq 0\]
where the sum $\sum _i$ is taken over all $w \in W_1$ such that 
$Lw=w^*_i$ and  $\sum '_i$ over all $w \in V_1$ such that 
$Lw=w^*_i$.  Substituting $x=x_u$, we have
\[p_i=p_i(x_u, u) =\hat{P}_1^u[w^*_i], \ q_i=q_i(x_u, u) ={\hat{P}}'^u_1[w^*_i].\]
Let $U_1$ be a set of single loops formed at $O$ on $F_1$: 
\[U_1= \{\ w=(w(0), w(1), \cdots , w(n) ):
\  w(0) =w(n)=O, \   w(i) \in G_1\setminus G_0,\ 1\leq i\leq n-1, \]  
\[\{w(i), w(i+1)\} \in E_1, \  0 \leq i \leq n-1,  
\ n \in {\mathbb N}\ \}.\]
and define $N(w)$ by \eqnu{N}.  Define 
\[\Xi =\Xi  (x, u)=\sum _{w \in U_1}  
u^{N(w)} \ x^{\ell (w)}, \ x , u \geq 0.
\]
We obtain the explicit form ($\Theta$ in \cite{HHH}) as follows: 
\[\Xi (x, u)=\frac{2ux^2}{(1+ux)(1-2ux)}\{1+2(1-u^2)x^2
-2(1-u)^2ux^3 \}.
\]
We show the explicit forms of $p_i(x,u)$ and $q_i(x,u)$ below.
Each factor in these expressions represents a particular part of 
paths.  The common factor $1/(1-2u\Xi )$ comes from 
the sum over all the possible loops formed at $O$.
In the  lengthy expression of $q'_2(x,u)$, 
the first term is related to those paths with loops that are formed at  $(1/2, 0)$
and include $b$.  The factor   
$\displaystyle  x^2 \left(1+ux+\frac{u^2x^2}{(1+ux)(1-2ux)}
\{2(u^2-u+1)x+3\}\right)$ 
represents the part from the last hit at $O$ followed 
immediately by a step 
to $(1/2, 0)$ then to the first hit of $b$.  The factor $1/(1-2u\Xi )$
 stands for the sum over all the possible loops formed at $b$.
$\displaystyle  x^2  \left(1+\frac{u^2x(1+2x)}{(1+ux)(1-2x)}\right)$
corresponds to the trip back from $b$ to  $(1/2,0)$ followed 
immediately by a step to $(1/4, \sqrt{3}/4)$,
and 
$\displaystyle  \left(1+\frac{u^3x^2}{1-u^2x^2}+\frac{u^3x^4}{(1-u^2x^2)^2}\frac{1}{1-u(\ell +\Xi ) } \right)$ concerns the loops formed at  $(1/4, \sqrt{3}/4)$. 
The second term is related to paths whose first hit to $b$ occurs 
in a loop formed at $(1/4, \sqrt{3}/4)$.

\[p_1 (x, u)=\frac{x}{1-2u\Xi }\ \left(1+\frac{u^2x^2\{(1-u)^2x+2\}}{(1+ux)(1-2ux)} \right).\]

\[p_2 (x, u)=\frac{ux^2 }{1-2u\Xi }\ \left(1+\frac{u(1+u)x^2}{(1+ux)(1-2ux)}\right)\
 \left(1+\frac{u^3x^2}{1-u^2x^2}\right).\]

\[p_3 (x, u)=\frac{ux^2}{1-2u\Xi }\ 
 \left(1+\frac{u^2(1+u)x^2}{(1+ux)(1-2ux)}\right)\
 \left(1+\frac{ux^2}{1-u^2x^2}\right).\]

\[p_4 (x, u)=\frac{u^3x^3 }{(1-2u\Xi )(1-u^2x^2) }\ \left(1+\frac{u(1+u)x^2}{(1+ux)(1-2ux)}\right).\]

\[p_5 (x, u)= \frac{u^2x^3 }{(1-2u\Xi )(1-u^2x^2) }\ \left(1+\frac{u^2(1+u)x^2}{(1+ux)(1-2ux)}\right) .\]

\[p_6 (x, u)= \frac{u^2x^3 }{(1-2u\Xi )(1-u^2x^2) }\ \left(1+\frac{u(1+u)x^2}{(1+ux)(1-2ux)}\right) .\]

\[p_7 (x, u)=\frac{x^2 }{1-2u\Xi }\  \left(1+\frac{u^2(1+u)x^2}{(1+ux)(1-2ux)}\right)\
\left(1+\frac{u^3x^2}{1-u^2x^2}\right).\]

\[p_8 (x, u)=p_9 (x, u)=p_{10} (x, u)=0.\]


Let 
\[\ell  (x, u)=ux^2+\frac{ux^4}{1-u^2x^2},\ \ \Sigma  (x, u)=\frac{2ux^2}{1-ux},\]
and  
\[q_i (x, u)=\frac{1}{x^2} \frac{1}{1-2u\Xi } \ q'_i (x, u), \  \ i=1, \cdots , 10.\]
Then 
\[q'_1 (x, u)= \frac{(1+u)^2x^6}{1-2u\Xi } \left(1+\frac{u^2x(1+2x)}{(1+ux)(1-2ux)}\right)^2.\]
\[q'_2(x,u)= x^5 \left(1+ux+\frac{u^2x^2}{(1+ux)(1-2ux)}
\{2(u^2-u+1)x+3\}\right)\ \frac{1}{1-2u\Xi }  \left(1+\frac{u^2x(1+2x)}{(1+ux)(1-2x)}\right)\]
\[\times 
\left(1+\frac{u^3x^2}{1-u^2x^2}+\frac{u^3x^4}{(1-u^2x^2)^2}\frac{1}{1-u(\ell +\Xi ) } \right)
\]
\[
+ \frac{u^3x^7}{(1-u^2x^2)^2(1-u(\ell +\Xi ))}
\left(1+\frac{u(1+u)x^2}{(1+ux)(1-2ux)}\right)
.\]

\[q'_3 (x, u)= 
 \frac{2u(1+u)x^7}{1-2u\Xi } \left( 1+ \frac{u^2x(1+2x)}{(1+ux)(1-2ux)} \right)^2
\]
\[\times
\left\{ u+\frac{u^2x^2}{1-u^2x^2} +\frac{u}{1-u(\Sigma +\Xi )}
\frac{x^2}{1-ux}\left(1+\frac{u^2x}{1-ux}\right)  \right\}
\]

\[
+\frac{x^5}{1-u(\Sigma +\Xi )} \left(1+\frac{u^2(1+u)x^2}{(1+ux)(1-2ux)}\right) \frac{1}{1-ux}\left(1+\frac{u^2x}{1-ux}\right) 
.\]

\[q'_4 (x, u)= \frac{u^2x^6}{1-2u\Xi } \left(1+ux+
\frac{u^2x^2}{(1+ux)(1-2ux)}\left\{2(u^2-u+1)x+3 \right\}
\right)\left(1+\frac{u^2x(1+2x)}{(1+ux)(1-2ux)}\right)
\]
\[\times \frac{1}{1-u^2x^2}\left(1+\frac{u^2x^4}{(1-u^2x^2)(1-u(\ell +\Xi) )}\right)
\left(1+ \frac{x^2}{1-u(ux^2+\Xi )}\right)
\]
\[+\frac{u^4x^8}{(1-u^2x^2)^2(1-u(\ell +\Xi) )}\left( 1+\frac{u(1+u)x^2}{(1+ux)(1-2ux)}
\right)
\left(1+\frac{x^2}{1-u(ux^2+\Xi )}\right)
\]
\[+\frac{u^2x^6}{(1-u^2x^2)(1-u(ux^2+\Xi ))}\left( 1+ \frac{u(1+u)x^2}{(1+ux)(1-2ux)}\right)
.\]

\[q'_5 (x, u)=  \frac{u^3x^6}{1-2u\Xi } \left(1+ux+
\frac{u^2x^2}{(1+ux)(1-2ux)}\left\{2(u^2-u+1)x+3 \right\}
\right)\]
\[\times \left(1+\frac{ux(1+2x)}{(1+ux)(1-2ux)}\right)
 \left\{ \left( 1+\frac{ux^2}{1-u^2x^2}\right) \frac{ux^2}{(1-u^2x^2)
(1-u(\ell +\Xi) )} +\frac{1}{1-u^2x^2}\right\}
\]
\[+u^2x^6 \left(1+\frac{u^2(1+u)x^2}{(1+ux)(1-2ux)}\right) 
\left(1+\frac{ux^2}{1-u^2x^2} \right) \frac{1}{1-u(\ell +\Xi)}\ \frac{1}{1-u^2x^2}
.\]

\[q'_6 (x, u)= 2u(1+u)x^8\left(1+\frac{u^2x(1+2x)}{(1+ux)(1-2ux)} \right)^2 \frac{1}{1-2u\Xi }
\]
\[\times \left\{\left(1+\frac{u^2x}{1-ux}\right)\frac{u}{1-u(\Sigma +\Xi )}\frac{ux^2}{1-ux} 
+\frac{u}{1-u^2x^2}
\right\}\left(1+\frac{u^2x^2}{1-u(ux^2+\Xi )}\right)\]
\[ +u^2x^6\left( 1+ \frac{u(1+u)x^2}{(1+ux)(1-2ux)}\right)\left(1+\frac{u^2x}{1-ux}\right)
\frac{1}{1-u(\Sigma +\Xi )}\frac{1}{1-ux}\left(1+\frac{u^2x^2}{1-u(ux^2+\Xi)}\right) 
\]
\[+u^3x^6
\left(1+\frac{u(1+u)x^2}{(1+ux)(1-2ux)}\right) \frac{1}{1-u^2x^2}
\frac{1}{1-u(ux^2+\Xi )}
.\]

\[q'_7 (x, u)=  \frac{ux^5}{1-2u\Xi } \left(1+ux+
\frac{u^2x^2}{(1+ux)(1-2ux)}\left\{2(u^2-u+1)x+3 \right\}
\right)\left(1+\frac{ux(1+2x)}{(1+ux)(1-2ux)}\right)
\]
\[\times \left\{\left( 1+\frac{ux^2}{1-u^2x^2}\right) \frac{u^2x^2}{1-u^2x^2}
\ \frac{1}{1-u(\ell +\Xi )} 
+\left(1+\frac{u^3x^2}{1-u^2x^2}\right) \right\}
\]
\[+ux^5 \left(1+\frac{u^2(1+u)x^2}{(1+ux)(1-2ux)}\right) 
\left(1+\frac{ux^2}{1-u^2x^2} \right) \frac{1}{1-u(\ell +\Xi)}\ 
\frac{1}{1-u^2x^2}
.\]

\[q'_8 (x, u)= \left\{u^2x^4  \left(1+ux +
\frac{u^2x^2}{(1+ux)(1-2ux)}\left\{2(u^2-u+1)x+3)\right\} \right) \frac{1}{1-2u\Xi } 
\right.\]
\[
 \times \left( 1+ux+x^2\frac{(4u^2-2u)x+u^2+2}{(1+ux)(1-2ux)}\right) 
\]
\[\left.
+x^2 \left(1+\frac{u^2x^2}{(1+ux)(1-2ux)}\left\{(1-u)^2x+2\right\} \right)
\right\} \frac{1}{1-u(\ell +\Xi)}\ \frac{x^2}{1-u^2x^2}
.\]

\[q'_9 (x, u)=uxq'_8 (x, u).\]

\[q'_{10} (x, u)= 
2u(1+u)x^6\left(1+\frac{u^2x(1+2x)}{(1+ux)(1-2ux)} \right)^2 \frac{1}{1-2u\Xi }
\]
\[ \times \left\{ x \left(1+\frac{u^2x}{1-ux}\right)\frac{1}{1-u(\Sigma +\Xi )}
\frac{u^2x^2}{1-ux}+x+\frac{u^3x^3}{1-u^2x^2}
\right\} \ u\ 
  \frac{x^2}{1-u(ux^2+\Xi )}
\]
\[+ux^2 \left(1+ \frac{u(1+u)x^2}{(1+ux)(1-2ux)}\right)
\left\{ x\left( 1+\frac{u^2x}{1-ux} \right)\frac{1}{1-u(\Sigma +\Xi )}
\frac{u^2x^2}{1-ux}+x+\frac{u^3x^3}{1-u^2x^2}
\right\}\]
\[\times \frac{x^2}{1-u(ux^2+\Xi )}.\]

\vspace{0.5cm}\par
Using MATHEMATICA, we have confirmed that  as functions of $x$ and $u$, the 
following holds:
\[\sum _{i=1}^{10} p_i(x, u)=\Phi (x, u)/x, \ \ \sum _{i=1}^{10} q_i(x, u)=\Phi (x, u)^2/x^2, \]
as  required by the 
definitions of $\hat{P}^u_1$ and ${\hat{P}}'^u_1$,
where $\Phi (x,u) $ is defined in \eqnu{phi}.


\end{document}

%% file: Fig2a.tex
{\unitlength 0.1in%
\begin{picture}( 40.7300, 23.0000)(  3.0700,-26.4000)%
%
\special{pn 8}%
\special{pa 2760 2375}%
\special{pa 2268 2375}%
\special{pa 2513 1949}%
\special{pa 2760 2375}%
\special{pa 2268 2375}%
\special{fp}%
%
\special{pn 8}%
\special{pa 2513 2375}%
\special{pa 2637 2160}%
\special{pa 2760 2375}%
\special{pa 2513 2375}%
\special{pa 2637 2160}%
\special{fp}%
%
\special{pn 8}%
\special{pa 2513 2375}%
\special{pa 2637 2160}%
\special{pa 2760 2375}%
\special{pa 2513 2375}%
\special{pa 2637 2160}%
\special{fp}%
%
\special{pn 8}%
\special{pa 2513 2375}%
\special{pa 2637 2160}%
\special{pa 2760 2375}%
\special{pa 2513 2375}%
\special{pa 2637 2160}%
\special{fp}%
%
\special{pn 8}%
\special{pa 2513 2375}%
\special{pa 2637 2160}%
\special{pa 2760 2375}%
\special{pa 2513 2375}%
\special{pa 2637 2160}%
\special{fp}%
%
\special{pn 8}%
\special{pa 2268 2375}%
\special{pa 2392 2160}%
\special{pa 2514 2375}%
\special{pa 2268 2375}%
\special{pa 2392 2160}%
\special{fp}%
%
\special{pn 8}%
\special{pa 2268 2375}%
\special{pa 2392 2160}%
\special{pa 2514 2375}%
\special{pa 2268 2375}%
\special{pa 2392 2160}%
\special{fp}%
%
\special{pn 8}%
\special{pa 2268 2375}%
\special{pa 2392 2160}%
\special{pa 2514 2375}%
\special{pa 2268 2375}%
\special{pa 2392 2160}%
\special{fp}%
%
\special{pn 8}%
\special{pa 2268 2375}%
\special{pa 2392 2160}%
\special{pa 2514 2375}%
\special{pa 2268 2375}%
\special{pa 2392 2160}%
\special{fp}%
%
\special{pn 8}%
\special{pa 2392 2160}%
\special{pa 2513 1948}%
\special{pa 2637 2160}%
\special{pa 2392 2160}%
\special{pa 2513 1948}%
\special{fp}%
%
\special{pn 8}%
\special{pa 3250 2375}%
\special{pa 2759 2375}%
\special{pa 3004 1948}%
\special{pa 3250 2375}%
\special{pa 2759 2375}%
\special{fp}%
%
\special{pn 8}%
\special{pa 3004 2375}%
\special{pa 3128 2161}%
\special{pa 3250 2375}%
\special{pa 3004 2375}%
\special{pa 3128 2161}%
\special{fp}%
%
\special{pn 8}%
\special{pa 3004 2375}%
\special{pa 3128 2161}%
\special{pa 3250 2375}%
\special{pa 3004 2375}%
\special{pa 3128 2161}%
\special{fp}%
%
\special{pn 8}%
\special{pa 3004 2375}%
\special{pa 3128 2161}%
\special{pa 3250 2375}%
\special{pa 3004 2375}%
\special{pa 3128 2161}%
\special{fp}%
%
\special{pn 8}%
\special{pa 3004 2375}%
\special{pa 3128 2161}%
\special{pa 3250 2375}%
\special{pa 3004 2375}%
\special{pa 3128 2161}%
\special{fp}%
%
\special{pn 8}%
\special{pa 2759 2375}%
\special{pa 2882 2161}%
\special{pa 3005 2375}%
\special{pa 2759 2375}%
\special{pa 2882 2161}%
\special{fp}%
%
\special{pn 8}%
\special{pa 2759 2375}%
\special{pa 2882 2161}%
\special{pa 3005 2375}%
\special{pa 2759 2375}%
\special{pa 2882 2161}%
\special{fp}%
%
\special{pn 8}%
\special{pa 2759 2375}%
\special{pa 2882 2161}%
\special{pa 3005 2375}%
\special{pa 2759 2375}%
\special{pa 2882 2161}%
\special{fp}%
%
\special{pn 8}%
\special{pa 2759 2375}%
\special{pa 2882 2161}%
\special{pa 3005 2375}%
\special{pa 2759 2375}%
\special{pa 2882 2161}%
\special{fp}%
%
\special{pn 8}%
\special{pa 2882 2161}%
\special{pa 3004 1948}%
\special{pa 3128 2161}%
\special{pa 2882 2161}%
\special{pa 3004 1948}%
\special{fp}%
%
\special{pn 8}%
\special{pa 3004 1948}%
\special{pa 2513 1948}%
\special{pa 2759 1523}%
\special{pa 3004 1948}%
\special{pa 2513 1948}%
\special{fp}%
%
\special{pn 8}%
\special{pa 2513 1948}%
\special{pa 2636 1735}%
\special{pa 2759 1948}%
\special{pa 2513 1948}%
\special{pa 2636 1735}%
\special{fp}%
%
\special{pn 8}%
\special{pa 2759 1948}%
\special{pa 2882 1735}%
\special{pa 3004 1948}%
\special{pa 2759 1948}%
\special{pa 2882 1735}%
\special{fp}%
%
\special{pn 8}%
\special{pa 2635 1735}%
\special{pa 2759 1521}%
\special{pa 2882 1735}%
\special{pa 2635 1735}%
\special{pa 2759 1521}%
\special{fp}%
%
\special{pn 8}%
\special{pa 3741 2375}%
\special{pa 3250 2375}%
\special{pa 3494 1949}%
\special{pa 3741 2375}%
\special{pa 3250 2375}%
\special{fp}%
%
\special{pn 8}%
\special{pa 3494 2375}%
\special{pa 3618 2160}%
\special{pa 3741 2375}%
\special{pa 3494 2375}%
\special{pa 3618 2160}%
\special{fp}%
%
\special{pn 8}%
\special{pa 3494 2375}%
\special{pa 3618 2160}%
\special{pa 3741 2375}%
\special{pa 3494 2375}%
\special{pa 3618 2160}%
\special{fp}%
%
\special{pn 8}%
\special{pa 3494 2375}%
\special{pa 3618 2160}%
\special{pa 3741 2375}%
\special{pa 3494 2375}%
\special{pa 3618 2160}%
\special{fp}%
%
\special{pn 8}%
\special{pa 3494 2375}%
\special{pa 3618 2160}%
\special{pa 3741 2375}%
\special{pa 3494 2375}%
\special{pa 3618 2160}%
\special{fp}%
%
\special{pn 8}%
\special{pa 3250 2375}%
\special{pa 3372 2160}%
\special{pa 3495 2375}%
\special{pa 3250 2375}%
\special{pa 3372 2160}%
\special{fp}%
%
\special{pn 8}%
\special{pa 3250 2375}%
\special{pa 3372 2160}%
\special{pa 3495 2375}%
\special{pa 3250 2375}%
\special{pa 3372 2160}%
\special{fp}%
%
\special{pn 8}%
\special{pa 3250 2375}%
\special{pa 3372 2160}%
\special{pa 3495 2375}%
\special{pa 3250 2375}%
\special{pa 3372 2160}%
\special{fp}%
%
\special{pn 8}%
\special{pa 3250 2375}%
\special{pa 3372 2160}%
\special{pa 3495 2375}%
\special{pa 3250 2375}%
\special{pa 3372 2160}%
\special{fp}%
%
\special{pn 8}%
\special{pa 3372 2160}%
\special{pa 3494 1948}%
\special{pa 3618 2160}%
\special{pa 3372 2160}%
\special{pa 3494 1948}%
\special{fp}%
%
\special{pn 8}%
\special{pa 4232 2375}%
\special{pa 3740 2375}%
\special{pa 3985 1948}%
\special{pa 4232 2375}%
\special{pa 3740 2375}%
\special{fp}%
%
\special{pn 8}%
\special{pa 3985 2375}%
\special{pa 4108 2161}%
\special{pa 4232 2375}%
\special{pa 3985 2375}%
\special{pa 4108 2161}%
\special{fp}%
%
\special{pn 8}%
\special{pa 3985 2375}%
\special{pa 4108 2161}%
\special{pa 4232 2375}%
\special{pa 3985 2375}%
\special{pa 4108 2161}%
\special{fp}%
%
\special{pn 8}%
\special{pa 3985 2375}%
\special{pa 4108 2161}%
\special{pa 4232 2375}%
\special{pa 3985 2375}%
\special{pa 4108 2161}%
\special{fp}%
%
\special{pn 8}%
\special{pa 3985 2375}%
\special{pa 4108 2161}%
\special{pa 4232 2375}%
\special{pa 3985 2375}%
\special{pa 4108 2161}%
\special{fp}%
%
\special{pn 8}%
\special{pa 3740 2375}%
\special{pa 3863 2161}%
\special{pa 3986 2375}%
\special{pa 3740 2375}%
\special{pa 3863 2161}%
\special{fp}%
%
\special{pn 8}%
\special{pa 3740 2375}%
\special{pa 3863 2161}%
\special{pa 3986 2375}%
\special{pa 3740 2375}%
\special{pa 3863 2161}%
\special{fp}%
%
\special{pn 8}%
\special{pa 3740 2375}%
\special{pa 3863 2161}%
\special{pa 3986 2375}%
\special{pa 3740 2375}%
\special{pa 3863 2161}%
\special{fp}%
%
\special{pn 8}%
\special{pa 3740 2375}%
\special{pa 3863 2161}%
\special{pa 3986 2375}%
\special{pa 3740 2375}%
\special{pa 3863 2161}%
\special{fp}%
%
\special{pn 8}%
\special{pa 3863 2161}%
\special{pa 3985 1948}%
\special{pa 4108 2161}%
\special{pa 3863 2161}%
\special{pa 3985 1948}%
\special{fp}%
%
\special{pn 8}%
\special{pa 3985 1948}%
\special{pa 3494 1948}%
\special{pa 3740 1523}%
\special{pa 3985 1948}%
\special{pa 3494 1948}%
\special{fp}%
%
\special{pn 8}%
\special{pa 3494 1948}%
\special{pa 3617 1735}%
\special{pa 3740 1948}%
\special{pa 3494 1948}%
\special{pa 3617 1735}%
\special{fp}%
%
\special{pn 8}%
\special{pa 3740 1948}%
\special{pa 3863 1735}%
\special{pa 3985 1948}%
\special{pa 3740 1948}%
\special{pa 3863 1735}%
\special{fp}%
%
\special{pn 8}%
\special{pa 3616 1735}%
\special{pa 3740 1521}%
\special{pa 3863 1735}%
\special{pa 3616 1735}%
\special{pa 3740 1521}%
\special{fp}%
%
\special{pn 8}%
\special{pa 3249 1521}%
\special{pa 2757 1521}%
\special{pa 3002 1095}%
\special{pa 3249 1521}%
\special{pa 2757 1521}%
\special{fp}%
%
\special{pn 8}%
\special{pa 3002 1521}%
\special{pa 3125 1306}%
\special{pa 3249 1521}%
\special{pa 3002 1521}%
\special{pa 3125 1306}%
\special{fp}%
%
\special{pn 8}%
\special{pa 3002 1521}%
\special{pa 3125 1306}%
\special{pa 3249 1521}%
\special{pa 3002 1521}%
\special{pa 3125 1306}%
\special{fp}%
%
\special{pn 8}%
\special{pa 3002 1521}%
\special{pa 3125 1306}%
\special{pa 3249 1521}%
\special{pa 3002 1521}%
\special{pa 3125 1306}%
\special{fp}%
%
\special{pn 8}%
\special{pa 3002 1521}%
\special{pa 3125 1306}%
\special{pa 3249 1521}%
\special{pa 3002 1521}%
\special{pa 3125 1306}%
\special{fp}%
%
\special{pn 8}%
\special{pa 2757 1521}%
\special{pa 2879 1306}%
\special{pa 3003 1521}%
\special{pa 2757 1521}%
\special{pa 2879 1306}%
\special{fp}%
%
\special{pn 8}%
\special{pa 2757 1521}%
\special{pa 2879 1306}%
\special{pa 3003 1521}%
\special{pa 2757 1521}%
\special{pa 2879 1306}%
\special{fp}%
%
\special{pn 8}%
\special{pa 2757 1521}%
\special{pa 2879 1306}%
\special{pa 3003 1521}%
\special{pa 2757 1521}%
\special{pa 2879 1306}%
\special{fp}%
%
\special{pn 8}%
\special{pa 2757 1521}%
\special{pa 2879 1306}%
\special{pa 3003 1521}%
\special{pa 2757 1521}%
\special{pa 2879 1306}%
\special{fp}%
%
\special{pn 8}%
\special{pa 2879 1306}%
\special{pa 3002 1094}%
\special{pa 3125 1306}%
\special{pa 2879 1306}%
\special{pa 3002 1094}%
\special{fp}%
%
\special{pn 8}%
\special{pa 3739 1521}%
\special{pa 3248 1521}%
\special{pa 3492 1094}%
\special{pa 3739 1521}%
\special{pa 3248 1521}%
\special{fp}%
%
\special{pn 8}%
\special{pa 3492 1521}%
\special{pa 3616 1307}%
\special{pa 3739 1521}%
\special{pa 3492 1521}%
\special{pa 3616 1307}%
\special{fp}%
%
\special{pn 8}%
\special{pa 3492 1521}%
\special{pa 3616 1307}%
\special{pa 3739 1521}%
\special{pa 3492 1521}%
\special{pa 3616 1307}%
\special{fp}%
%
\special{pn 8}%
\special{pa 3492 1521}%
\special{pa 3616 1307}%
\special{pa 3739 1521}%
\special{pa 3492 1521}%
\special{pa 3616 1307}%
\special{fp}%
%
\special{pn 8}%
\special{pa 3492 1521}%
\special{pa 3616 1307}%
\special{pa 3739 1521}%
\special{pa 3492 1521}%
\special{pa 3616 1307}%
\special{fp}%
%
\special{pn 8}%
\special{pa 3248 1521}%
\special{pa 3370 1307}%
\special{pa 3493 1521}%
\special{pa 3248 1521}%
\special{pa 3370 1307}%
\special{fp}%
%
\special{pn 8}%
\special{pa 3248 1521}%
\special{pa 3370 1307}%
\special{pa 3493 1521}%
\special{pa 3248 1521}%
\special{pa 3370 1307}%
\special{fp}%
%
\special{pn 8}%
\special{pa 3248 1521}%
\special{pa 3370 1307}%
\special{pa 3493 1521}%
\special{pa 3248 1521}%
\special{pa 3370 1307}%
\special{fp}%
%
\special{pn 8}%
\special{pa 3248 1521}%
\special{pa 3370 1307}%
\special{pa 3493 1521}%
\special{pa 3248 1521}%
\special{pa 3370 1307}%
\special{fp}%
%
\special{pn 8}%
\special{pa 3370 1307}%
\special{pa 3492 1094}%
\special{pa 3616 1307}%
\special{pa 3370 1307}%
\special{pa 3492 1094}%
\special{fp}%
%
\special{pn 8}%
\special{pa 3002 1094}%
\special{pa 3124 881}%
\special{pa 3248 1094}%
\special{pa 3002 1094}%
\special{pa 3124 881}%
\special{fp}%
%
\special{pn 8}%
\special{pa 3248 1094}%
\special{pa 3370 881}%
\special{pa 3492 1094}%
\special{pa 3248 1094}%
\special{pa 3370 881}%
\special{fp}%
%
\special{pn 8}%
\special{pa 3124 880}%
\special{pa 3248 666}%
\special{pa 3370 880}%
\special{pa 3124 880}%
\special{pa 3248 666}%
\special{fp}%
%
\special{pn 8}%
\special{pa 798 2375}%
\special{pa 307 2375}%
\special{pa 552 1949}%
\special{pa 798 2375}%
\special{pa 307 2375}%
\special{fp}%
%
\special{pn 8}%
\special{pa 552 2375}%
\special{pa 676 2160}%
\special{pa 798 2375}%
\special{pa 552 2375}%
\special{pa 676 2160}%
\special{fp}%
%
\special{pn 8}%
\special{pa 552 2375}%
\special{pa 676 2160}%
\special{pa 798 2375}%
\special{pa 552 2375}%
\special{pa 676 2160}%
\special{fp}%
%
\special{pn 8}%
\special{pa 552 2375}%
\special{pa 676 2160}%
\special{pa 798 2375}%
\special{pa 552 2375}%
\special{pa 676 2160}%
\special{fp}%
%
\special{pn 8}%
\special{pa 552 2375}%
\special{pa 676 2160}%
\special{pa 798 2375}%
\special{pa 552 2375}%
\special{pa 676 2160}%
\special{fp}%
%
\special{pn 8}%
\special{pa 307 2375}%
\special{pa 430 2160}%
\special{pa 552 2375}%
\special{pa 307 2375}%
\special{pa 430 2160}%
\special{fp}%
%
\special{pn 8}%
\special{pa 307 2375}%
\special{pa 430 2160}%
\special{pa 552 2375}%
\special{pa 307 2375}%
\special{pa 430 2160}%
\special{fp}%
%
\special{pn 8}%
\special{pa 307 2375}%
\special{pa 430 2160}%
\special{pa 552 2375}%
\special{pa 307 2375}%
\special{pa 430 2160}%
\special{fp}%
%
\special{pn 8}%
\special{pa 307 2375}%
\special{pa 430 2160}%
\special{pa 552 2375}%
\special{pa 307 2375}%
\special{pa 430 2160}%
\special{fp}%
%
\special{pn 8}%
\special{pa 430 2160}%
\special{pa 552 1948}%
\special{pa 676 2160}%
\special{pa 430 2160}%
\special{pa 552 1948}%
\special{fp}%
%
\special{pn 8}%
\special{pa 1290 2375}%
\special{pa 798 2375}%
\special{pa 1043 1948}%
\special{pa 1290 2375}%
\special{pa 798 2375}%
\special{fp}%
%
\special{pn 8}%
\special{pa 1043 2375}%
\special{pa 1166 2161}%
\special{pa 1290 2375}%
\special{pa 1043 2375}%
\special{pa 1166 2161}%
\special{fp}%
%
\special{pn 8}%
\special{pa 1043 2375}%
\special{pa 1166 2161}%
\special{pa 1290 2375}%
\special{pa 1043 2375}%
\special{pa 1166 2161}%
\special{fp}%
%
\special{pn 8}%
\special{pa 1043 2375}%
\special{pa 1166 2161}%
\special{pa 1290 2375}%
\special{pa 1043 2375}%
\special{pa 1166 2161}%
\special{fp}%
%
\special{pn 8}%
\special{pa 1043 2375}%
\special{pa 1166 2161}%
\special{pa 1290 2375}%
\special{pa 1043 2375}%
\special{pa 1166 2161}%
\special{fp}%
%
\special{pn 8}%
\special{pa 798 2375}%
\special{pa 921 2161}%
\special{pa 1044 2375}%
\special{pa 798 2375}%
\special{pa 921 2161}%
\special{fp}%
%
\special{pn 8}%
\special{pa 798 2375}%
\special{pa 921 2161}%
\special{pa 1044 2375}%
\special{pa 798 2375}%
\special{pa 921 2161}%
\special{fp}%
%
\special{pn 8}%
\special{pa 798 2375}%
\special{pa 921 2161}%
\special{pa 1044 2375}%
\special{pa 798 2375}%
\special{pa 921 2161}%
\special{fp}%
%
\special{pn 8}%
\special{pa 798 2375}%
\special{pa 921 2161}%
\special{pa 1044 2375}%
\special{pa 798 2375}%
\special{pa 921 2161}%
\special{fp}%
%
\special{pn 8}%
\special{pa 921 2161}%
\special{pa 1043 1948}%
\special{pa 1166 2161}%
\special{pa 921 2161}%
\special{pa 1043 1948}%
\special{fp}%
%
\special{pn 8}%
\special{pa 1043 1948}%
\special{pa 552 1948}%
\special{pa 798 1523}%
\special{pa 1043 1948}%
\special{pa 552 1948}%
\special{fp}%
%
\special{pn 8}%
\special{pa 552 1948}%
\special{pa 675 1735}%
\special{pa 798 1948}%
\special{pa 552 1948}%
\special{pa 675 1735}%
\special{fp}%
%
\special{pn 8}%
\special{pa 798 1948}%
\special{pa 921 1735}%
\special{pa 1043 1948}%
\special{pa 798 1948}%
\special{pa 921 1735}%
\special{fp}%
%
\special{pn 8}%
\special{pa 674 1735}%
\special{pa 798 1521}%
\special{pa 921 1735}%
\special{pa 674 1735}%
\special{pa 798 1521}%
\special{fp}%
%
\special{pn 8}%
\special{pa 1779 2375}%
\special{pa 1288 2375}%
\special{pa 1533 1949}%
\special{pa 1779 2375}%
\special{pa 1288 2375}%
\special{fp}%
%
\special{pn 8}%
\special{pa 1533 2375}%
\special{pa 1656 2160}%
\special{pa 1779 2375}%
\special{pa 1533 2375}%
\special{pa 1656 2160}%
\special{fp}%
%
\special{pn 8}%
\special{pa 1533 2375}%
\special{pa 1656 2160}%
\special{pa 1779 2375}%
\special{pa 1533 2375}%
\special{pa 1656 2160}%
\special{fp}%
%
\special{pn 8}%
\special{pa 1533 2375}%
\special{pa 1656 2160}%
\special{pa 1779 2375}%
\special{pa 1533 2375}%
\special{pa 1656 2160}%
\special{fp}%
%
\special{pn 8}%
\special{pa 1533 2375}%
\special{pa 1656 2160}%
\special{pa 1779 2375}%
\special{pa 1533 2375}%
\special{pa 1656 2160}%
\special{fp}%
%
\special{pn 8}%
\special{pa 1288 2375}%
\special{pa 1410 2160}%
\special{pa 1534 2375}%
\special{pa 1288 2375}%
\special{pa 1410 2160}%
\special{fp}%
%
\special{pn 8}%
\special{pa 1288 2375}%
\special{pa 1410 2160}%
\special{pa 1534 2375}%
\special{pa 1288 2375}%
\special{pa 1410 2160}%
\special{fp}%
%
\special{pn 8}%
\special{pa 1288 2375}%
\special{pa 1410 2160}%
\special{pa 1534 2375}%
\special{pa 1288 2375}%
\special{pa 1410 2160}%
\special{fp}%
%
\special{pn 8}%
\special{pa 1288 2375}%
\special{pa 1410 2160}%
\special{pa 1534 2375}%
\special{pa 1288 2375}%
\special{pa 1410 2160}%
\special{fp}%
%
\special{pn 8}%
\special{pa 1410 2160}%
\special{pa 1533 1948}%
\special{pa 1656 2160}%
\special{pa 1410 2160}%
\special{pa 1533 1948}%
\special{fp}%
%
\special{pn 8}%
\special{pa 2270 2375}%
\special{pa 1778 2375}%
\special{pa 2023 1948}%
\special{pa 2270 2375}%
\special{pa 1778 2375}%
\special{fp}%
%
\special{pn 8}%
\special{pa 2023 2375}%
\special{pa 2147 2161}%
\special{pa 2270 2375}%
\special{pa 2023 2375}%
\special{pa 2147 2161}%
\special{fp}%
%
\special{pn 8}%
\special{pa 2023 2375}%
\special{pa 2147 2161}%
\special{pa 2270 2375}%
\special{pa 2023 2375}%
\special{pa 2147 2161}%
\special{fp}%
%
\special{pn 8}%
\special{pa 2023 2375}%
\special{pa 2147 2161}%
\special{pa 2270 2375}%
\special{pa 2023 2375}%
\special{pa 2147 2161}%
\special{fp}%
%
\special{pn 8}%
\special{pa 2023 2375}%
\special{pa 2147 2161}%
\special{pa 2270 2375}%
\special{pa 2023 2375}%
\special{pa 2147 2161}%
\special{fp}%
%
\special{pn 8}%
\special{pa 1778 2375}%
\special{pa 1901 2161}%
\special{pa 2024 2375}%
\special{pa 1778 2375}%
\special{pa 1901 2161}%
\special{fp}%
%
\special{pn 8}%
\special{pa 1778 2375}%
\special{pa 1901 2161}%
\special{pa 2024 2375}%
\special{pa 1778 2375}%
\special{pa 1901 2161}%
\special{fp}%
%
\special{pn 8}%
\special{pa 1778 2375}%
\special{pa 1901 2161}%
\special{pa 2024 2375}%
\special{pa 1778 2375}%
\special{pa 1901 2161}%
\special{fp}%
%
\special{pn 8}%
\special{pa 1778 2375}%
\special{pa 1901 2161}%
\special{pa 2024 2375}%
\special{pa 1778 2375}%
\special{pa 1901 2161}%
\special{fp}%
%
\special{pn 8}%
\special{pa 1901 2161}%
\special{pa 2023 1948}%
\special{pa 2147 2161}%
\special{pa 1901 2161}%
\special{pa 2023 1948}%
\special{fp}%
%
\special{pn 8}%
\special{pa 2023 1948}%
\special{pa 1533 1948}%
\special{pa 1778 1523}%
\special{pa 2023 1948}%
\special{pa 1533 1948}%
\special{fp}%
%
\special{pn 8}%
\special{pa 1533 1948}%
\special{pa 1655 1735}%
\special{pa 1778 1948}%
\special{pa 1533 1948}%
\special{pa 1655 1735}%
\special{fp}%
%
\special{pn 8}%
\special{pa 1778 1948}%
\special{pa 1901 1735}%
\special{pa 2023 1948}%
\special{pa 1778 1948}%
\special{pa 1901 1735}%
\special{fp}%
%
\special{pn 8}%
\special{pa 1655 1735}%
\special{pa 1778 1521}%
\special{pa 1901 1735}%
\special{pa 1655 1735}%
\special{pa 1778 1521}%
\special{fp}%
%
\special{pn 8}%
\special{pa 1287 1521}%
\special{pa 795 1521}%
\special{pa 1040 1095}%
\special{pa 1287 1521}%
\special{pa 795 1521}%
\special{fp}%
%
\special{pn 8}%
\special{pa 1040 1521}%
\special{pa 1164 1306}%
\special{pa 1287 1521}%
\special{pa 1040 1521}%
\special{pa 1164 1306}%
\special{fp}%
%
\special{pn 8}%
\special{pa 1040 1521}%
\special{pa 1164 1306}%
\special{pa 1287 1521}%
\special{pa 1040 1521}%
\special{pa 1164 1306}%
\special{fp}%
%
\special{pn 8}%
\special{pa 1040 1521}%
\special{pa 1164 1306}%
\special{pa 1287 1521}%
\special{pa 1040 1521}%
\special{pa 1164 1306}%
\special{fp}%
%
\special{pn 8}%
\special{pa 1040 1521}%
\special{pa 1164 1306}%
\special{pa 1287 1521}%
\special{pa 1040 1521}%
\special{pa 1164 1306}%
\special{fp}%
%
\special{pn 8}%
\special{pa 795 1521}%
\special{pa 919 1306}%
\special{pa 1041 1521}%
\special{pa 795 1521}%
\special{pa 919 1306}%
\special{fp}%
%
\special{pn 8}%
\special{pa 795 1521}%
\special{pa 919 1306}%
\special{pa 1041 1521}%
\special{pa 795 1521}%
\special{pa 919 1306}%
\special{fp}%
%
\special{pn 8}%
\special{pa 795 1521}%
\special{pa 919 1306}%
\special{pa 1041 1521}%
\special{pa 795 1521}%
\special{pa 919 1306}%
\special{fp}%
%
\special{pn 8}%
\special{pa 795 1521}%
\special{pa 919 1306}%
\special{pa 1041 1521}%
\special{pa 795 1521}%
\special{pa 919 1306}%
\special{fp}%
%
\special{pn 8}%
\special{pa 919 1306}%
\special{pa 1040 1094}%
\special{pa 1164 1306}%
\special{pa 919 1306}%
\special{pa 1040 1094}%
\special{fp}%
%
\special{pn 8}%
\special{pa 1778 1521}%
\special{pa 1286 1521}%
\special{pa 1531 1094}%
\special{pa 1778 1521}%
\special{pa 1286 1521}%
\special{fp}%
%
\special{pn 8}%
\special{pa 1531 1521}%
\special{pa 1655 1307}%
\special{pa 1778 1521}%
\special{pa 1531 1521}%
\special{pa 1655 1307}%
\special{fp}%
%
\special{pn 8}%
\special{pa 1531 1521}%
\special{pa 1655 1307}%
\special{pa 1778 1521}%
\special{pa 1531 1521}%
\special{pa 1655 1307}%
\special{fp}%
%
\special{pn 8}%
\special{pa 1531 1521}%
\special{pa 1655 1307}%
\special{pa 1778 1521}%
\special{pa 1531 1521}%
\special{pa 1655 1307}%
\special{fp}%
%
\special{pn 8}%
\special{pa 1531 1521}%
\special{pa 1655 1307}%
\special{pa 1778 1521}%
\special{pa 1531 1521}%
\special{pa 1655 1307}%
\special{fp}%
%
\special{pn 8}%
\special{pa 1286 1521}%
\special{pa 1409 1307}%
\special{pa 1532 1521}%
\special{pa 1286 1521}%
\special{pa 1409 1307}%
\special{fp}%
%
\special{pn 8}%
\special{pa 1286 1521}%
\special{pa 1409 1307}%
\special{pa 1532 1521}%
\special{pa 1286 1521}%
\special{pa 1409 1307}%
\special{fp}%
%
\special{pn 8}%
\special{pa 1286 1521}%
\special{pa 1409 1307}%
\special{pa 1532 1521}%
\special{pa 1286 1521}%
\special{pa 1409 1307}%
\special{fp}%
%
\special{pn 8}%
\special{pa 1286 1521}%
\special{pa 1409 1307}%
\special{pa 1532 1521}%
\special{pa 1286 1521}%
\special{pa 1409 1307}%
\special{fp}%
%
\special{pn 8}%
\special{pa 1409 1307}%
\special{pa 1531 1094}%
\special{pa 1655 1307}%
\special{pa 1409 1307}%
\special{pa 1531 1094}%
\special{fp}%
%
\special{pn 8}%
\special{pa 1531 1094}%
\special{pa 1040 1094}%
\special{pa 1286 669}%
\special{pa 1531 1094}%
\special{pa 1040 1094}%
\special{fp}%
%
\special{pn 8}%
\special{pa 1040 1094}%
\special{pa 1163 881}%
\special{pa 1286 1094}%
\special{pa 1040 1094}%
\special{pa 1163 881}%
\special{fp}%
%
\special{pn 8}%
\special{pa 1286 1094}%
\special{pa 1409 881}%
\special{pa 1531 1094}%
\special{pa 1286 1094}%
\special{pa 1409 881}%
\special{fp}%
%
\special{pn 8}%
\special{pa 1163 880}%
\special{pa 1286 666}%
\special{pa 1409 880}%
\special{pa 1163 880}%
\special{pa 1286 666}%
\special{fp}%
\put(22.1300,-24.2500){\makebox(0,0)[lt]{$O$}}%
\put(12.2000,-27.7000){\makebox(0,0)[lb]{$F_3$}}%
\put(33.0000,-4.7000){\makebox(0,0)[lb]{$a$}}%
\put(43.8000,-24.5000){\makebox(0,0)[lb]{$b$}}%
\end{picture}}%

%% file: decompositionsmall2.tex
{\unitlength 0.1in%
\begin{picture}( 41.2900, 24.1000)(  7.1900,-24.3000)%
%
\special{pn 8}%
\special{pa 790 2385}%
\special{pa 1021 1980}%
\special{pa 1257 2385}%
\special{pa 790 2385}%
\special{pa 1021 1980}%
\special{fp}%
%
\special{pn 8}%
\special{pa 790 2385}%
\special{pa 1021 1980}%
\special{pa 1257 2385}%
\special{pa 790 2385}%
\special{pa 1021 1980}%
\special{fp}%
%
\special{pn 8}%
\special{pa 790 2385}%
\special{pa 1021 1980}%
\special{pa 1257 2385}%
\special{pa 790 2385}%
\special{pa 1021 1980}%
\special{fp}%
%
\special{pn 8}%
\special{pa 790 2385}%
\special{pa 1021 1980}%
\special{pa 1257 2385}%
\special{pa 790 2385}%
\special{pa 1021 1980}%
\special{fp}%
%
\special{pn 8}%
\special{pa 1021 1980}%
\special{pa 1255 1577}%
\special{pa 1489 1980}%
\special{pa 1021 1980}%
\special{pa 1255 1577}%
\special{fp}%
\put(17.0900,-24.3000){\makebox(0,0)[lt]{$b$}}%
\put(11.6600,-14.5700){\makebox(0,0)[lb]{$a$}}%
\put(7.4400,-24.3000){\makebox(0,0)[lt]{$O$}}%
%
\special{pn 8}%
\special{pa 974 1064}%
\special{pa 1091 862}%
\special{pa 1208 1064}%
\special{pa 974 1064}%
\special{pa 1091 862}%
\special{fp}%
%
\special{pn 8}%
\special{pa 974 1064}%
\special{pa 1091 862}%
\special{pa 1208 1064}%
\special{pa 974 1064}%
\special{pa 1091 862}%
\special{fp}%
%
\special{pn 8}%
\special{pa 974 1064}%
\special{pa 1091 862}%
\special{pa 1208 1064}%
\special{pa 974 1064}%
\special{pa 1091 862}%
\special{fp}%
%
\special{pn 8}%
\special{pa 974 1064}%
\special{pa 1091 862}%
\special{pa 1208 1064}%
\special{pa 974 1064}%
\special{pa 1091 862}%
\special{fp}%
%
\special{pn 8}%
\special{pa 743 1064}%
\special{pa 858 862}%
\special{pa 975 1064}%
\special{pa 743 1064}%
\special{pa 858 862}%
\special{fp}%
%
\special{pn 8}%
\special{pa 743 1064}%
\special{pa 858 862}%
\special{pa 975 1064}%
\special{pa 743 1064}%
\special{pa 858 862}%
\special{fp}%
%
\special{pn 8}%
\special{pa 743 1064}%
\special{pa 858 862}%
\special{pa 975 1064}%
\special{pa 743 1064}%
\special{pa 858 862}%
\special{fp}%
%
\special{pn 8}%
\special{pa 743 1064}%
\special{pa 858 862}%
\special{pa 975 1064}%
\special{pa 743 1064}%
\special{pa 858 862}%
\special{fp}%
\put(7.1900,-10.8700){\makebox(0,0)[lt]{$O$}}%
%
\special{pn 8}%
\special{pa 1673 1063}%
\special{pa 1208 1063}%
\special{pa 1440 660}%
\special{pa 1673 1063}%
\special{pa 1208 1063}%
\special{fp}%
%
\special{pn 8}%
\special{pa 1440 1063}%
\special{pa 1557 861}%
\special{pa 1673 1063}%
\special{pa 1440 1063}%
\special{pa 1557 861}%
\special{fp}%
%
\special{pn 8}%
\special{pa 1440 1063}%
\special{pa 1557 861}%
\special{pa 1673 1063}%
\special{pa 1440 1063}%
\special{pa 1557 861}%
\special{fp}%
%
\special{pn 8}%
\special{pa 1440 1063}%
\special{pa 1557 861}%
\special{pa 1673 1063}%
\special{pa 1440 1063}%
\special{pa 1557 861}%
\special{fp}%
%
\special{pn 8}%
\special{pa 1440 1063}%
\special{pa 1557 861}%
\special{pa 1673 1063}%
\special{pa 1440 1063}%
\special{pa 1557 861}%
\special{fp}%
%
\special{pn 8}%
\special{pa 1208 1063}%
\special{pa 1324 861}%
\special{pa 1440 1063}%
\special{pa 1208 1063}%
\special{pa 1324 861}%
\special{fp}%
%
\special{pn 8}%
\special{pa 1208 1063}%
\special{pa 1324 861}%
\special{pa 1440 1063}%
\special{pa 1208 1063}%
\special{pa 1324 861}%
\special{fp}%
%
\special{pn 8}%
\special{pa 1208 1063}%
\special{pa 1324 861}%
\special{pa 1440 1063}%
\special{pa 1208 1063}%
\special{pa 1324 861}%
\special{fp}%
%
\special{pn 8}%
\special{pa 1208 1063}%
\special{pa 1324 861}%
\special{pa 1440 1063}%
\special{pa 1208 1063}%
\special{pa 1324 861}%
\special{fp}%
%
\special{pn 8}%
\special{pa 1324 861}%
\special{pa 1440 660}%
\special{pa 1557 861}%
\special{pa 1324 861}%
\special{pa 1440 660}%
\special{fp}%
\put(16.6700,-10.8600){\makebox(0,0)[lt]{$b$}}%
%
\special{pn 8}%
\special{pa 1440 660}%
\special{pa 974 660}%
\special{pa 1208 256}%
\special{pa 1440 660}%
\special{pa 974 660}%
\special{fp}%
%
\special{pn 8}%
\special{pa 974 660}%
\special{pa 1091 457}%
\special{pa 1208 660}%
\special{pa 974 660}%
\special{pa 1091 457}%
\special{fp}%
%
\special{pn 8}%
\special{pa 1207 660}%
\special{pa 1324 457}%
\special{pa 1440 660}%
\special{pa 1207 660}%
\special{pa 1324 457}%
\special{fp}%
%
\special{pn 8}%
\special{pa 1091 457}%
\special{pa 1208 254}%
\special{pa 1324 457}%
\special{pa 1091 457}%
\special{pa 1208 254}%
\special{fp}%
\put(11.6600,-2.0000){\makebox(0,0)[lb]{$a$}}%
%
\special{pn 20}%
\special{pa 3490 20}%
\special{pa 3490 20}%
\special{fp}%
%
\special{pn 20}%
\special{pa 1205 1062}%
\special{pa 1322 861}%
\special{fp}%
%
\special{pn 20}%
\special{pa 1322 861}%
\special{pa 1548 861}%
\special{fp}%
%
\special{pn 20}%
\special{pa 1548 861}%
\special{pa 1435 655}%
\special{fp}%
%
\special{pn 20}%
\special{pa 1435 655}%
\special{pa 1210 258}%
\special{fp}%
%
\special{pn 8}%
\special{pa 2960 1043}%
\special{pa 3193 639}%
\special{pa 3426 1043}%
\special{pa 2960 1043}%
\special{pa 3193 639}%
\special{fp}%
%
\special{pn 8}%
\special{pa 2960 1043}%
\special{pa 3193 639}%
\special{pa 3426 1043}%
\special{pa 2960 1043}%
\special{pa 3193 639}%
\special{fp}%
%
\special{pn 8}%
\special{pa 2960 1043}%
\special{pa 3193 639}%
\special{pa 3426 1043}%
\special{pa 2960 1043}%
\special{pa 3193 639}%
\special{fp}%
%
\special{pn 8}%
\special{pa 2495 1043}%
\special{pa 2726 639}%
\special{pa 2961 1043}%
\special{pa 2495 1043}%
\special{pa 2726 639}%
\special{fp}%
%
\special{pn 8}%
\special{pa 2495 1043}%
\special{pa 2726 639}%
\special{pa 2961 1043}%
\special{pa 2495 1043}%
\special{pa 2726 639}%
\special{fp}%
%
\special{pn 8}%
\special{pa 2495 1043}%
\special{pa 2726 639}%
\special{pa 2961 1043}%
\special{pa 2495 1043}%
\special{pa 2726 639}%
\special{fp}%
%
\special{pn 8}%
\special{pa 2495 1043}%
\special{pa 2726 639}%
\special{pa 2961 1043}%
\special{pa 2495 1043}%
\special{pa 2726 639}%
\special{fp}%
%
\special{pn 8}%
\special{pa 2726 639}%
\special{pa 2960 235}%
\special{pa 3193 639}%
\special{pa 2726 639}%
\special{pa 2960 235}%
\special{fp}%
\put(34.1300,-10.8800){\makebox(0,0)[lt]{$b$}}%
\put(29.1800,-1.7100){\makebox(0,0)[lb]{$a$}}%
\put(24.4900,-10.8800){\makebox(0,0)[lt]{$O$}}%
%
\special{pn 20}%
\special{pa 2478 1036}%
\special{pa 2723 634}%
\special{fp}%
%
\special{pn 20}%
\special{pa 2723 634}%
\special{pa 2948 1036}%
\special{fp}%
%
\special{pn 20}%
\special{pa 2948 1036}%
\special{pa 3184 634}%
\special{fp}%
%
\special{pn 20}%
\special{pa 3184 634}%
\special{pa 2958 223}%
\special{fp}%
%
\special{pn 20}%
\special{pa 857 861}%
\special{pa 1082 861}%
\special{fp}%
%
\special{pn 20}%
\special{pa 1092 454}%
\special{pa 1317 454}%
\special{fp}%
%
\special{pn 8}%
\special{pa 2813 2365}%
\special{pa 3047 1961}%
\special{pa 3280 2365}%
\special{pa 2813 2365}%
\special{pa 3047 1961}%
\special{fp}%
%
\special{pn 8}%
\special{pa 2813 2365}%
\special{pa 3047 1961}%
\special{pa 3280 2365}%
\special{pa 2813 2365}%
\special{pa 3047 1961}%
\special{fp}%
%
\special{pn 8}%
\special{pa 2813 2365}%
\special{pa 3047 1961}%
\special{pa 3280 2365}%
\special{pa 2813 2365}%
\special{pa 3047 1961}%
\special{fp}%
%
\special{pn 8}%
\special{pa 2349 2365}%
\special{pa 2580 1961}%
\special{pa 2815 2365}%
\special{pa 2349 2365}%
\special{pa 2580 1961}%
\special{fp}%
%
\special{pn 8}%
\special{pa 2349 2365}%
\special{pa 2580 1961}%
\special{pa 2815 2365}%
\special{pa 2349 2365}%
\special{pa 2580 1961}%
\special{fp}%
%
\special{pn 8}%
\special{pa 2349 2365}%
\special{pa 2580 1961}%
\special{pa 2815 2365}%
\special{pa 2349 2365}%
\special{pa 2580 1961}%
\special{fp}%
%
\special{pn 8}%
\special{pa 2349 2365}%
\special{pa 2580 1961}%
\special{pa 2815 2365}%
\special{pa 2349 2365}%
\special{pa 2580 1961}%
\special{fp}%
%
\special{pn 8}%
\special{pa 2580 1961}%
\special{pa 2813 1558}%
\special{pa 3047 1961}%
\special{pa 2580 1961}%
\special{pa 2813 1558}%
\special{fp}%
\put(32.6700,-24.1000){\makebox(0,0)[lt]{$b$}}%
\put(27.4300,-14.5700){\makebox(0,0)[lb]{$a$}}%
\put(23.0200,-24.1000){\makebox(0,0)[lt]{$O$}}%
%
\special{pn 8}%
\special{pa 4381 2365}%
\special{pa 4615 1961}%
\special{pa 4848 2365}%
\special{pa 4381 2365}%
\special{pa 4615 1961}%
\special{fp}%
%
\special{pn 8}%
\special{pa 4381 2365}%
\special{pa 4615 1961}%
\special{pa 4848 2365}%
\special{pa 4381 2365}%
\special{pa 4615 1961}%
\special{fp}%
%
\special{pn 8}%
\special{pa 4381 2365}%
\special{pa 4615 1961}%
\special{pa 4848 2365}%
\special{pa 4381 2365}%
\special{pa 4615 1961}%
\special{fp}%
%
\special{pn 8}%
\special{pa 3917 2365}%
\special{pa 4148 1961}%
\special{pa 4383 2365}%
\special{pa 3917 2365}%
\special{pa 4148 1961}%
\special{fp}%
%
\special{pn 8}%
\special{pa 3917 2365}%
\special{pa 4148 1961}%
\special{pa 4383 2365}%
\special{pa 3917 2365}%
\special{pa 4148 1961}%
\special{fp}%
%
\special{pn 8}%
\special{pa 3917 2365}%
\special{pa 4148 1961}%
\special{pa 4383 2365}%
\special{pa 3917 2365}%
\special{pa 4148 1961}%
\special{fp}%
%
\special{pn 8}%
\special{pa 3917 2365}%
\special{pa 4148 1961}%
\special{pa 4383 2365}%
\special{pa 3917 2365}%
\special{pa 4148 1961}%
\special{fp}%
%
\special{pn 8}%
\special{pa 4148 1961}%
\special{pa 4381 1558}%
\special{pa 4615 1961}%
\special{pa 4148 1961}%
\special{pa 4381 1558}%
\special{fp}%
\put(48.3500,-24.1000){\makebox(0,0)[lt]{$b$}}%
\put(43.2100,-14.5700){\makebox(0,0)[lb]{$a$}}%
\put(38.7000,-24.1000){\makebox(0,0)[lt]{$O$}}%
%
\special{pn 20}%
\special{pa 1254 2368}%
\special{pa 1489 1966}%
\special{fp}%
%
\special{pn 20}%
\special{pa 2351 2358}%
\special{pa 2577 1957}%
\special{fp}%
%
\special{pn 20}%
\special{pa 2577 1957}%
\special{pa 3037 1957}%
\special{fp}%
%
\special{pn 20}%
\special{pa 3037 1957}%
\special{pa 2802 1555}%
\special{fp}%
%
\special{pn 20}%
\special{pa 3910 2358}%
\special{pa 4370 1555}%
\special{fp}%
%
\special{pn 20}%
\special{pa 4145 1957}%
\special{pa 4596 1957}%
\special{fp}%
\put(14.5000,-4.2500){\makebox(0,0)[lb]{$w$}}%
\put(34.2900,-4.3800){\makebox(0,0)[lb]{$\tilde{w}$}}%
\put(16.0700,-17.4100){\makebox(0,0)[lb]{$w_1$}}%
\put(32.3300,-17.1200){\makebox(0,0)[lb]{$w_2$}}%
\put(47.8200,-17.1200){\makebox(0,0)[lb]{$w_3$}}%
%
\special{pn 8}%
\special{pa 1250 2385}%
\special{pa 1481 1980}%
\special{pa 1717 2385}%
\special{pa 1250 2385}%
\special{pa 1481 1980}%
\special{fp}%
%
\special{pn 8}%
\special{pa 1250 2385}%
\special{pa 1481 1980}%
\special{pa 1717 2385}%
\special{pa 1250 2385}%
\special{pa 1481 1980}%
\special{fp}%
%
\special{pn 8}%
\special{pa 1250 2385}%
\special{pa 1481 1980}%
\special{pa 1717 2385}%
\special{pa 1250 2385}%
\special{pa 1481 1980}%
\special{fp}%
%
\special{pn 8}%
\special{pa 1250 2385}%
\special{pa 1481 1980}%
\special{pa 1717 2385}%
\special{pa 1250 2385}%
\special{pa 1481 1980}%
\special{fp}%
%
\special{pn 20}%
\special{pa 790 2390}%
\special{pa 780 2380}%
\special{fp}%
\special{pa 1240 2390}%
\special{pa 1240 2390}%
\special{fp}%
%
\special{pn 20}%
\special{pa 790 2380}%
\special{pa 1230 2380}%
\special{fp}%
%
\special{pn 8}%
\special{pa 1090 850}%
\special{pa 970 670}%
\special{fp}%
%
\special{pn 8}%
\special{pa 970 660}%
\special{pa 840 860}%
\special{fp}%
%
\special{pn 20}%
\special{pa 1080 870}%
\special{pa 1200 1060}%
\special{fp}%
%
\special{pn 20}%
\special{pa 1250 1590}%
\special{pa 1470 1980}%
\special{fp}%
%
\special{pn 20}%
\special{pa 840 860}%
\special{pa 840 860}%
\special{fp}%
\special{pa 840 860}%
\special{pa 740 1070}%
\special{fp}%
\end{picture}}%

%% file: Fig4star.tex
{\unitlength 0.1in%
\begin{picture}( 42.0200, 19.0900)(  6.0000,-25.8900)%
%
\special{pn 4}%
\special{pa 1454 1472}%
\special{pa 728 1472}%
\special{pa 1091 844}%
\special{pa 1454 1472}%
\special{pa 728 1472}%
\special{fp}%
%
\special{pn 4}%
\special{pa 1091 1472}%
\special{pa 1272 1157}%
\special{pa 1454 1472}%
\special{pa 1091 1472}%
\special{pa 1272 1157}%
\special{fp}%
%
\special{pn 4}%
\special{pa 1091 1472}%
\special{pa 1272 1157}%
\special{pa 1454 1472}%
\special{pa 1091 1472}%
\special{pa 1272 1157}%
\special{fp}%
%
\special{pn 4}%
\special{pa 1091 1472}%
\special{pa 1272 1157}%
\special{pa 1454 1472}%
\special{pa 1091 1472}%
\special{pa 1272 1157}%
\special{fp}%
%
\special{pn 4}%
\special{pa 1091 1472}%
\special{pa 1272 1157}%
\special{pa 1454 1472}%
\special{pa 1091 1472}%
\special{pa 1272 1157}%
\special{fp}%
%
\special{pn 4}%
\special{pa 1091 1472}%
\special{pa 1272 1157}%
\special{pa 1454 1472}%
\special{pa 1091 1472}%
\special{pa 1272 1157}%
\special{fp}%
%
\special{pn 4}%
\special{pa 728 1472}%
\special{pa 910 1157}%
\special{pa 1091 1472}%
\special{pa 728 1472}%
\special{pa 910 1157}%
\special{fp}%
%
\special{pn 4}%
\special{pa 728 1472}%
\special{pa 910 1157}%
\special{pa 1091 1472}%
\special{pa 728 1472}%
\special{pa 910 1157}%
\special{fp}%
%
\special{pn 4}%
\special{pa 728 1472}%
\special{pa 910 1157}%
\special{pa 1091 1472}%
\special{pa 728 1472}%
\special{pa 910 1157}%
\special{fp}%
%
\special{pn 4}%
\special{pa 728 1472}%
\special{pa 910 1157}%
\special{pa 1091 1472}%
\special{pa 728 1472}%
\special{pa 910 1157}%
\special{fp}%
%
\special{pn 4}%
\special{pa 910 1157}%
\special{pa 1091 843}%
\special{pa 1272 1157}%
\special{pa 910 1157}%
\special{pa 1091 843}%
\special{fp}%
\put(10.9000,-8.1000){\makebox(0,0)[rb]{$a$}}%
\put(14.1000,-15.8000){\makebox(0,0){$b$}}%
%
\special{pn 4}%
\special{pa 2295 1471}%
\special{pa 1570 1471}%
\special{pa 1933 843}%
\special{pa 2295 1471}%
\special{pa 1570 1471}%
\special{fp}%
%
\special{pn 4}%
\special{pa 1933 1471}%
\special{pa 2115 1157}%
\special{pa 2295 1471}%
\special{pa 1933 1471}%
\special{pa 2115 1157}%
\special{fp}%
%
\special{pn 4}%
\special{pa 1933 1471}%
\special{pa 2115 1157}%
\special{pa 2295 1471}%
\special{pa 1933 1471}%
\special{pa 2115 1157}%
\special{fp}%
%
\special{pn 4}%
\special{pa 1933 1471}%
\special{pa 2115 1157}%
\special{pa 2295 1471}%
\special{pa 1933 1471}%
\special{pa 2115 1157}%
\special{fp}%
%
\special{pn 4}%
\special{pa 1933 1471}%
\special{pa 2115 1157}%
\special{pa 2295 1471}%
\special{pa 1933 1471}%
\special{pa 2115 1157}%
\special{fp}%
%
\special{pn 4}%
\special{pa 1933 1471}%
\special{pa 2115 1157}%
\special{pa 2295 1471}%
\special{pa 1933 1471}%
\special{pa 2115 1157}%
\special{fp}%
%
\special{pn 4}%
\special{pa 1570 1471}%
\special{pa 1751 1157}%
\special{pa 1933 1471}%
\special{pa 1570 1471}%
\special{pa 1751 1157}%
\special{fp}%
%
\special{pn 4}%
\special{pa 1570 1471}%
\special{pa 1751 1157}%
\special{pa 1933 1471}%
\special{pa 1570 1471}%
\special{pa 1751 1157}%
\special{fp}%
%
\special{pn 4}%
\special{pa 1570 1471}%
\special{pa 1751 1157}%
\special{pa 1933 1471}%
\special{pa 1570 1471}%
\special{pa 1751 1157}%
\special{fp}%
%
\special{pn 4}%
\special{pa 1570 1471}%
\special{pa 1751 1157}%
\special{pa 1933 1471}%
\special{pa 1570 1471}%
\special{pa 1751 1157}%
\special{fp}%
%
\special{pn 4}%
\special{pa 1751 1157}%
\special{pa 1933 843}%
\special{pa 2115 1157}%
\special{pa 1751 1157}%
\special{pa 1933 843}%
\special{fp}%
\put(6.0000,-16.3000){\makebox(0,0)[lb]{$O$}}%
%
\special{pn 20}%
\special{pa 728 1471}%
\special{pa 1089 844}%
\special{fp}%
%
\special{pn 4}%
\special{pa 3131 1471}%
\special{pa 2406 1471}%
\special{pa 2768 843}%
\special{pa 3131 1471}%
\special{pa 2406 1471}%
\special{fp}%
%
\special{pn 4}%
\special{pa 2768 1471}%
\special{pa 2949 1157}%
\special{pa 3131 1471}%
\special{pa 2768 1471}%
\special{pa 2949 1157}%
\special{fp}%
%
\special{pn 4}%
\special{pa 2768 1471}%
\special{pa 2949 1157}%
\special{pa 3131 1471}%
\special{pa 2768 1471}%
\special{pa 2949 1157}%
\special{fp}%
%
\special{pn 4}%
\special{pa 2768 1471}%
\special{pa 2949 1157}%
\special{pa 3131 1471}%
\special{pa 2768 1471}%
\special{pa 2949 1157}%
\special{fp}%
%
\special{pn 4}%
\special{pa 2768 1471}%
\special{pa 2949 1157}%
\special{pa 3131 1471}%
\special{pa 2768 1471}%
\special{pa 2949 1157}%
\special{fp}%
%
\special{pn 4}%
\special{pa 2768 1471}%
\special{pa 2949 1157}%
\special{pa 3131 1471}%
\special{pa 2768 1471}%
\special{pa 2949 1157}%
\special{fp}%
%
\special{pn 4}%
\special{pa 2406 1471}%
\special{pa 2587 1157}%
\special{pa 2768 1471}%
\special{pa 2406 1471}%
\special{pa 2587 1157}%
\special{fp}%
%
\special{pn 4}%
\special{pa 2406 1471}%
\special{pa 2587 1157}%
\special{pa 2768 1471}%
\special{pa 2406 1471}%
\special{pa 2587 1157}%
\special{fp}%
%
\special{pn 4}%
\special{pa 2406 1471}%
\special{pa 2587 1157}%
\special{pa 2768 1471}%
\special{pa 2406 1471}%
\special{pa 2587 1157}%
\special{fp}%
%
\special{pn 4}%
\special{pa 2406 1471}%
\special{pa 2587 1157}%
\special{pa 2768 1471}%
\special{pa 2406 1471}%
\special{pa 2587 1157}%
\special{fp}%
%
\special{pn 4}%
\special{pa 2587 1157}%
\special{pa 2768 843}%
\special{pa 2949 1157}%
\special{pa 2587 1157}%
\special{pa 2768 843}%
\special{fp}%
%
\special{pn 4}%
\special{pa 3098 2589}%
\special{pa 2373 2589}%
\special{pa 2735 1961}%
\special{pa 3098 2589}%
\special{pa 2373 2589}%
\special{fp}%
%
\special{pn 4}%
\special{pa 2735 2589}%
\special{pa 2916 2274}%
\special{pa 3098 2589}%
\special{pa 2735 2589}%
\special{pa 2916 2274}%
\special{fp}%
%
\special{pn 4}%
\special{pa 2735 2589}%
\special{pa 2916 2274}%
\special{pa 3098 2589}%
\special{pa 2735 2589}%
\special{pa 2916 2274}%
\special{fp}%
%
\special{pn 4}%
\special{pa 2735 2589}%
\special{pa 2916 2274}%
\special{pa 3098 2589}%
\special{pa 2735 2589}%
\special{pa 2916 2274}%
\special{fp}%
%
\special{pn 4}%
\special{pa 2735 2589}%
\special{pa 2916 2274}%
\special{pa 3098 2589}%
\special{pa 2735 2589}%
\special{pa 2916 2274}%
\special{fp}%
%
\special{pn 4}%
\special{pa 2735 2589}%
\special{pa 2916 2274}%
\special{pa 3098 2589}%
\special{pa 2735 2589}%
\special{pa 2916 2274}%
\special{fp}%
%
\special{pn 4}%
\special{pa 2373 2589}%
\special{pa 2554 2274}%
\special{pa 2735 2589}%
\special{pa 2373 2589}%
\special{pa 2554 2274}%
\special{fp}%
%
\special{pn 4}%
\special{pa 2373 2589}%
\special{pa 2554 2274}%
\special{pa 2735 2589}%
\special{pa 2373 2589}%
\special{pa 2554 2274}%
\special{fp}%
%
\special{pn 4}%
\special{pa 2373 2589}%
\special{pa 2554 2274}%
\special{pa 2735 2589}%
\special{pa 2373 2589}%
\special{pa 2554 2274}%
\special{fp}%
%
\special{pn 4}%
\special{pa 2373 2589}%
\special{pa 2554 2274}%
\special{pa 2735 2589}%
\special{pa 2373 2589}%
\special{pa 2554 2274}%
\special{fp}%
%
\special{pn 4}%
\special{pa 2554 2274}%
\special{pa 2735 1961}%
\special{pa 2916 2274}%
\special{pa 2554 2274}%
\special{pa 2735 1961}%
\special{fp}%
%
\special{pn 4}%
\special{pa 2262 2589}%
\special{pa 1537 2589}%
\special{pa 1900 1961}%
\special{pa 2262 2589}%
\special{pa 1537 2589}%
\special{fp}%
%
\special{pn 4}%
\special{pa 1900 2589}%
\special{pa 2082 2274}%
\special{pa 2262 2589}%
\special{pa 1900 2589}%
\special{pa 2082 2274}%
\special{fp}%
%
\special{pn 4}%
\special{pa 1900 2589}%
\special{pa 2082 2274}%
\special{pa 2262 2589}%
\special{pa 1900 2589}%
\special{pa 2082 2274}%
\special{fp}%
%
\special{pn 4}%
\special{pa 1900 2589}%
\special{pa 2082 2274}%
\special{pa 2262 2589}%
\special{pa 1900 2589}%
\special{pa 2082 2274}%
\special{fp}%
%
\special{pn 4}%
\special{pa 1900 2589}%
\special{pa 2082 2274}%
\special{pa 2262 2589}%
\special{pa 1900 2589}%
\special{pa 2082 2274}%
\special{fp}%
%
\special{pn 4}%
\special{pa 1900 2589}%
\special{pa 2082 2274}%
\special{pa 2262 2589}%
\special{pa 1900 2589}%
\special{pa 2082 2274}%
\special{fp}%
%
\special{pn 4}%
\special{pa 1537 2589}%
\special{pa 1718 2274}%
\special{pa 1900 2589}%
\special{pa 1537 2589}%
\special{pa 1718 2274}%
\special{fp}%
%
\special{pn 4}%
\special{pa 1537 2589}%
\special{pa 1718 2274}%
\special{pa 1900 2589}%
\special{pa 1537 2589}%
\special{pa 1718 2274}%
\special{fp}%
%
\special{pn 4}%
\special{pa 1537 2589}%
\special{pa 1718 2274}%
\special{pa 1900 2589}%
\special{pa 1537 2589}%
\special{pa 1718 2274}%
\special{fp}%
%
\special{pn 4}%
\special{pa 1537 2589}%
\special{pa 1718 2274}%
\special{pa 1900 2589}%
\special{pa 1537 2589}%
\special{pa 1718 2274}%
\special{fp}%
%
\special{pn 4}%
\special{pa 1718 2274}%
\special{pa 1900 1961}%
\special{pa 2082 2274}%
\special{pa 1718 2274}%
\special{pa 1900 1961}%
\special{fp}%
%
\special{pn 4}%
\special{pa 1427 2589}%
\special{pa 702 2589}%
\special{pa 1064 1961}%
\special{pa 1427 2589}%
\special{pa 702 2589}%
\special{fp}%
%
\special{pn 4}%
\special{pa 1064 2589}%
\special{pa 1246 2274}%
\special{pa 1428 2589}%
\special{pa 1064 2589}%
\special{pa 1246 2274}%
\special{fp}%
%
\special{pn 4}%
\special{pa 1064 2589}%
\special{pa 1246 2274}%
\special{pa 1428 2589}%
\special{pa 1064 2589}%
\special{pa 1246 2274}%
\special{fp}%
%
\special{pn 4}%
\special{pa 1064 2589}%
\special{pa 1246 2274}%
\special{pa 1428 2589}%
\special{pa 1064 2589}%
\special{pa 1246 2274}%
\special{fp}%
%
\special{pn 4}%
\special{pa 1064 2589}%
\special{pa 1246 2274}%
\special{pa 1428 2589}%
\special{pa 1064 2589}%
\special{pa 1246 2274}%
\special{fp}%
%
\special{pn 4}%
\special{pa 1064 2589}%
\special{pa 1246 2274}%
\special{pa 1428 2589}%
\special{pa 1064 2589}%
\special{pa 1246 2274}%
\special{fp}%
%
\special{pn 4}%
\special{pa 702 2589}%
\special{pa 883 2274}%
\special{pa 1065 2589}%
\special{pa 702 2589}%
\special{pa 883 2274}%
\special{fp}%
%
\special{pn 4}%
\special{pa 702 2589}%
\special{pa 883 2274}%
\special{pa 1065 2589}%
\special{pa 702 2589}%
\special{pa 883 2274}%
\special{fp}%
%
\special{pn 4}%
\special{pa 702 2589}%
\special{pa 883 2274}%
\special{pa 1065 2589}%
\special{pa 702 2589}%
\special{pa 883 2274}%
\special{fp}%
%
\special{pn 4}%
\special{pa 702 2589}%
\special{pa 883 2274}%
\special{pa 1065 2589}%
\special{pa 702 2589}%
\special{pa 883 2274}%
\special{fp}%
%
\special{pn 4}%
\special{pa 883 2274}%
\special{pa 1064 1961}%
\special{pa 1246 2274}%
\special{pa 883 2274}%
\special{pa 1064 1961}%
\special{fp}%
%
\special{pn 4}%
\special{pa 3966 1471}%
\special{pa 3241 1471}%
\special{pa 3604 843}%
\special{pa 3966 1471}%
\special{pa 3241 1471}%
\special{fp}%
%
\special{pn 4}%
\special{pa 3604 1471}%
\special{pa 3785 1157}%
\special{pa 3966 1471}%
\special{pa 3604 1471}%
\special{pa 3785 1157}%
\special{fp}%
%
\special{pn 4}%
\special{pa 3604 1471}%
\special{pa 3785 1157}%
\special{pa 3966 1471}%
\special{pa 3604 1471}%
\special{pa 3785 1157}%
\special{fp}%
%
\special{pn 4}%
\special{pa 3604 1471}%
\special{pa 3785 1157}%
\special{pa 3966 1471}%
\special{pa 3604 1471}%
\special{pa 3785 1157}%
\special{fp}%
%
\special{pn 4}%
\special{pa 3604 1471}%
\special{pa 3785 1157}%
\special{pa 3966 1471}%
\special{pa 3604 1471}%
\special{pa 3785 1157}%
\special{fp}%
%
\special{pn 4}%
\special{pa 3604 1471}%
\special{pa 3785 1157}%
\special{pa 3966 1471}%
\special{pa 3604 1471}%
\special{pa 3785 1157}%
\special{fp}%
%
\special{pn 4}%
\special{pa 3241 1471}%
\special{pa 3421 1157}%
\special{pa 3604 1471}%
\special{pa 3241 1471}%
\special{pa 3421 1157}%
\special{fp}%
%
\special{pn 4}%
\special{pa 3241 1471}%
\special{pa 3421 1157}%
\special{pa 3604 1471}%
\special{pa 3241 1471}%
\special{pa 3421 1157}%
\special{fp}%
%
\special{pn 4}%
\special{pa 3241 1471}%
\special{pa 3421 1157}%
\special{pa 3604 1471}%
\special{pa 3241 1471}%
\special{pa 3421 1157}%
\special{fp}%
%
\special{pn 4}%
\special{pa 3241 1471}%
\special{pa 3421 1157}%
\special{pa 3604 1471}%
\special{pa 3241 1471}%
\special{pa 3421 1157}%
\special{fp}%
%
\special{pn 4}%
\special{pa 3421 1157}%
\special{pa 3604 843}%
\special{pa 3785 1157}%
\special{pa 3421 1157}%
\special{pa 3604 843}%
\special{fp}%
%
\special{pn 4}%
\special{pa 4801 1471}%
\special{pa 4075 1471}%
\special{pa 4438 843}%
\special{pa 4801 1471}%
\special{pa 4075 1471}%
\special{fp}%
%
\special{pn 4}%
\special{pa 4438 1471}%
\special{pa 4620 1157}%
\special{pa 4802 1471}%
\special{pa 4438 1471}%
\special{pa 4620 1157}%
\special{fp}%
%
\special{pn 4}%
\special{pa 4438 1471}%
\special{pa 4620 1157}%
\special{pa 4802 1471}%
\special{pa 4438 1471}%
\special{pa 4620 1157}%
\special{fp}%
%
\special{pn 4}%
\special{pa 4438 1471}%
\special{pa 4620 1157}%
\special{pa 4802 1471}%
\special{pa 4438 1471}%
\special{pa 4620 1157}%
\special{fp}%
%
\special{pn 4}%
\special{pa 4438 1471}%
\special{pa 4620 1157}%
\special{pa 4802 1471}%
\special{pa 4438 1471}%
\special{pa 4620 1157}%
\special{fp}%
%
\special{pn 4}%
\special{pa 4438 1471}%
\special{pa 4620 1157}%
\special{pa 4802 1471}%
\special{pa 4438 1471}%
\special{pa 4620 1157}%
\special{fp}%
%
\special{pn 4}%
\special{pa 4075 1471}%
\special{pa 4257 1157}%
\special{pa 4439 1471}%
\special{pa 4075 1471}%
\special{pa 4257 1157}%
\special{fp}%
%
\special{pn 4}%
\special{pa 4075 1471}%
\special{pa 4257 1157}%
\special{pa 4439 1471}%
\special{pa 4075 1471}%
\special{pa 4257 1157}%
\special{fp}%
%
\special{pn 4}%
\special{pa 4075 1471}%
\special{pa 4257 1157}%
\special{pa 4439 1471}%
\special{pa 4075 1471}%
\special{pa 4257 1157}%
\special{fp}%
%
\special{pn 4}%
\special{pa 4075 1471}%
\special{pa 4257 1157}%
\special{pa 4439 1471}%
\special{pa 4075 1471}%
\special{pa 4257 1157}%
\special{fp}%
%
\special{pn 4}%
\special{pa 4257 1157}%
\special{pa 4438 843}%
\special{pa 4620 1157}%
\special{pa 4257 1157}%
\special{pa 4438 843}%
\special{fp}%
%
\special{pn 4}%
\special{pa 4768 2589}%
\special{pa 4042 2589}%
\special{pa 4405 1961}%
\special{pa 4768 2589}%
\special{pa 4042 2589}%
\special{fp}%
%
\special{pn 4}%
\special{pa 4405 2589}%
\special{pa 4587 2274}%
\special{pa 4769 2589}%
\special{pa 4405 2589}%
\special{pa 4587 2274}%
\special{fp}%
%
\special{pn 4}%
\special{pa 4405 2589}%
\special{pa 4587 2274}%
\special{pa 4769 2589}%
\special{pa 4405 2589}%
\special{pa 4587 2274}%
\special{fp}%
%
\special{pn 4}%
\special{pa 4405 2589}%
\special{pa 4587 2274}%
\special{pa 4769 2589}%
\special{pa 4405 2589}%
\special{pa 4587 2274}%
\special{fp}%
%
\special{pn 4}%
\special{pa 4405 2589}%
\special{pa 4587 2274}%
\special{pa 4769 2589}%
\special{pa 4405 2589}%
\special{pa 4587 2274}%
\special{fp}%
%
\special{pn 4}%
\special{pa 4405 2589}%
\special{pa 4587 2274}%
\special{pa 4769 2589}%
\special{pa 4405 2589}%
\special{pa 4587 2274}%
\special{fp}%
%
\special{pn 4}%
\special{pa 4042 2589}%
\special{pa 4224 2274}%
\special{pa 4406 2589}%
\special{pa 4042 2589}%
\special{pa 4224 2274}%
\special{fp}%
%
\special{pn 4}%
\special{pa 4042 2589}%
\special{pa 4224 2274}%
\special{pa 4406 2589}%
\special{pa 4042 2589}%
\special{pa 4224 2274}%
\special{fp}%
%
\special{pn 4}%
\special{pa 4042 2589}%
\special{pa 4224 2274}%
\special{pa 4406 2589}%
\special{pa 4042 2589}%
\special{pa 4224 2274}%
\special{fp}%
%
\special{pn 4}%
\special{pa 4042 2589}%
\special{pa 4224 2274}%
\special{pa 4406 2589}%
\special{pa 4042 2589}%
\special{pa 4224 2274}%
\special{fp}%
%
\special{pn 4}%
\special{pa 4224 2274}%
\special{pa 4405 1961}%
\special{pa 4587 2274}%
\special{pa 4224 2274}%
\special{pa 4405 1961}%
\special{fp}%
%
\special{pn 4}%
\special{pa 3933 2589}%
\special{pa 3208 2589}%
\special{pa 3571 1961}%
\special{pa 3933 2589}%
\special{pa 3208 2589}%
\special{fp}%
%
\special{pn 4}%
\special{pa 3571 2589}%
\special{pa 3752 2274}%
\special{pa 3933 2589}%
\special{pa 3571 2589}%
\special{pa 3752 2274}%
\special{fp}%
%
\special{pn 4}%
\special{pa 3571 2589}%
\special{pa 3752 2274}%
\special{pa 3933 2589}%
\special{pa 3571 2589}%
\special{pa 3752 2274}%
\special{fp}%
%
\special{pn 4}%
\special{pa 3571 2589}%
\special{pa 3752 2274}%
\special{pa 3933 2589}%
\special{pa 3571 2589}%
\special{pa 3752 2274}%
\special{fp}%
%
\special{pn 4}%
\special{pa 3571 2589}%
\special{pa 3752 2274}%
\special{pa 3933 2589}%
\special{pa 3571 2589}%
\special{pa 3752 2274}%
\special{fp}%
%
\special{pn 4}%
\special{pa 3571 2589}%
\special{pa 3752 2274}%
\special{pa 3933 2589}%
\special{pa 3571 2589}%
\special{pa 3752 2274}%
\special{fp}%
%
\special{pn 4}%
\special{pa 3208 2589}%
\special{pa 3388 2274}%
\special{pa 3571 2589}%
\special{pa 3208 2589}%
\special{pa 3388 2274}%
\special{fp}%
%
\special{pn 4}%
\special{pa 3208 2589}%
\special{pa 3388 2274}%
\special{pa 3571 2589}%
\special{pa 3208 2589}%
\special{pa 3388 2274}%
\special{fp}%
%
\special{pn 4}%
\special{pa 3208 2589}%
\special{pa 3388 2274}%
\special{pa 3571 2589}%
\special{pa 3208 2589}%
\special{pa 3388 2274}%
\special{fp}%
%
\special{pn 4}%
\special{pa 3208 2589}%
\special{pa 3388 2274}%
\special{pa 3571 2589}%
\special{pa 3208 2589}%
\special{pa 3388 2274}%
\special{fp}%
%
\special{pn 4}%
\special{pa 3388 2274}%
\special{pa 3571 1961}%
\special{pa 3752 2274}%
\special{pa 3388 2274}%
\special{pa 3571 1961}%
\special{fp}%
%
\special{pn 20}%
\special{pa 1570 1471}%
\special{pa 1933 1471}%
\special{fp}%
\special{pa 1933 1471}%
\special{pa 1750 1156}%
\special{fp}%
\special{pa 1750 1156}%
\special{pa 1931 840}%
\special{fp}%
%
\special{pn 20}%
\special{pa 2405 1471}%
\special{pa 2584 1156}%
\special{fp}%
\special{pa 2584 1156}%
\special{pa 2948 1156}%
\special{fp}%
\special{pa 2948 1156}%
\special{pa 2766 842}%
\special{fp}%
%
\special{pn 20}%
\special{pa 3238 1471}%
\special{pa 3602 1471}%
\special{fp}%
\special{pa 3602 1471}%
\special{pa 3420 1156}%
\special{fp}%
\special{pa 3420 1156}%
\special{pa 3784 1156}%
\special{fp}%
\special{pa 3784 1156}%
\special{pa 3602 842}%
\special{fp}%
%
\special{pn 20}%
\special{pa 4075 1471}%
\special{pa 4437 1471}%
\special{fp}%
\special{pa 4437 1471}%
\special{pa 4618 1156}%
\special{fp}%
\special{pa 4618 1156}%
\special{pa 4255 1156}%
\special{fp}%
\special{pa 4255 1156}%
\special{pa 4437 842}%
\special{fp}%
%
\special{pn 20}%
\special{pa 700 2589}%
\special{pa 881 2274}%
\special{fp}%
\special{pa 881 2274}%
\special{pa 1065 2589}%
\special{fp}%
\special{pa 1065 2589}%
\special{pa 1242 2274}%
\special{fp}%
\special{pa 1242 2274}%
\special{pa 1062 1960}%
\special{fp}%
%
\special{pn 20}%
\special{pa 1534 2589}%
\special{pa 1898 2589}%
\special{fp}%
\special{pa 1898 2589}%
\special{pa 2079 2274}%
\special{fp}%
\special{pa 2079 2274}%
\special{pa 1898 1960}%
\special{fp}%
%
\special{pn 20}%
\special{pa 2370 2589}%
\special{pa 3097 2589}%
\special{fp}%
\special{pa 3097 2589}%
\special{pa 2733 1960}%
\special{fp}%
%
\special{pn 20}%
\special{pa 3205 2589}%
\special{pa 3932 2589}%
\special{fp}%
\special{pa 3932 2589}%
\special{pa 3751 2274}%
\special{fp}%
\special{pa 3751 2274}%
\special{pa 3387 2274}%
\special{fp}%
\special{pa 3387 2274}%
\special{pa 3569 1960}%
\special{fp}%
%
\special{pn 20}%
\special{pa 4041 2589}%
\special{pa 4222 2274}%
\special{fp}%
\special{pa 4222 2274}%
\special{pa 4406 2589}%
\special{fp}%
\special{pa 4406 2589}%
\special{pa 4768 2589}%
\special{fp}%
\special{pa 4768 2589}%
\special{pa 4404 1960}%
\special{fp}%
\put(12.7100,-9.4700){\makebox(0,0)[lb]{$w^*_1$}}%
\put(21.0600,-9.4700){\makebox(0,0)[lb]{$w^*_2$}}%
\put(29.4200,-9.4700){\makebox(0,0)[lb]{$w^*_3$}}%
\put(37.7800,-9.4700){\makebox(0,0)[lb]{$w^*_4$}}%
\put(46.1200,-9.4700){\makebox(0,0)[lb]{$w^*_5$}}%
\put(12.3800,-21.0700){\makebox(0,0)[lb]{$w^*_6$}}%
\put(20.7300,-21.0700){\makebox(0,0)[lb]{$w^*_7$}}%
\put(29.0900,-21.0700){\makebox(0,0)[lb]{$w^*_8$}}%
\put(37.4500,-21.0700){\makebox(0,0)[lb]{$w^*_9$}}%
\put(45.7900,-21.0700){\makebox(0,0)[lb]{$w^*_{10}$}}%
\end{picture}}%

%% file: Fig3a.tex
\unitlength 0.1in
\begin{picture}( 41.1000, 38.0700)(  7.4000,-42.0700)
%
{\color[named]{Black}{%
\special{pn 4}%
\special{pa 2562 1904}%
\special{pa 890 1904}%
\special{pa 1726 518}%
\special{pa 2562 1904}%
\special{pa 890 1904}%
\special{fp}%
}}%
%
{\color[named]{Black}{%
\special{pn 4}%
\special{pa 1726 1904}%
\special{pa 2144 1212}%
\special{pa 2562 1904}%
\special{pa 1726 1904}%
\special{pa 2144 1212}%
\special{fp}%
}}%
%
{\color[named]{Black}{%
\special{pn 4}%
\special{pa 1726 1904}%
\special{pa 2144 1212}%
\special{pa 2562 1904}%
\special{pa 1726 1904}%
\special{pa 2144 1212}%
\special{fp}%
}}%
%
{\color[named]{Black}{%
\special{pn 4}%
\special{pa 1726 1904}%
\special{pa 2144 1212}%
\special{pa 2562 1904}%
\special{pa 1726 1904}%
\special{pa 2144 1212}%
\special{fp}%
}}%
%
{\color[named]{Black}{%
\special{pn 4}%
\special{pa 1726 1904}%
\special{pa 2144 1212}%
\special{pa 2562 1904}%
\special{pa 1726 1904}%
\special{pa 2144 1212}%
\special{fp}%
}}%
%
{\color[named]{Black}{%
\special{pn 4}%
\special{pa 1726 1904}%
\special{pa 2144 1212}%
\special{pa 2562 1904}%
\special{pa 1726 1904}%
\special{pa 2144 1212}%
\special{fp}%
}}%
%
{\color[named]{Black}{%
\special{pn 4}%
\special{pa 890 1904}%
\special{pa 1308 1212}%
\special{pa 1726 1904}%
\special{pa 890 1904}%
\special{pa 1308 1212}%
\special{fp}%
}}%
%
{\color[named]{Black}{%
\special{pn 4}%
\special{pa 890 1904}%
\special{pa 1308 1212}%
\special{pa 1726 1904}%
\special{pa 890 1904}%
\special{pa 1308 1212}%
\special{fp}%
}}%
%
{\color[named]{Black}{%
\special{pn 4}%
\special{pa 890 1904}%
\special{pa 1308 1212}%
\special{pa 1726 1904}%
\special{pa 890 1904}%
\special{pa 1308 1212}%
\special{fp}%
}}%
%
{\color[named]{Black}{%
\special{pn 4}%
\special{pa 890 1904}%
\special{pa 1308 1212}%
\special{pa 1726 1904}%
\special{pa 890 1904}%
\special{pa 1308 1212}%
\special{fp}%
}}%
%
{\color[named]{Black}{%
\special{pn 4}%
\special{pa 1308 1212}%
\special{pa 1726 518}%
\special{pa 2144 1212}%
\special{pa 1308 1212}%
\special{pa 1726 518}%
\special{fp}%
}}%
\put(17.9000,-5.3000){\makebox(0,0)[lb]{$a$}}%
%
{\color[named]{Black}{%
\special{pn 4}%
\special{pa 4702 1902}%
\special{pa 3030 1902}%
\special{pa 3866 516}%
\special{pa 4702 1902}%
\special{pa 3030 1902}%
\special{fp}%
}}%
%
{\color[named]{Black}{%
\special{pn 4}%
\special{pa 3866 1902}%
\special{pa 4284 1208}%
\special{pa 4702 1902}%
\special{pa 3866 1902}%
\special{pa 4284 1208}%
\special{fp}%
}}%
%
{\color[named]{Black}{%
\special{pn 4}%
\special{pa 3866 1902}%
\special{pa 4284 1208}%
\special{pa 4702 1902}%
\special{pa 3866 1902}%
\special{pa 4284 1208}%
\special{fp}%
}}%
%
{\color[named]{Black}{%
\special{pn 4}%
\special{pa 3866 1902}%
\special{pa 4284 1208}%
\special{pa 4702 1902}%
\special{pa 3866 1902}%
\special{pa 4284 1208}%
\special{fp}%
}}%
%
{\color[named]{Black}{%
\special{pn 4}%
\special{pa 3866 1902}%
\special{pa 4284 1208}%
\special{pa 4702 1902}%
\special{pa 3866 1902}%
\special{pa 4284 1208}%
\special{fp}%
}}%
%
{\color[named]{Black}{%
\special{pn 4}%
\special{pa 3866 1902}%
\special{pa 4284 1208}%
\special{pa 4702 1902}%
\special{pa 3866 1902}%
\special{pa 4284 1208}%
\special{fp}%
}}%
%
{\color[named]{Black}{%
\special{pn 4}%
\special{pa 3030 1902}%
\special{pa 3448 1208}%
\special{pa 3868 1902}%
\special{pa 3030 1902}%
\special{pa 3448 1208}%
\special{fp}%
}}%
%
{\color[named]{Black}{%
\special{pn 4}%
\special{pa 3030 1902}%
\special{pa 3448 1208}%
\special{pa 3868 1902}%
\special{pa 3030 1902}%
\special{pa 3448 1208}%
\special{fp}%
}}%
%
{\color[named]{Black}{%
\special{pn 4}%
\special{pa 3030 1902}%
\special{pa 3448 1208}%
\special{pa 3868 1902}%
\special{pa 3030 1902}%
\special{pa 3448 1208}%
\special{fp}%
}}%
%
{\color[named]{Black}{%
\special{pn 4}%
\special{pa 3030 1902}%
\special{pa 3448 1208}%
\special{pa 3868 1902}%
\special{pa 3030 1902}%
\special{pa 3448 1208}%
\special{fp}%
}}%
%
{\color[named]{Black}{%
\special{pn 4}%
\special{pa 3448 1208}%
\special{pa 3866 514}%
\special{pa 4284 1208}%
\special{pa 3448 1208}%
\special{pa 3866 514}%
\special{fp}%
}}%
%
{\color[named]{Black}{%
\special{pn 8}%
\special{pa 2558 1902}%
\special{pa 2536 1880}%
\special{pa 2512 1858}%
\special{pa 2486 1838}%
\special{pa 2458 1826}%
\special{pa 2434 1848}%
\special{pa 2414 1870}%
\special{pa 2392 1856}%
\special{pa 2364 1826}%
\special{pa 2332 1802}%
\special{pa 2308 1808}%
\special{pa 2296 1838}%
\special{pa 2280 1868}%
\special{pa 2258 1866}%
\special{pa 2236 1838}%
\special{pa 2218 1800}%
\special{pa 2202 1772}%
\special{pa 2168 1768}%
\special{pa 2146 1792}%
\special{pa 2128 1792}%
\special{pa 2108 1756}%
\special{pa 2074 1740}%
\special{pa 2044 1752}%
\special{pa 2034 1782}%
\special{pa 2042 1814}%
\special{pa 2020 1826}%
\special{pa 1978 1818}%
\special{pa 1948 1796}%
\special{pa 1920 1806}%
\special{pa 1908 1838}%
\special{pa 1906 1856}%
\special{pa 1870 1838}%
\special{pa 1846 1818}%
\special{pa 1812 1834}%
\special{pa 1806 1866}%
\special{pa 1778 1872}%
\special{pa 1750 1886}%
\special{pa 1724 1904}%
\special{pa 1720 1908}%
\special{fp}%
}}%
\put(7.9200,-19.4900){\makebox(0,0)[lt]{$O$}}%
%
{\color[named]{Black}{%
\special{pn 4}%
\special{pa 754 1902}%
\special{pa 910 1902}%
\special{fp}%
}}%
%
{\color[named]{Black}{%
\special{pn 4}%
\special{pa 2906 1902}%
\special{pa 3062 1902}%
\special{fp}%
}}%
%
{\color[named]{Black}{%
\special{pn 4}%
\special{pa 886 1902}%
\special{pa 772 1714}%
\special{fp}%
\special{pa 788 1738}%
\special{pa 798 1738}%
\special{fp}%
}}%
%
{\color[named]{Black}{%
\special{pn 4}%
\special{pa 3028 1902}%
\special{pa 2916 1714}%
\special{fp}%
\special{pa 2932 1738}%
\special{pa 2940 1738}%
\special{fp}%
}}%
%
{\color[named]{Black}{%
\special{pn 8}%
\special{pa 2560 1902}%
\special{pa 2542 1874}%
\special{pa 2526 1848}%
\special{pa 2508 1822}%
\special{pa 2486 1798}%
\special{pa 2462 1776}%
\special{pa 2436 1758}%
\special{pa 2404 1742}%
\special{pa 2368 1730}%
\special{pa 2332 1722}%
\special{pa 2304 1710}%
\special{pa 2298 1688}%
\special{pa 2310 1656}%
\special{pa 2324 1622}%
\special{pa 2328 1588}%
\special{pa 2322 1554}%
\special{pa 2310 1524}%
\special{pa 2288 1496}%
\special{pa 2262 1474}%
\special{pa 2230 1454}%
\special{pa 2200 1438}%
\special{pa 2180 1416}%
\special{pa 2182 1386}%
\special{pa 2186 1352}%
\special{pa 2178 1322}%
\special{pa 2160 1294}%
\special{pa 2148 1264}%
\special{pa 2144 1232}%
\special{pa 2144 1208}%
\special{fp}%
}}%
%
{\color[named]{Black}{%
\special{pn 8}%
\special{pa 2138 1208}%
\special{pa 2118 1184}%
\special{pa 2098 1162}%
\special{pa 2072 1140}%
\special{pa 2042 1120}%
\special{pa 2010 1100}%
\special{pa 1990 1080}%
\special{pa 1998 1058}%
\special{pa 2026 1036}%
\special{pa 2032 1032}%
\special{fp}%
}}%
%
{\color[named]{Black}{%
\special{pn 8}%
\special{pa 2138 1208}%
\special{pa 2126 1238}%
\special{pa 2116 1268}%
\special{pa 2110 1298}%
\special{pa 2104 1330}%
\special{pa 2104 1364}%
\special{pa 2110 1398}%
\special{pa 2104 1428}%
\special{pa 2074 1442}%
\special{pa 2046 1458}%
\special{pa 2042 1490}%
\special{pa 2048 1524}%
\special{pa 2040 1550}%
\special{pa 2008 1562}%
\special{pa 1978 1576}%
\special{pa 1970 1602}%
\special{pa 1976 1638}%
\special{pa 1972 1668}%
\special{pa 1950 1686}%
\special{pa 1916 1696}%
\special{pa 1884 1710}%
\special{pa 1860 1734}%
\special{pa 1844 1760}%
\special{pa 1828 1788}%
\special{pa 1808 1814}%
\special{pa 1784 1836}%
\special{pa 1762 1860}%
\special{pa 1742 1884}%
\special{pa 1724 1902}%
\special{fp}%
}}%
%
{\color[named]{Black}{%
\special{pn 4}%
\special{pa 740 4026}%
\special{pa 896 4026}%
\special{fp}%
}}%
%
{\color[named]{Black}{%
\special{pn 4}%
\special{pa 884 4024}%
\special{pa 770 3834}%
\special{fp}%
\special{pa 788 3862}%
\special{pa 796 3862}%
\special{fp}%
}}%
\put(16.3200,-21.1900){\makebox(0,0)[lt]{(a)}}%
%
{\color[named]{Black}{%
\special{pn 20}%
\special{pa 3028 1902}%
\special{pa 4702 1902}%
\special{fp}%
\special{pa 4702 1902}%
\special{pa 4280 1208}%
\special{fp}%
\special{pa 4280 1208}%
\special{pa 3864 1902}%
\special{fp}%
\special{pa 3864 1902}%
\special{pa 3446 1208}%
\special{fp}%
\special{pa 3446 1208}%
\special{pa 3864 518}%
\special{fp}%
}}%
\put(25.9000,-19.4000){\makebox(0,0)[lt]{$b$}}%
%
{\color[named]{Black}{%
\special{pn 8}%
\special{pa 880 1908}%
\special{pa 898 1884}%
\special{pa 922 1862}%
\special{pa 950 1846}%
\special{pa 982 1842}%
\special{pa 1012 1844}%
\special{pa 1046 1842}%
\special{pa 1068 1858}%
\special{pa 1098 1862}%
\special{pa 1136 1872}%
\special{pa 1166 1870}%
\special{pa 1168 1838}%
\special{pa 1182 1810}%
\special{pa 1206 1786}%
\special{pa 1232 1764}%
\special{pa 1258 1746}%
\special{pa 1290 1734}%
\special{pa 1310 1752}%
\special{pa 1310 1784}%
\special{pa 1326 1810}%
\special{pa 1358 1828}%
\special{pa 1392 1834}%
\special{pa 1424 1834}%
\special{pa 1456 1842}%
\special{pa 1488 1840}%
\special{pa 1518 1832}%
\special{pa 1550 1826}%
\special{pa 1576 1844}%
\special{pa 1594 1880}%
\special{pa 1608 1906}%
\special{pa 1628 1896}%
\special{pa 1652 1868}%
\special{pa 1684 1868}%
\special{pa 1706 1888}%
\special{pa 1706 1892}%
\special{fp}%
}}%
%
{\color[named]{Black}{%
\special{pn 8}%
\special{pa 1310 1228}%
\special{pa 1310 1258}%
\special{pa 1308 1288}%
\special{pa 1298 1328}%
\special{pa 1292 1364}%
\special{pa 1314 1376}%
\special{pa 1340 1396}%
\special{pa 1352 1426}%
\special{pa 1352 1490}%
\special{pa 1346 1524}%
\special{pa 1352 1554}%
\special{pa 1370 1580}%
\special{pa 1394 1592}%
\special{pa 1422 1604}%
\special{pa 1438 1628}%
\special{pa 1448 1664}%
\special{pa 1470 1686}%
\special{pa 1498 1706}%
\special{pa 1528 1718}%
\special{pa 1558 1722}%
\special{pa 1586 1742}%
\special{pa 1602 1774}%
\special{pa 1606 1810}%
\special{pa 1638 1808}%
\special{pa 1654 1832}%
\special{pa 1684 1842}%
\special{pa 1680 1862}%
\special{pa 1656 1858}%
\special{fp}%
}}%
%
{\color[named]{Black}{%
\special{pn 8}%
\special{pa 1302 1194}%
\special{pa 1334 1186}%
\special{pa 1362 1172}%
\special{pa 1390 1156}%
\special{pa 1418 1142}%
\special{pa 1448 1130}%
\special{pa 1480 1122}%
\special{pa 1512 1112}%
\special{pa 1538 1096}%
\special{pa 1562 1074}%
\special{pa 1586 1050}%
\special{pa 1606 1024}%
\special{pa 1614 996}%
\special{pa 1608 966}%
\special{pa 1600 934}%
\special{pa 1598 902}%
\special{pa 1598 870}%
\special{pa 1606 838}%
\special{pa 1624 810}%
\special{pa 1654 796}%
\special{pa 1690 794}%
\special{pa 1710 776}%
\special{pa 1724 744}%
\special{pa 1728 716}%
\special{pa 1718 686}%
\special{pa 1710 654}%
\special{pa 1712 626}%
\special{pa 1736 596}%
\special{pa 1730 574}%
\special{pa 1724 546}%
\special{pa 1724 530}%
\special{fp}%
}}%
%
{\color[named]{Black}{%
\special{pn 8}%
\special{pa 1706 784}%
\special{pa 1706 816}%
\special{pa 1708 846}%
\special{pa 1708 910}%
\special{pa 1740 920}%
\special{pa 1764 942}%
\special{pa 1800 946}%
\special{pa 1826 930}%
\special{pa 1824 896}%
\special{pa 1834 864}%
\special{pa 1826 838}%
\special{pa 1798 816}%
\special{pa 1776 794}%
\special{pa 1748 784}%
\special{pa 1718 776}%
\special{pa 1698 784}%
\special{fp}%
}}%
%
{\color[named]{Black}{%
\special{pn 4}%
\special{pa 2558 4020}%
\special{pa 896 4020}%
\special{pa 1726 2634}%
\special{pa 2558 4020}%
\special{pa 896 4020}%
\special{fp}%
}}%
%
{\color[named]{Black}{%
\special{pn 4}%
\special{pa 1726 4020}%
\special{pa 2144 3326}%
\special{pa 2560 4020}%
\special{pa 1726 4020}%
\special{pa 2144 3326}%
\special{fp}%
}}%
%
{\color[named]{Black}{%
\special{pn 4}%
\special{pa 1726 4020}%
\special{pa 2144 3326}%
\special{pa 2560 4020}%
\special{pa 1726 4020}%
\special{pa 2144 3326}%
\special{fp}%
}}%
%
{\color[named]{Black}{%
\special{pn 4}%
\special{pa 1726 4020}%
\special{pa 2144 3326}%
\special{pa 2560 4020}%
\special{pa 1726 4020}%
\special{pa 2144 3326}%
\special{fp}%
}}%
%
{\color[named]{Black}{%
\special{pn 4}%
\special{pa 1726 4020}%
\special{pa 2144 3326}%
\special{pa 2560 4020}%
\special{pa 1726 4020}%
\special{pa 2144 3326}%
\special{fp}%
}}%
%
{\color[named]{Black}{%
\special{pn 4}%
\special{pa 1726 4020}%
\special{pa 2144 3326}%
\special{pa 2560 4020}%
\special{pa 1726 4020}%
\special{pa 2144 3326}%
\special{fp}%
}}%
%
{\color[named]{Black}{%
\special{pn 4}%
\special{pa 896 4020}%
\special{pa 1312 3326}%
\special{pa 1728 4020}%
\special{pa 896 4020}%
\special{pa 1312 3326}%
\special{fp}%
}}%
%
{\color[named]{Black}{%
\special{pn 4}%
\special{pa 896 4020}%
\special{pa 1312 3326}%
\special{pa 1728 4020}%
\special{pa 896 4020}%
\special{pa 1312 3326}%
\special{fp}%
}}%
%
{\color[named]{Black}{%
\special{pn 4}%
\special{pa 896 4020}%
\special{pa 1312 3326}%
\special{pa 1728 4020}%
\special{pa 896 4020}%
\special{pa 1312 3326}%
\special{fp}%
}}%
%
{\color[named]{Black}{%
\special{pn 4}%
\special{pa 896 4020}%
\special{pa 1312 3326}%
\special{pa 1728 4020}%
\special{pa 896 4020}%
\special{pa 1312 3326}%
\special{fp}%
}}%
%
{\color[named]{Black}{%
\special{pn 4}%
\special{pa 1312 3326}%
\special{pa 1726 2634}%
\special{pa 2144 3326}%
\special{pa 1312 3326}%
\special{pa 1726 2634}%
\special{fp}%
}}%
%
{\color[named]{Black}{%
\special{pn 20}%
\special{pa 892 4018}%
\special{pa 1726 4018}%
\special{fp}%
\special{pa 1726 4018}%
\special{pa 1312 3326}%
\special{fp}%
\special{pa 1312 3326}%
\special{pa 1726 2634}%
\special{fp}%
}}%
\put(37.8800,-22.1600){\makebox(0,0)[lb]{(b)}}%
\put(31.2200,-20.8500){\makebox(0,0)[lb]{$O$}}%
\put(39.5000,-5.3000){\makebox(0,0)[lb]{$a$}}%
\put(47.4000,-20.6000){\makebox(0,0)[lb]{$b$}}%
\put(16.6000,-43.3700){\makebox(0,0)[lb]{(c)}}%
\put(8.1000,-42.3700){\makebox(0,0)[lb]{$O$}}%
\put(18.0000,-26.3000){\makebox(0,0)[lb]{$a$}}%
\put(26.1000,-41.5000){\makebox(0,0)[lb]{$b$}}%
%
{\color[named]{Black}{%
\special{pn 8}%
\special{pa 1350 1560}%
\special{pa 1324 1540}%
\special{pa 1294 1530}%
\special{pa 1264 1532}%
\special{pa 1236 1536}%
\special{pa 1230 1570}%
\special{pa 1220 1600}%
\special{pa 1222 1630}%
\special{pa 1256 1630}%
\special{pa 1276 1642}%
\special{pa 1308 1640}%
\special{pa 1340 1636}%
\special{pa 1358 1608}%
\special{pa 1360 1580}%
\special{fp}%
}}%
%
{\color[named]{Black}{%
\special{pn 4}%
\special{pa 4798 4048}%
\special{pa 3126 4048}%
\special{pa 3962 2662}%
\special{pa 4798 4048}%
\special{pa 3126 4048}%
\special{fp}%
}}%
%
{\color[named]{Black}{%
\special{pn 4}%
\special{pa 3962 4048}%
\special{pa 4380 3354}%
\special{pa 4798 4048}%
\special{pa 3962 4048}%
\special{pa 4380 3354}%
\special{fp}%
}}%
%
{\color[named]{Black}{%
\special{pn 4}%
\special{pa 3962 4048}%
\special{pa 4380 3354}%
\special{pa 4798 4048}%
\special{pa 3962 4048}%
\special{pa 4380 3354}%
\special{fp}%
}}%
%
{\color[named]{Black}{%
\special{pn 4}%
\special{pa 3962 4048}%
\special{pa 4380 3354}%
\special{pa 4798 4048}%
\special{pa 3962 4048}%
\special{pa 4380 3354}%
\special{fp}%
}}%
%
{\color[named]{Black}{%
\special{pn 4}%
\special{pa 3962 4048}%
\special{pa 4380 3354}%
\special{pa 4798 4048}%
\special{pa 3962 4048}%
\special{pa 4380 3354}%
\special{fp}%
}}%
%
{\color[named]{Black}{%
\special{pn 4}%
\special{pa 3962 4048}%
\special{pa 4380 3354}%
\special{pa 4798 4048}%
\special{pa 3962 4048}%
\special{pa 4380 3354}%
\special{fp}%
}}%
%
{\color[named]{Black}{%
\special{pn 4}%
\special{pa 3126 4048}%
\special{pa 3544 3354}%
\special{pa 3962 4048}%
\special{pa 3126 4048}%
\special{pa 3544 3354}%
\special{fp}%
}}%
%
{\color[named]{Black}{%
\special{pn 4}%
\special{pa 3126 4048}%
\special{pa 3544 3354}%
\special{pa 3962 4048}%
\special{pa 3126 4048}%
\special{pa 3544 3354}%
\special{fp}%
}}%
%
{\color[named]{Black}{%
\special{pn 4}%
\special{pa 3126 4048}%
\special{pa 3544 3354}%
\special{pa 3962 4048}%
\special{pa 3126 4048}%
\special{pa 3544 3354}%
\special{fp}%
}}%
%
{\color[named]{Black}{%
\special{pn 4}%
\special{pa 3126 4048}%
\special{pa 3544 3354}%
\special{pa 3962 4048}%
\special{pa 3126 4048}%
\special{pa 3544 3354}%
\special{fp}%
}}%
%
{\color[named]{Black}{%
\special{pn 4}%
\special{pa 3544 3354}%
\special{pa 3962 2660}%
\special{pa 4380 3354}%
\special{pa 3544 3354}%
\special{pa 3962 2660}%
\special{fp}%
}}%
\put(40.2000,-26.6000){\makebox(0,0)[lb]{$a$}}%
\put(30.2800,-40.9200){\makebox(0,0)[lt]{$O$}}%
%
{\color[named]{Black}{%
\special{pn 4}%
\special{pa 2990 4046}%
\special{pa 3146 4046}%
\special{fp}%
}}%
%
{\color[named]{Black}{%
\special{pn 4}%
\special{pa 3122 4046}%
\special{pa 3008 3856}%
\special{fp}%
\special{pa 3024 3882}%
\special{pa 3034 3882}%
\special{fp}%
}}%
\put(39.2000,-41.9700){\makebox(0,0)[lt]{(d)}}%
\put(48.5000,-41.2000){\makebox(0,0)[lt]{$b$}}%
%
{\color[named]{Black}{%
\special{pn 8}%
\special{pa 3116 4050}%
\special{pa 3134 4026}%
\special{pa 3158 4004}%
\special{pa 3186 3990}%
\special{pa 3218 3986}%
\special{pa 3248 3988}%
\special{pa 3282 3984}%
\special{pa 3304 4002}%
\special{pa 3334 4004}%
\special{pa 3372 4014}%
\special{pa 3402 4012}%
\special{pa 3404 3980}%
\special{pa 3418 3952}%
\special{pa 3442 3928}%
\special{pa 3468 3908}%
\special{pa 3494 3890}%
\special{pa 3526 3878}%
\special{pa 3546 3894}%
\special{pa 3546 3928}%
\special{pa 3562 3952}%
\special{pa 3594 3970}%
\special{pa 3628 3978}%
\special{pa 3660 3978}%
\special{pa 3692 3984}%
\special{pa 3724 3984}%
\special{pa 3754 3974}%
\special{pa 3786 3968}%
\special{pa 3812 3986}%
\special{pa 3830 4024}%
\special{pa 3844 4050}%
\special{pa 3864 4038}%
\special{pa 3888 4012}%
\special{pa 3920 4010}%
\special{pa 3942 4030}%
\special{pa 3942 4034}%
\special{fp}%
}}%
%
{\color[named]{Black}{%
\special{pn 8}%
\special{pa 3546 3370}%
\special{pa 3546 3402}%
\special{pa 3544 3432}%
\special{pa 3534 3470}%
\special{pa 3528 3508}%
\special{pa 3550 3518}%
\special{pa 3576 3540}%
\special{pa 3588 3570}%
\special{pa 3588 3634}%
\special{pa 3582 3668}%
\special{pa 3588 3698}%
\special{pa 3606 3724}%
\special{pa 3630 3736}%
\special{pa 3658 3748}%
\special{pa 3674 3772}%
\special{pa 3684 3808}%
\special{pa 3706 3830}%
\special{pa 3734 3848}%
\special{pa 3764 3862}%
\special{pa 3794 3864}%
\special{pa 3822 3884}%
\special{pa 3838 3916}%
\special{pa 3842 3952}%
\special{pa 3874 3950}%
\special{pa 3890 3974}%
\special{pa 3920 3984}%
\special{pa 3916 4006}%
\special{pa 3892 4002}%
\special{fp}%
}}%
%
{\color[named]{Black}{%
\special{pn 8}%
\special{pa 3538 3338}%
\special{pa 3570 3330}%
\special{pa 3598 3314}%
\special{pa 3626 3300}%
\special{pa 3654 3284}%
\special{pa 3684 3272}%
\special{pa 3748 3256}%
\special{pa 3774 3240}%
\special{pa 3822 3192}%
\special{pa 3842 3168}%
\special{pa 3850 3140}%
\special{pa 3844 3108}%
\special{pa 3836 3076}%
\special{pa 3834 3044}%
\special{pa 3834 3012}%
\special{pa 3842 2980}%
\special{pa 3860 2954}%
\special{pa 3890 2940}%
\special{pa 3926 2936}%
\special{pa 3946 2918}%
\special{pa 3960 2888}%
\special{pa 3964 2860}%
\special{pa 3954 2830}%
\special{pa 3946 2798}%
\special{pa 3948 2770}%
\special{pa 3972 2738}%
\special{pa 3966 2718}%
\special{pa 3960 2690}%
\special{pa 3960 2672}%
\special{fp}%
}}%
%
{\color[named]{Black}{%
\special{pn 8}%
\special{pa 3942 2926}%
\special{pa 3942 2958}%
\special{pa 3944 2990}%
\special{pa 3944 3052}%
\special{pa 3976 3064}%
\special{pa 4000 3086}%
\special{pa 4036 3090}%
\special{pa 4062 3072}%
\special{pa 4060 3040}%
\special{pa 4070 3006}%
\special{pa 4062 2980}%
\special{pa 4034 2960}%
\special{pa 4012 2938}%
\special{pa 3984 2928}%
\special{pa 3954 2918}%
\special{pa 3934 2926}%
\special{fp}%
}}%
%
{\color[named]{Black}{%
\special{pn 8}%
\special{pa 3586 3704}%
\special{pa 3560 3682}%
\special{pa 3530 3674}%
\special{pa 3500 3674}%
\special{pa 3472 3678}%
\special{pa 3466 3714}%
\special{pa 3456 3744}%
\special{pa 3458 3774}%
\special{pa 3492 3772}%
\special{pa 3512 3784}%
\special{pa 3544 3784}%
\special{pa 3576 3780}%
\special{pa 3594 3752}%
\special{pa 3596 3724}%
\special{fp}%
}}%
%
{\color[named]{Black}{%
\special{pn 20}%
\special{pa 3110 4048}%
\special{pa 3950 4048}%
\special{fp}%
}}%
%
{\color[named]{Black}{%
\special{pn 8}%
\special{pa 3940 4038}%
\special{pa 3950 4028}%
\special{fp}%
}}%
%
{\color[named]{Black}{%
\special{pn 20}%
\special{pa 3940 4048}%
\special{pa 3930 4048}%
\special{fp}%
\special{pa 3550 3328}%
\special{pa 3550 3328}%
\special{fp}%
}}%
%
{\color[named]{Black}{%
\special{pn 20}%
\special{pa 3950 4058}%
\special{pa 3950 4058}%
\special{fp}%
\special{pa 3940 4048}%
\special{pa 3530 3348}%
\special{fp}%
}}%
%
{\color[named]{Black}{%
\special{pn 20}%
\special{pa 3530 3348}%
\special{pa 3530 3348}%
\special{fp}%
\special{pa 3540 3348}%
\special{pa 3960 2658}%
\special{fp}%
}}%
%
\put(47.4000,-20.6000){\makebox(0,0)[lb]{}}%
\end{picture}%